%% file: regpol.tex
\overfullrule0pt
\def\firstpageno{1}
\pageno=\firstpageno
\magnification=\magstep1
\input psfig.sty
\input imsmark.tex
\SBIMSMark{1998/3}{April 1998}{}
\nopagenumbers

\headline={\ifnum\pageno>\firstpageno{
        \ifodd\pageno\oddheadline \else\evenheadline\fi}
        \else\hfil\fi}
\def\evenheadline{\sevenrm\folio\hfil ERIC BEDFORD AND MATTIAS JONSSON
    \hfil
}
\def\oddheadline{\sevenrm
    \hfil REGULAR POLYNOMIAL ENDOMORPHISMS\hfil\folio}
\outer\def\give#1. {\medbreak
             \noindent{\sl#1. }}                     
\outer\def\section #1\par{\bigbreak\centerline{\S
     {\bf#1}}\nobreak\smallskip\noindent}
\def\contract{{\,{\vrule height5pt width0.4pt depth 0pt}
    {\vrule height0.4pt width6pt depth 0pt}\,}}  
\def\sqr#1#2{{\vcenter{\hrule height.#2pt              
      \hbox{\vrule width.#2pt height#1pt\kern#1pt
        \vrule width.#2pt}
      \hrule height.#2pt}}}
\def\square{\mathchoice\sqr{5.5}4\sqr{5.0}4\sqr{4.8}3\sqr{4.8}3}
\def\qed{\hskip4pt plus1fill\ $\square$\par\medbreak}

\def\cC{{\cal C}}
\def\cE{{\cal E}}

\def\cW{{\cal W}}
\def\C{{\bf C}}
\def\D{{\bf D}}
\def\N{{\bf N}}
\def\P{{\bf P}}
\def\R{{\bf R}}
\def\Z{{\bf Z}}
\def\cx#1{{\C}^{#1}}     
\def\px#1{{\P}^{#1}}     
\def\supp{{\rm supp}}
\def\deg{d}
\def\dim{k}
\def\id{{\rm id}}
\def\sumlyap{\Lambda}

%
%
%
\centerline{\bf REGULAR POLYNOMIAL ENDOMORPHISMS OF $\cx\dim$}
\bigskip
\bigskip

\def\efn{Supported in part by the NSF.}
\def\mfn{Supported in part by the Swedish Natural Science Council.}
\centerline{\tenrm%
  ERIC BEDFORD\footnote{*}{\efn} AND MATTIAS JONSSON\footnote{**}{\mfn}}
\bigskip
\bigskip
\bigskip
%
%
%

\section 0 Introduction

\noindent Let $f=(f_1,\dots,f_\dim):\cx\dim\to\cx\dim$ be
a mapping such that each
$f_j$ is a polynomial of degree $\deg$.  We consider the behavior of $f$ as a
dynamical system.
 That is, we consider the behavior of the iterates $f^n=f\circ\cdots\circ f$
as $n\to\infty$.  The points of most interest are those whose forward orbits
show recurrent behavior.  Thus we are led to focus on the set $K$ of points
whose forward orbits are bounded.  The point of this paper is to
present an approach to studying $K$ by working instead at the ``points at
infinity'' and then making a descent to $K$ via ``external rays.'' 

To illustrate this, let us consider the case $\dim=1$.  The recurrent dynamics of a
polynomial mapping $p:\C\to\C$ is carried by the set $K$, and the chaotic
dynamics take place on the Julia set 
$J:=\partial K$.  A polynomial mapping may be characterized as a holomorphic
mapping which extends continuously to the Riemann sphere.  The point at
infinity is completely invariant, and it is possible to find a holomorphic
function $\varphi$ defined in a neighborhood of $\infty$, called the B\"ottcher
coordinate, such that
$\varphi$ conjugates $p$ to the model mapping $\sigma(\zeta)=\zeta^\deg$ in a
neighborhood of infinity.  If the set
$K$ is connected, then $\varphi$ has an analytic continuation to a
conformal equivalence $\varphi:\C-K\to\C-\bar \D$ with the complement of the
closed unit disk.  Thus the dynamics of the restriction $p|(\C-K)$ is
conjugate to
$\sigma$ on all of $\C-\bar \D$.  If the inverse $\varphi^{-1}$ can be
extended continuously to $\partial\D$, then $p|J$ is represented as a
quotient of $\sigma|\partial\D$.  

Another point of view is to consider the point at infinity as the pole
for the Green function for $\C-K$.  This serves as the starting point for the
use of potential-theoretic methods in the study of polynomial mappings, as was
introduced by Brolin [Bro] and further developed by Sibony (see [CG]) and
Tortrat [T].  A natural way to descend from infinity to the set $J$ is to
follow the gradient lines of the Green function, which are also known as
``external rays.''  External rays were
introduced by Douady and Hubbard and have been developed into a powerful tool
for studying the relationship between the mappings
$p|(\C-K)$ and $p|J$. 

In our approach, we view $\cx\dim$ as an affine coordinate chart in $\px\dim$. 
Thus $\Pi:=\px\dim-\cx\dim$, the hyperplane at infinity, is isomorphic
to $\cx{\dim-1}$.  We study polynomial mappings $f$ of degree $\deg\ge2$ which
are regular, which means that $f$ extends
continuously (and thus holomorphically) to $\px\dim$.  It follows that $\Pi$
is completely invariant, and we let
$f_\Pi$ denote the induced dynamical system at infinity.  Further, $\Pi$ is
(super)-attracting in the normal direction, so the basin $A$ of points which
are attracted to $\Pi$ in forward time is an open set containing $\Pi$.  This
gives a completely invariant partition $\px\dim=K\cup A$.

We will make use of methods of the dynamics of holomorphic mappings of
$\px\dim$ and $\Pi\simeq\px{\dim-1}$.  In particular, we will make use of the approach introduced in
[HP] and developed more generally and systematically  in [FS1--3].  Namely,
there is an invariant current $T$ on $\px\dim$, and the exterior powers
$T^l:=T\wedge\cdots\wedge T$, $1\le l\le\dim$, are well defined, positive, closed
currents of bidegree $(l,l)$.
The supports
$J_l:=\supp(T^l)$ serve as a family of intermediate Julia sets.  
In the maximal degree we have measures  $\mu:=T^\dim$ on $\px\dim$ and
$\mu_\Pi:=T_\Pi^{\dim-1}$ (corresponding to $f_\Pi$) on $\Pi$.  We will denote their supports
by $J$ and $J_\Pi$, respectively.  The importance of 
$\mu$ is shown by the fact that it is ergodic (see Forn\ae
ss-Sibony [FS2--3]); it is balanced and has maximal entropy
($=\dim\log\deg$); and
$\mu$ describes the distribution of periodic points (see Briend [Bri]).  
The corresponding properties hold for $\mu_\Pi$ and $f_\Pi$.

Closely related to the current $T$ is the function $G$,
defined in~(1.1),  which measures the superexponential rate at which an orbit
approaches $\Pi$. The function $G$ coincides with the plurisubharmonic Green
function in $\cx\dim$ for the set $K$, and  $G$ has a logarithmic singularity along
$\Pi$.  In our passage from the case $\dim=1$ to $\dim>1$, we replace the point at
infinity by the set $J_\Pi$.  We replace the one-dimensional set $\C-K$ by the
set $A\cap\supp(T^{\dim-1})$.  The usual critical locus of $f$ is the
set $\cC$ where the Jacobian determinant of $f$ vanishes.  For the purpose of
dynamical study, we consider the critical measure 
$\mu_c:=[\cC]\wedge T^{\dim-1}$.  
This is motivated by another critical measure defined in [BS2].  A major theme
of this paper is to explore the interplay between the dynamical
measures $\mu$ and $\mu_\Pi$, the critical measure $\mu_c$, and the current
$T^{\dim-1}$.

As a first example we consider the Lyapunov exponents, which measure the
rate of expansion of the mapping $f$ with respect to a measure. 
We concern ourselves here with $\sumlyap(f)=\sumlyap(f,\mu)$, which
is defined as the sum of
all the Lyapunov exponents of $f$ with respect to $\mu$, and thus
$\sumlyap(f)$ measures the infinitesimal exponential rate of Lebesgue volume
expansion.  We show (Theorem~3.2) that the Lyapunov exponents of $(f,\mu)$
and $(f_\Pi,\mu_\Pi)$ are related by the formula
$$\sumlyap(f)=\sumlyap(f_\Pi)+\log\deg + \int G\,\mu_c.$$
This generalizes formulas of Przytycki [Pr] and Jonsson [J1] and is analogous to
a formula in [BS5].  We are grateful to B. Berndtsson for showing us how to
greatly simplify our previous proof.

The connection between $J_\Pi$ and $J$ is mediated by the stable manifolds
$$W^s(a)=\{x\in\px\dim:\lim_{n\to+\infty}\hbox{\rm dist}(f^na,f^nx)=0\}$$
and
$$W^s_{{\rm loc}}(a)=\{x\in\px\dim:\hbox{\rm
dist}(f^na,f^nx)<\delta,n\ge0\}$$ for $a\in J_\Pi$.  By Pesin Theory there
is a $\delta=\delta(a)>0$ for $\mu_\Pi$ a.e.~$a$ such that
$W^s_{{\rm loc}}(a)$ is a complex disk.  
In using the potential-theoretic approach, we work with a local stable
manifold as  a current of integration
$[W^s_{\rm loc}(a)]$.  (In general, a component of the global stable manifold
is not locally closed and does not have locally bounded area, so
$W^s(a)$ does not define a current of integration.)  A crucial property of
$T^{\dim-1}$ is that it has, in a measure-theoretic sense, a laminar structure on
$A$, made up of currents of integration over pieces of stable manifolds.  
The strongest example of laminarity is Theorem~5.1: if
$f_\Pi$ is uniformly expanding on $J_\Pi$, then
$T^{\dim-1}$ has the uniform laminar structure~(5.1) in a neighborhood
of $\Pi$.  (This result was also obtained by Peng [Pe].)  For a general map,
we establish two sorts of laminar structure (see Theorems~6.4 and~6.10). 

The set $W^s(a)$ is a manifold (with
singularities), which contains the Pesin disk $W^s_{\rm loc}(a)$, as
well as $W^s_{\rm loc}(b)$ whenever $b\in J_\Pi$ and
$f_\Pi^na=f_\Pi^nb$ for some $n\ge 1$.
The restriction $G|W^s(a)$ is harmonic for $\mu_\Pi$-a.e $a$ and
has no critical points on $W^s_{\rm loc}(a)$.
We define the set $\cE_a$ as the set of gradient lines of
$G|W^s(a)$ emanating from $a$ (the point at infinity). Thus $\cE_a$ is 
parametrized by the angle a gradient line makes inside
$W^s(a)$ at $a$. Let $\cE$ be the union of all $\cE_a$ over
$\mu_\Pi$-a.e.\ $a$. The elements of $\cE$ are called external rays.
It follows that the measure $\nu:=\mu_\Pi\otimes{d\theta\over2\pi}$
is well defined on $\cE$.
Using the (non-uniform) laminar structure in the general case, we are able to
show that for $\nu$ a.e.\ ray
$\gamma\in\cE$, there is a well defined endpoint
$e(\gamma)\in J$.  The basic connection between $\nu$ and $\mu$ is given by
transport along the external rays (Theorem~7.3):
$e_*\nu=\mu$.

In case the critical measure $\mu_c$ puts no mass near $\Pi$,
and $\dim=2$, the family of Pesin disks in fact extends to a Riemann
surface lamination near $\Pi$ (Theorem~8.8). 
To obtain results with better control over the geometry, we assume that 
$f_\Pi$ is uniformly expanding on $J_\Pi$; thus by the Stable Manifold Theorem,
$\cW^s_{\rm loc}(J_\Pi)=\{W^s_{{\rm loc}}(a):a\in J_\Pi\}$ is a Riemann
surface lamination for a uniform 
$\delta>0$.  We use this
lamination to obtain a local conjugacy with a model mapping.  Let
$C(J_\Pi)$ denote the set of complex lines through the origin in $\cx\dim$
corresponding to points of
$J_\Pi$. Our canonical model is given by the restriction of
$f_h$, the homogeneous part of $f$, to $C(J_\Pi)$.  If $f_\Pi$ is uniformly
expanding on $J_\Pi$, then the restriction $f_h|C(J_\Pi)$ is conjugate to the
restriction
$f|\cW^s_{\rm loc}(J_\Pi)$ in a neighborhood of $J_\Pi$ (Theorem~4.3).  This
conjugacy, which was also found by G.\ Peng [Pe], can be viewed
as a local B\"ottcher coordinate for $f$. A global
version is given in Theorem~7.4. 

The general approach of starting with
$f_\Pi$ to study $f$ is dual in spirit to that of Hubbard and Papadopol
[HP], who study a superattracting fixed point by working on the fiber
$\px{\dim-1}$ over the blowup.  From this point of view, Theorem~4.3 is
analogous to Theorem~9.3 of [HP].  A different approach was given by
S.\ Heinemann [H], who works
directly on $K$ without appeal to $\Pi$ or $f_\Pi$.

The final part of this paper addresses the question of when the
landing map $e:\cE\to J$ is continuous.  When possible, we would like to
replace statements about almost every ray with statements about every ray. 
In the case $\dim=1$, the continuity of $e$ is equivalent to local
connectedness of
$J$. In dimension
$\dim=2$, $\Pi$ is the Riemann sphere, and $f_\Pi$ is a rational map. If
$f_\Pi$ is hyperbolic, then $\cE$ is homeomorphic to $J_\Pi\times S^1$. 
In order to obtain the continuity of $e$, we assume that $f$ is Axiom A.
Thus the nonwandering set is the closure of the set of periodic points
and can be written as the union $S_0\cup S_1\cup S_2$,
where $S_j$ is a (uniformly) hyperbolic set with unstable index
$j$.  Our main result about the landing map is Theorem~10.2: { \sl If
$\dim=2$, if
$f$ is Axiom A, if
$f^{-1}S_2=S_2$, and if
$\mu_c(A)=0$, then $e:\cE\to J$ is a continuous surjection. }

The contents of this paper are as follows: Sections~1 and~2 recall some
basic notation and results, mostly on currents and potential theory. \S3
is devoted to the proof of the formula for the sum of the Lyapunov
exponents of $\mu$.  In~\S4 we discuss the Stable Manifold Theorem in the context of stable disks
through points of
$J_\Pi$.  We show that $f|W^s(J_\Pi)$ is conjugate to the canonical model
in a neighborhood of
$J_\Pi$.  In~\S5 we show that $T^{\dim-1}$ is laminar in a
neighborhood of $\Pi$ when
$f_\Pi$ is uniformly expanding on $J_\Pi$. This serves as a prototype for our
results in~\S6, where we discuss laminarity in the general case.  The
difficulties of~\S6 arise from two sources.  First is the fact that in
general (without uniform expansion on
$J_\Pi$) we need to use Pesin Theory to obtain our stable manifolds,
and the geometry of the Pesin disks is not controlled.  The second source
of difficulty, which is present already in the uniformly expanding case, 
comes from working globally in
$A$, since in general the global stable manifolds do not define currents of
integration.  In~\S7, we define external rays and show that $\nu$ almost every external
ray has a well-defined
landing point, and the endpoint map pushes the measure $\nu$ forward to
$\mu$.  We also give a global version (Theorem~7.4) of the local conjugacy
result of~\S4.  \S8 discusses the structure of the support
of $T^{\dim-1}$ inside $A$ for general maps $f$.  First we show that for
$\mu_\Pi$ a.e.~$a$, the global stable manifold $W^s(a)$ is dense in $A\cap
\supp(T^{\dim-1})$.  Then we show (Theorem~8.8) that if
$\mu_c(A)=0$, then the family of local Pesin disks in $A$ is
contained in a Riemann surface lamination by disks which are proper in
$A$.  The property Axiom A is discussed for endomorphisms of
$\cx2$ in~\S9, and~\S10 is devoted to the proof of Theorem~10.2, which gives the continuous
landing of the external rays.  Appendix A analyzes the behavior of the
homogeneous (canonical) model.  Appendix B serves as a reference for some
of the hyperbolicity results that are used in~\S8 and~9.

\par{\bigbreak\centerline{
    {\bf List of notation}}\nobreak\smallskip\noindent}

\newdimen\columnlength
\columnlength=30mm
\hbox{\hbox to \columnlength{
    $f$
    \hfill}\hbox{
    regular polynomial endomorphism of $\cx\dim$ of degree $\deg$.
    }}\smallskip\noindent
\hbox{\hbox to \columnlength{
    $f_h$
    \hfill}\hbox{
    homogeneous part of $f$ of degree $\deg$.
    }}\smallskip\noindent
\hbox{\hbox to \columnlength{
    $\Pi$
    \hfill}\hbox{
    hyperplane at infinity.
    }}\smallskip\noindent
\hbox{\hbox to \columnlength{
    $f_\Pi$
    \hfill}\hbox{
    restriction of $f$ to $\Pi$.
    }}\smallskip\noindent
\hbox{\hbox to \columnlength{
    $\cC$ 
    \hfill}\hbox{
    critical set of $f$.
    }}\smallskip\noindent
\hbox{\hbox to \columnlength{
    $\cC_\Pi$ 
    \hfill}\hbox{
    critical set of $f_\Pi$.
    }}\smallskip\noindent
\hbox{\hbox to \columnlength{
    $\cx\dim_*$
    \hfill}\hbox{
    $\cx\dim-\{0\}$.
    }}\smallskip\noindent
\hbox{\hbox to \columnlength{
    $\pi$
    \hfill}\hbox{
    projection of $\cx\dim_*$ on $\Pi$
    or $\cx{\dim+1}_*$ on $\px\dim$.
    }}\smallskip\noindent
\hbox{\hbox to \columnlength{
    $A$
    \hfill}\hbox{
    basin of $\Pi$.
    }}\smallskip\noindent
\hbox{\hbox to \columnlength{
    $K$
    \hfill}\hbox{
    complement of $A$.
    }}\smallskip\noindent
\hbox{\hbox to \columnlength{
    $G$
    \hfill}\hbox{
    Green function for $f$.
    }}\smallskip\noindent
\hbox{\hbox to \columnlength{
    $G_h$
    \hfill}\hbox{
    homogeneous Green function for $f_h$.
    }}\smallskip\noindent
\hbox{\hbox to \columnlength{
    $\tilde f$
    \hfill}\hbox{
    homogeneous map on $\cx{\dim+1}_*$.
    }}\smallskip\noindent
\hbox{\hbox to \columnlength{
    $\tilde G$
    \hfill}\hbox{
    homogeneous Green function for $\tilde f$.
    }}\smallskip\noindent
\hbox{\hbox to \columnlength{
    $\rho_G$
    \hfill}\hbox{
    Robin function for $G$.
    }}\smallskip\noindent
\hbox{\hbox to \columnlength{
    $T$
    \hfill}\hbox{
    invariant current for $f$.
    }}\smallskip\noindent
\hbox{\hbox to \columnlength{
    $T_h$
    \hfill}\hbox{
    homogeneous invariant current for $f_h$.
    }}\smallskip\noindent
\hbox{\hbox to \columnlength{
    $T_\Pi$
    \hfill}\hbox{
    invariant current for $f_\Pi$.
    }}\smallskip\noindent
\hbox{\hbox to \columnlength{
    $\mu$
    \hfill}\hbox{
    $T^\dim$.
    }}\smallskip\noindent
\hbox{\hbox to \columnlength{
    $\mu_\Pi$
    \hfill}\hbox{
    $T_\Pi^{\dim-1}$.
    }}\smallskip\noindent
\hbox{\hbox to \columnlength{
    $J$
    \hfill}\hbox{
    support of $\mu$.
    }}\smallskip\noindent
\hbox{\hbox to \columnlength{
    $J_\Pi$
    \hfill}\hbox{
    support of $\mu_\Pi$.
    }}\smallskip\noindent
\hbox{\hbox to \columnlength{
    $L_a$
    \hfill}\hbox{
    line in $\px\dim$ through 0 associated with $a\in\Pi$.
    }}\smallskip\noindent
\hbox{\hbox to \columnlength{
    $\sumlyap(f)$
    \hfill}\hbox{
    sum of Lyapunov exponents of $f$ with respect to $\mu$.
    }}\smallskip\noindent
\hbox{\hbox to \columnlength{
    $\sumlyap(f_\Pi)$
    \hfill}\hbox{
    sum of Lyapunov exponents of $f_\Pi$ with respect to $\mu_\Pi$.
    }}\smallskip\noindent
\hbox{\hbox to \columnlength{
    $\mu_c$
    \hfill}\hbox{
    critical measure: $\mu_c=[\cC]\wedge T^{\dim-1}$.
    }}\smallskip\noindent
\hbox{\hbox to \columnlength{
    $W^s_{\rm loc}(a)$
    \hfill}\hbox{
    local stable manifold at $a$.
    }}\smallskip\noindent
\hbox{\hbox to \columnlength{
    $A_0$
    \hfill}\hbox{
    subset of $A$ where $G>R_0$.
    }}\smallskip\noindent
\hbox{\hbox to \columnlength{
    $A_n$
    \hfill}\hbox{
    $f^{-n}(A_0)$.
    }}\smallskip\noindent
\hbox{\hbox to \columnlength{
    $W^s(a)$ 
    \hfill}\hbox{
    global stable manifold of $a$.
    }}\smallskip\noindent
\hbox{\hbox to \columnlength{
    $W^s(a,f_\Pi)$ 
    \hfill}\hbox{
    global stable manifold of $a$ with respect to $f_\Pi$.
    }}\smallskip\noindent
\hbox{\hbox to \columnlength{
    $W^s(J_\Pi)$
    \hfill}\hbox{
    stable set of $J_\Pi$ for $f$.
    }}\smallskip\noindent
\hbox{\hbox to \columnlength{
    $W^s(J_\Pi,f_h)$
    \hfill}\hbox{
    stable set of $J_\Pi$ for $f_h$.
    }}\smallskip\noindent
\hbox{\hbox to \columnlength{
    $W^s_0(a)$
    \hfill}\hbox{
    local stable disk for $f$ at $a\in J_\Pi$.
    }}\smallskip\noindent
\hbox{\hbox to \columnlength{
    $W^s_{0}(a,f_h)$
    \hfill}\hbox{
    local stable disk for $f_h$ at $a$ (subset of complex line).
    }}\smallskip\noindent
\hbox{\hbox to \columnlength{
    $\cW^s(J_\Pi)$ 
    \hfill}\hbox{
    stable lamination for $f$.
    }}\smallskip\noindent
\hbox{\hbox to \columnlength{
    $C(J_\Pi)$
    \hfill}\hbox{
    complex homogeneous cone over $J_\Pi$.
    }}\smallskip\noindent
\hbox{\hbox to \columnlength{
    $\Psi$ 
    \hfill}\hbox{
    conjugation of $f_h|W^s(J_\Pi,f_h)$ to $f|W^s(J_\Pi)$.
    }}\smallskip\noindent
\hbox{\hbox to \columnlength{
    $W^s_{-m}(a)$
    \hfill}\hbox{
    $W^s_{\rm loc}(a)\cap A_{-m}$.
    }}\smallskip\noindent
\hbox{\hbox to \columnlength{
    $Z_{a,n}$
    \hfill}\hbox{
    component of $W^s(a)\cap A_n$ containing $a$.
    }}\smallskip\noindent
\hbox{\hbox to \columnlength{
    $N_n(a)$
    \hfill}\hbox{
    number of points in $Z_{a,n}\cap\Pi$.
    }}\smallskip\noindent
\hbox{\hbox to \columnlength{
    $\cC_n$
    \hfill}\hbox{
    $\bigcup_{0\le j\le n-1}f^{-j}(\cC)$ (critical set of $f^n$).
    }}\smallskip\noindent
\hbox{\hbox to \columnlength{
    $\cC_{\infty}$
    \hfill}\hbox{
    $\bigcup_{n\ge 0}\cC_n$.
    }}\smallskip\noindent
\hbox{\hbox to \columnlength{
    $S_{a,n}$
    \hfill}\hbox{
    union of incomplete gradient lines in $Z_{a,n}$.
    }}\smallskip\noindent
\hbox{\hbox to \columnlength{
    $W_{a,n}$
    \hfill}\hbox{
    component of $Z_{a,n}-S_{a,n}$ containing $a$.
    }}\smallskip\noindent
\hbox{\hbox to \columnlength{
    $W_a$
    \hfill}\hbox{
    union of $W_{a,n}$ over $n\ge 0$.
    }}\smallskip\noindent
\hbox{\hbox to \columnlength{
    $\varphi_a$
    \hfill}\hbox{
    uniformizing map on $W_a$.
    }}\smallskip\noindent
\hbox{\hbox to \columnlength{
    $H_a$
    \hfill}\hbox{
    range of $\varphi_a$ (hedgehog domain).
    }}\smallskip\noindent
\hbox{\hbox to \columnlength{
    $\psi_a$
    \hfill}\hbox{
    inverse of $\varphi_a$.
    }}\smallskip\noindent
\hbox{\hbox to \columnlength{
    $\cE_a$
    \hfill}\hbox{
    set of external rays in $W_a$.
    }}\smallskip\noindent
\hbox{\hbox to \columnlength{
    $\cE$
    \hfill}\hbox{
    set of external rays.
    }}\smallskip\noindent
\hbox{\hbox to \columnlength{
    $\sigma$
    \hfill}\hbox{
    map on $\cE$.
    }}\smallskip\noindent
\hbox{\hbox to \columnlength{
    $\nu$
    \hfill}\hbox{
    invariant measure on $\cE$.
    }}\smallskip\noindent
\hbox{\hbox to \columnlength{
    $e_r$
    \hfill}\hbox{
    endpoint map on $\cE$ to level $G=r$.
    }}\smallskip\noindent
\hbox{\hbox to \columnlength{
    $e$
    \hfill}\hbox{
    endpoint map on $\cE$.
    }}\smallskip\noindent
\hbox{\hbox to \columnlength{
    $G_r$
    \hfill}\hbox{
    $\max(G,r)$.
    }}\smallskip\noindent
\hbox{\hbox to \columnlength{
    $S$
    \hfill}\hbox{
    union of incomplete gradient lines.
    }}\smallskip\noindent
\hbox{\hbox to \columnlength{
    $S_h$
    \hfill}\hbox{
    union of rays in $W^s(J_\Pi,f_h)$ corresponding to $S$.
    }}\smallskip\noindent
\hbox{\hbox to \columnlength{
    $\hat q$
    \hfill}\hbox{
    history of a point $q$, $\hat q=(q_i)_{i\le 0}$.
    }}\smallskip\noindent
\hbox{\hbox to \columnlength{
    $\hat f$
    \hfill}\hbox{
    shift map:$\hat f((q_i))=(q_{i+1})$.
    }}\smallskip\noindent
\hbox{\hbox to \columnlength{
    $S_i$
    \hfill}\hbox{
    union of basic sets of unstable index $i$ in $K$.
    }}\smallskip\noindent
\hbox{\hbox to \columnlength{
    $W^u_{\rm loc}(\hat q)$
    \hfill}\hbox{
    local unstable manifold at $\hat q$.
    }}\smallskip\noindent
\hbox{\hbox to \columnlength{
    $W^u(\hat q)$
    \hfill}\hbox{
    global unstable manifold at $\hat q$.
    }}\smallskip\noindent
\hbox{\hbox to \columnlength{
    $W^u(J)$
    \hfill}\hbox{
    backwards attracting basin for $J$.
    }}\smallskip\noindent
\hbox{\hbox to \columnlength{
    $W^u(S_1)$
    \hfill}\hbox{
    unstable set of $S_1$.
    }}\smallskip\noindent

%
%

\section 1 Regular polynomial endomorphisms and their Green functions

In the following two sections we summarize several basic results that
we will use. Additional details may be found in [HP], [FS1--3], and
[U]. We recommend the unified treatment in [FS3].  Throughout this
paper, we will let $f$ be a {\it regular polynomial endomorphism\/} of
$\cx\dim$ of degree $\deg\ge2$.  This means that the components of $f$ are
polynomials of degree $\deg$, and the homogeneous part $f_h$ of degree
$\deg$ satisfies $f^{-1}_h(0)=\{0\}$.  Alternatively, $f$ is regular if
and only if $\liminf|f(z)|/|z|^\deg>0$ as $|z|\to\infty$.  Such mappings are
also called strict polynomials by Heinemann [H].

Let $z=(z_1,\dots,z_\dim)$ denote (inhomogeneous) coordinates on $\cx\dim$,
and let $[z:t]=[z_1:\dots:z_\dim:t]$ denote homogeneous coordinates on
$\px\dim$. We fix the embedding of $\cx\dim$ into $\px\dim$ given
by $z\mapsto [z:1]$. Thus $\Pi=\{t=0\}$ corresponds to the hyperplane
at infinity, and $\Pi$ may be identified with with $\px{\dim-1}$
using homogeneous coordinates $[z]=[z_1:\dots:z_\dim]$. We equip 
$\px\dim$ with the Fubini-Study metric and measure distances and volumes
in that metric unless otherwise stated.

A regular polynomial endomorphism $f$ extends to a holomorphic endomorphism
of
$\px\dim$, still denoted by $f$, which may be defined by the formula
$f[z:t]=[t^\deg f(z/t):t^\deg]$. In fact, a holomorphic 
endomorphism of $\px\dim$ has
a completely invariant hyperplane exactly when it is conjugate to a
regular polynomial endomorphism of $\cx\dim$. We let $f_\Pi$ denote
the restriction of $f$ to $\Pi$. Under the identification
$\Pi\cong\px{\dim-1}$, $f_\Pi$ is given by $[z]\mapsto[f_h(z)]$.
Let $\cC$, $\cC_{\px\dim}$ and
$\cC_\Pi$ denote the critical sets of $f$ as a map of 
$\cx\dim$, $\px\dim$ and $\Pi$, respectively. Then we have
$\cC_{\px\dim}=\cC\cup\Pi$ and $\cC_\Pi=\overline{\cC}\cap\Pi$.

The model for our study of regular polynomial automorphisms is the
case when $f=f_h$ is a homogeneous mapping of $\cx\dim$. In this case
we write
$\cx\dim_*=\cx\dim-\{0\}$, and we let $\pi:\cx\dim_*\to\Pi$ be the
projection given by $\pi(z)=[z]$.  It is evident that $\pi$ gives a
semiconjugacy from $f_h$ to $f_\Pi$: $\pi\circ f_h=f_\Pi\circ\pi$.
In fact, $f_h$ is essentially a skew product over $f_\Pi$, as is shown
in Appendix A.

In the general case, we let $K$ be the compact set of points in $\cx\dim$
with bounded 
forward orbits and define $A:=\px\dim-K$. Thus $A$ is the basin of
attraction of $\Pi$. The function
$$
G(z)=\lim_{n\to\infty}\deg^{-n}\log^+|f^n(z)|\eqno(1.1)
$$ 
gives the (super-exponential) rate at which the orbit of $z\in\cx\dim$
approaches $\Pi$.  This is continuous and 
plurisubharmonic (psh) on $\cx\dim$ and
coincides with the pluri-complex Green function of $K$.
We will therefore also call
$G$ the Green function of
$f$. The homogeneous Green function for the homogeneous part $f_h$ of
$f$ of maximal degree $\deg$ is defined in an analogous way,
namely as
$$
G_h(z)=\lim_{n\to\infty}\deg^{-n}\log|f_h^n(z)|.
$$ 
The functions $G$ and $G_h$ are continuous on 
$\cx\dim$ and $\cx\dim_*$, respectively. We use $\log$ instead of
$\log^+$ so that $G_h$ is logarithmically homogeneous.

We may also define a homogeneous map $\tilde{f}$
on $\cx{\dim+1}$ by $\tilde{f}(z,t)=(t^\deg f(z/t),t^\deg)$.
The pair $(\tilde f,f)$ has properties analogous to those of $(f_h,f_\Pi)$.
The projection $\pi:\cx{\dim+1}_*\to\px\dim$ given by
$\pi(z,t)=[z:t]$ semiconjugates $\tilde{f}$ to $f$:
$f\circ\pi=\pi\circ\tilde{f}$.
We define the homogeneous Green function $\tilde G$ for $\tilde{f}$ by
$$
\tilde G(z,t):=\lim_{n\to\infty}\deg^{-n}\log|\tilde f^n(z,t)|\eqno(1.2)
$$
for $(z,t)\in\cx{\dim+1}_*$. 
The connection between $\tilde G$, $G$ and $G_h$ is
$\tilde G(z,1)=G(z)$
and $\tilde G(z,0)=G_h(z)$.
This leads us to the following asymptotic formulas for $G$ and $G_h$ 
near $\Pi$. Here $\rho_G$ denotes the Robin function of $G$ (cf.~[BT2]).
\proclaim Lemma~1.1. The asymptotics of $G$ and $G_h$ at $\Pi$ are given by
$$
\eqalign{
G_h(z)&=\log|z|+\rho_G[z]\cr
G(z)&=\log|z|+\rho_G[z]+o(1),}
$$
where $\rho_G$ is continuous on $\Pi$. Here 
$[z]=\pi(z)$ is the projection of $z$ on $\Pi$ defined above.

\give Proof. Since $G_h$ is homogeneous we have
$$
G_h(z)=\log|z|+G_h(z/|z|).
$$ 
Here the second term is continuous in $z$
and depends only on the projection $[z]$ of $z$ on $\Pi$. Hence
there exists a continuous function $\rho_G$ on $\Pi$ such that
$G_h(z/|z|)=\rho_G[z]$. This proves the first formula. 
To prove the second we write
$$
\eqalign{G(z)
&=\tilde{G}(z,1)\cr
&=\log|z|+\tilde{G}(z/|z|,0)+(\tilde{G}(z/|z|,1/|z|)-\tilde{G}(z/|z|,0))\cr
&=\log|z|+\rho_G[z]+o(1),}
$$
where the last line follows from the continuity of $\tilde{G}$ on
$\cx{\dim+1}_*$.
\qed  
The next result, implicitly contained in~[FS2], is crucial.
\proclaim Lemma~1.2.  If $M\subset\cx\dim$ is a complex manifold, and if
the iterates of $f$, restricted to $M$, are a normal family, then $G|M$ is
pluriharmonic on $M$.

\give Proof.  We may assume that $M\subset A$.  Passing to a subsequence
of the iterates $f^n=(f^n_{(1)},\dots,f^n_{(\dim)})$,  we may assume that
on $M$ we have 
$|f^n_{(j)}/f^n_{(1)}|$ bounded and $f^n_{(1)}\to\infty$ as
$j\to\infty$.  Thus $\log|f^n|=\log|f^n_{(1)}|+{1\over
2}\log\sum|f^n_{(j)}/f^n_{(1)}|^2$, so $\log|f^n|$ is written as
a pluriharmonic function plus something bounded.  Dividing by $\deg^n$ and
letting $n\to\infty$, we have that $G|M$ is pluriharmonic. \qed

%
%

\section 2 Invariant Currents

Here we assemble some basic facts about currents and give the definitions of
the invariant currents that are defined in terms of the Green
functions. Let $\Omega$ be a complex Hermitian manifold. If $M$ is a
positive current of bidimension $(p,p)$ on $\Omega$, then $M$ is
representable by integration.  This means that there is a total variation
measure $\Vert M\Vert$ and a measurable family of $(p,p)$-vectors $\vec m$ of
unit length with respect to the Hermitian metric, such
that we have the polar decomposition $M=\vec m \Vert M\Vert$.
In terms of a test form $\phi$, this means that $\langle
M,\phi\rangle=\int\langle\vec m(x),\phi(x)\rangle_x\,\Vert
M\Vert(x)$, where $\langle\cdot,\cdot\rangle_x$ denotes the pointwise pairing
between vectors and covectors.  This representation allows us to treat
positive currents as measures.

Let $Z$ is a closed subset of $\Omega$, and let $M$ be a positive
current of bidimension on $\Omega-Z$. If the total variation measure
$\Vert M\Vert$ has locally bounded mass near $Z$, then we may make the
{\it trivial extension\/} of
$M$ to $\Omega$ by extending the domain of definition of
$\Vert M\Vert$ to $\Omega$, and setting $\Vert M\Vert(Z)=0$.

For a general current $M$, we may define $M\contract\beta$, the contraction
with a smooth form
$\beta$, as the current which acts on a test form $\phi$ according to 
$\langle M\contract\beta,\phi\rangle=\langle M,\beta\wedge\phi\rangle$.  For
a positive current $M$ and a Borel set $S$, we may also define the restriction
$M\contract S$ by restricting $\Vert M\Vert$ to $S$.  This is the same as the
contraction of 
$M$ by the function which is 1 on $S$ and 0 elsewhere.

Let $\Omega$ and $\Omega'$ be complex manifolds and 
$g:\Omega'\to\Omega$ a holomorphic mapping. We define the pullback
$g^*S$ of a positive closed current $S$ on $\Omega$ in two cases.
The first is when $g$ is a submersion (eg.\ $g=\pi$):
$g^*S$ is then defined by integrating over the fibers
of $g$. The second case is when $g$ is a finite branched cover (eg.\ $g$ is a regular polynomial
endomorphism). If $S$ puts no mass on the critical image $g\cC$, then
$S$ coincides with the trivial extension to $\Omega$ of
$S\contract(\Omega-g\cC)$.
We define $g^*S$ to be the trivial extension to $\Omega'$ of
$(g|_{(\Omega'-g^{-1}g\cC)})^*(S\contract(\Omega-g\cC))$.

The two sorts of currents for which we will take pullbacks are as follows. 
First, if $[M]$ is a current of integration over a complex manifold, then our
definition gives $g^*[M]=[g^{-1}M]$ both in the case where $g$ is a
submersion and when $M\cap g\cC$ does not contain an open subset of $M$.
Second, if locally $S\le(dd^cu)^j$ for a bounded psh function $u$, then $S$
puts no mass on any complex analytic set and thus no mass on $g\cC$.
Similarly, $(dd^c(u\circ g))^j$ puts no mass on $g^{-1}g\cC$,
and $g^*(dd^cu)^j=(dd^c(u\circ g))^j$. Thus $g^*S$ is well defined.

Before we define currents on $\px\dim$, we recall the structure of
$\px\dim$ as a complex manifold.  If 
$U_j$ is the open subset of
$\px\dim$ where
$z_j\ne 0$, then 
$$
[z_1:\dots:z_\dim:t]\to(z_1/z_j,\dots,1,\dots,z_{\dim}/z_j,t/z_j) =
(\xi_1,\dots,1,\dots,\xi_{\dim+1})
$$
is a section of the bundle
$\pi:\cx{\dim+1}_*\to\px\dim$ over
$U_j$ and $\xi=(\xi_1,\dots,\xi_{\dim+1})$ (with the $j$th coordinate
missing) defines a biholomorphism between $U_j$ and $\cx\dim$. 

We define invariant currents on $\cx\dim$ by $T_{\cx\dim}:={1\over2\pi}dd^cG$ and
$T_{h,\cx\dim}:={1\over2\pi}dd^cG_h$.
To define the invariant currents on $\px\dim$, we define $g_j$ on the coordinate chart $U_j$ in
terms of the
$\xi$-coordinates by
$$G(\xi)=\log{1\over|\xi_{\dim+1}|}+g_j(\xi).$$
By Lemma~1.1, $g_j$ has a continuous extension from $U_j-\Pi$ to $U_j$.  A
property of the operator $dd^c$ is that $dd^cg_j$ can put no mass on a
pluripolar set if $g_j$ is bounded and psh. Since 
$dd^cG=dd^cg_j$ on $U_j-\Pi$, it follows that this formula defines a positive, closed current
on $U_j$, which coincides with the trivial extension of $T_{\cx\dim}$ from $\cx\dim\cap U_j$ to
$U_j$, and these definitions on $U_j$ fit together to give $T_{\px\dim}$. A similar formula
serves to define 
$T_{h,\px\dim}$ as a positive, closed current on the affine coordinate chart
$U_j$, and this coincides with the trivial extension of $T_{h,\cx\dim}$.
Since $T_{\px\dim}$ and $T_{h,\px\dim}$ are the trivial extensions of
$T_{\cx\dim}$ and $T_{h,\cx\dim}$, respectively, we will just denote these
currents as $T$ and $T_h$.

We recall that if $S=dd^cu$ is a positive closed current
of bidegree $(1,1)$ on a complex manifold with continuous potential
$u$, and $M$ is a complex submanifold, then the slice $S|_M$ is well-defined
and is equal to the current on $M$ defined by
$S|_M=dd^c|_M(u|_M)$.  We have seen that $T=T_{\px\dim}$ has a
local, continuous psh potential everywhere on $\px\dim$.  Thus we may define
the  current
$T_\Pi:=T|_\Pi$ as the slice current.   By Lemma~1.1, we also see
that $T_\Pi$ is also given as the slice of the
homogeneous current $T_\Pi=T_h|_\Pi$.

The positive closed
currents of bidegree $(1,1)$ on $\px\dim$ are characterized in terms
of logarithmically homogeneous psh functions on $\cx{\dim+1}$.  In
terms of this characterization, we have
$\pi^*T_{\px\dim}={1\over2\pi}dd^c\tilde G$, with $\tilde G$
from~(1.2), and  $T_\Pi$ is the unique positive,
closed current $T_\Pi$ on $\Pi$ such that
$\pi^*T_\Pi=T_{h}$.

Since the current $T$ has a locally defined continuous potential, we may 
define the exterior powers $T^l$ for $1\le l\le\dim$.  Using
the fact that for a bounded, psh function $g_j$ on $U_j$, $(dd^cg_j)^l$
puts no mass on a pluripolar set (and thus no mass on $\Pi$)  
for $1\le l\le\dim$, we see that the trivial extension of
$T_{\cx\dim}^l$ to $\px\dim$ is given by
$T_{\px\dim}^l$. Thus we may denote the exterior powers of our currents
simply as $T^l$ without ambiguity. The currents $T^l$,
$1\le l\le\dim$, are positive and closed on
$\px\dim$ and satisfy $f^*T^l=\deg^lT^l$.

The same arguments and properties, e.g.\ $f_h^*T_h^l=\deg^lT_h^l$, apply to $T_h$, with the
small complication that the potential $G_h=\log|z|+O(1)$ has a logarithmic singularity at the
origin.  Let us observe, however, that by a familiar calculation $(dd^c\log|z|)^l$ is equal to
$(2\pi)^\dim$ times the point mass at the origin if
$l=\dim$; and is a positive, closed current which is absolutely continuous with respect to
Lebesgue measure if
$l\le\dim-1$.  Now the Comparison Theorem (see~[BT1]) may be applied in the standard way to $G_h$
and
$\log|z|$, to conclude that $(dd^cG_h)^l$ is $(2\pi)^\dim$ times the point mass at the origin if
$l=\dim$, and puts no mass on the origin if $l<\dim$.   

Most important for us will be the currents 
$T^{\dim-1}$, $T_h^{\dim-1}$, of bidimension~(1,1), and $\mu:=T^\dim$ 
and $\mu_\Pi:=T_\Pi^{\dim-1}$, of bidimension (0,0).
Note that $\mu$ and $\mu_\Pi$ are represented by
probability measures on $\cx\dim$
and $\Pi$, respectively. We will denote their supports by
$J:=\supp(\mu)$ and $J_\Pi:=\supp(\mu_\Pi)$.

\give Remark. In the notation of [HP] the latter two sets would be called
$J_\dim$ and $J_{\Pi,\dim-1}$, respectively. We use $J$ and $J_\Pi$
for brevity, as we will not be using the other intermediate Julia sets.

\medskip
Currents that appear in complex dynamics often have a laminar structure.
Let $\Omega$ be a complex manifold with a Hermitian metric. Let $(A,\nu)$
be a measure space, and let $a\mapsto M_a$
denote a measurable family of positive currents on $\Omega$ with the
property that for every relatively compact domain
$\Omega_0\subset\Omega$ we have
$$
\int_{a\in A}\nu(a)\,\Vert M_a\Vert(\Omega_0)<\infty.\eqno(2.1)
$$
It follows from~(2.1) that we may define a positive current
$S=\int_{a\in A}\nu(a)\,M_a$, where the action on a test form $\phi$
is given by
$$
\langle S,\phi\rangle:=\int_{a\in A}\nu(a)\,\langle
M_a,\phi\rangle.
$$
We refer to $S$ as {\it laminar\/} if for almost every $a$ and $b$, either
$M_a=M_b$, or the supports of $M_a$ and $M_b$ are disjoint.
If in addition $a\to M_a$ is continuous, then the we say that
$S$ is {\it uniformly laminar\/}.
We note that if the currents $M_a$ are closed in $\Omega$, then so is $S$. 

A consequence of positivity is that the total variation measure $\Vert
M_a\Vert$ is equivalent to the measure $M_a\contract\beta^p$, where
$\beta$ is any  strictly positive~(1,1)-form, and $(p,p)$
is the bidimension of
$M_a$. Equivalent here means that the two measures are bounded above and
below by each other on compact subsets of $\Omega$, with the constant
depending only on $\beta$. Since $S\contract \beta^p=\int\,\nu(a)
M_a\contract\beta^p$, we conclude that the total variation measure of
$S$ is equivalent to the integral of the total variation measures of
the currents $M_a$:
$$
\Vert S\Vert\sim\int\nu(a)\,\Vert M_a\Vert.\eqno(2.2)
$$

The case that will appear in the sequel is where $M_a$
is the current of integration over a 1-dimensional complex variety in
$\Omega$. An example is the following result.

\proclaim Proposition~2.1. The following holds on $\px\dim$: 
$$
T^{\dim-1}_h=\int[L_a]\,\mu_\Pi(a),\eqno(2.3)
$$
where $L_a=\overline{\pi^{-1}(a)}$ is the complex line in $\px\dim$
passing through $a\in\Pi$ and the origin in $\cx\dim$.

\give Proof.
We know that $T_h^l$ puts no mass on the origin for $l<\dim$, and no mass on $\Pi$ in any case.
Since taking the wedge product commutes with taking pullback, the identity $T_h=\pi^*T_\Pi$
gives us 
$T_h^l=\pi^*(T_\Pi^l)$ on $\cx\dim_*$, and hence on $\px\dim$ if $l<\dim$.  Hence,  by the
definition of $\pi^*$ as integration over the fibers of
$\pi$, we have
$$
T_{h}^{\dim-1}
=\int\left[\pi^{-1}(a)\right]\,\mu_\Pi(a).
$$
Since $\{a\}$ and $\{0\}$ are sets of measure zero with 
respect to $L_a$, it follows that $[L_a]$ and 
$\left[\pi^{-1}(a)\right]$ define the same current on $\px\dim$.
 Therefore,
the equation above yields~(2.3).
\qed

Next we show that the laminar structure of a current is
preserved under wedge products.
Let $S$ and $X$ be positive closed currents on a complex Hermitian
manifold $\Omega$ and suppose that $X$ is of bidegree $(1,1)$.
For any ball
$\Omega_0$, there is a psh function $\psi$ on $\Omega_0$ such that we
may write $X=dd^c\psi$ on $\Omega_0$. In order to define $X\wedge S$,
it suffices to define the action on any test form $\phi$ on $\Omega_0$
as
$$
\langle X\wedge S,\phi\rangle:=\langle
\psi\,S,dd^c\phi\rangle.\eqno(2.4)
$$
If we have
$$
\int_{\Omega_0}|\psi|\,\Vert S\Vert<\infty,
$$
then~(2.4)
defines a positive, closed current. Note that, in this case, if
$\psi_m$ is a sequence of smooth, psh functions decreasing to $\psi$,
then $dd^c\psi_m\wedge S$ converges in the sense of currents to
$X\wedge S=dd^c\psi\wedge S$.

Let us recall the estimate of Alexander and Taylor [AT].
If $u$ is a bounded, psh function on
$\Omega$, then for $1\le j\le\dim$, the current $(dd^c u)^j$ has the
property that
$$
\int_{\Omega_0}|\psi| \,\Vert(dd^cu)^j\Vert<\infty
$$
for
any relatively compact domain $\Omega_0\subset\Omega$ and any psh
function $\psi$ on $\Omega$. Thus, if $S$ is a positive, closed
current such that locally there exists a bounded, psh function $u$
such that
$$
S\le (dd^cu)^j,\eqno(2.5)
$$
then the integral in formula~(2.4)
converges and defines $X\wedge S$ as a positive, closed current on
$\Omega$.

\proclaim Lemma~2.2. Let $M_a$, $a\in A$, be a measurable
family of positive, closed currents on $\Omega$, and let $\nu$ be a
measure on $A$ such that~(2.1) holds. Suppose, too, that the current
$S=\int_{a\in A}\nu(a)\, M_a$ satisfies~(2.5). If $X$ is a positive,
closed current of bidegree~(1,1) with local potential $\psi$, then for
$\nu$ almost every $a$, $X\wedge M_a$ is a well-defined positive
closed current on $\Omega$. Further, we have
$$
X\wedge S = \int_{a\in A} \nu(a)\,(X\wedge M_a),
$$
where $X\wedge S$ is defined according to~(2.4).

\give Proof. Let $\Omega'\subset\Omega$ be a domain where there is a
psh function $\psi$ with $dd^c\psi=X$. By~(2.2) and~(2.5), we
have
$$
\int\nu(a)\,(\int_{\Omega_0}|\psi|\,\Vert
M_a\Vert)\sim\int_{\Omega_0} |\psi|\,\Vert S\Vert<\infty\eqno(2.6)
$$
for every relatively compact $\Omega_0\subset\Omega'$. It follows that for
$\nu$ almost every $a$ we have $\int_{\Omega_0}|\psi|\,\Vert
M_a\Vert<\infty$ for all relatively compact $\Omega_0$, and thus
$X\wedge M_a$ is well defined.

Let $\psi_m$ denote a sequence of smooth, psh functions decreasing to
$\psi$. If $\phi$ is a test form on $\Omega_0$, then the smooth
current $X_m:=dd^c\psi_m$ satisfies
$$
\langle X_m\wedge S,\phi\rangle = \langle \psi_m S,dd^c\phi\rangle =
\int_{a\in A}\nu(a)\langle \psi_mM_a,dd^c\phi\rangle.
$$
The left hand side
converges to $X\wedge S$ as $m\to\infty$. For fixed $a$, the
integrand on the right hand side converges to $\langle \psi M_a,
dd^c\phi\rangle=\langle X\wedge M_a,\phi\rangle$ as $m\to\infty$.
Further, we have
$$
|\langle \psi_mM_a,dd^c\phi\rangle|\le C_\phi\int|\psi|\Vert M_a\Vert,
$$
where $C_\phi$ does not depend on $m$ or $a$, so
the lemma follows from~(2.6) and the Dominated Convergence Theorem.
\qed

We observe that taking pullbacks respects the laminarity of a
positive current $S$.

\proclaim Lemma~2.3. Let $\Omega,\Omega'$ be complex
manifolds and let $g:\Omega'\to\Omega$ be a branched covering.
Let $M_a$, $a\in A$, be a measurable family of positive
currents on $\Omega$, and let $\nu$ be a
measure on $A$ such that~(2.1) holds. Then
$$
g^*S = \int_{a\in A} \nu(a)\,g^*M_a.
$$

\give Proof. This follows immediately from the definition of $g^*S$ and
$g^*M_a$.
\qed

The object dual to the pullback is the pushforward. If $S$ is a positive,
closed current on $\Omega'$, then $g_*S$ is the current on $\Omega$
defined by $\langle g_*S,\varphi\rangle=\langle S,g^*\varphi\rangle$ for
all test forms $\varphi$. If $g$ is a finite branched covering,
then the pushforward may be thought of, locally, as the inverse of
the pullback so that $g_*S=\sum (g^{-1}_j)^*S$,
where the sum is taken over all the branches $g^{-1}_j$ of
$g^{-1}$. Applying this to the $\deg^\dim$ local inverses of a regular polynomial
endomorphism, we get that
$$
f_*T^l=\deg^{\dim-l}T,
{\rm\ \ and\ \ }f_*(GT^{\dim-1})=GT^{\dim-1}.\eqno(2.7)
$$

%
%

\section 3 Lyapunov exponents.

In this section we prove a formula for the sum of the Lyapunov
exponents of a regular polynomial endomorphism $f$ of $\cx\dim$.
The proof below, which depends significantly on ideas of Berndtsson [Be],
simplifies the argument in a previous version of this paper [J2].
A special case of~(3.1) below was proved in [J1]. 

Let us recall the notion of Lyapunov exponents. For more
details we refer to [Y]. The sum 
of the Lyapunov exponents of $f$ with respect to $\mu$ 
is the number $\sumlyap(f)=\sumlyap(f,\mu)$ given by
$$
\lim_{n\to\infty}{1\over n}\log|\det Df^n(x)|=\sumlyap(f),
$$
for $\mu$-a.e.\ $x\in\px\dim$. That this is well-defined is part of
the statement of Oseledec's Theorem. Hence $\sumlyap(f)$ measures
asymptotic infinitesimal Lebesgue volume growth of $f^n$
at $\mu$-a.e.\ point.
The individual Lyapunov exponents measure the asymptotic growth of the
derivative of $f^n$ in different directions; we will not give the
precise definition since we do not need it.

Our formula for $\sumlyap(f)$
will involve the integral of the Green function
against a critical measure so we begin by defining the latter 
measure as
$$
\mu_c
:=[\cC]\wedge T^{\dim-1}
={1\over 2\pi}dd^cH\wedge T^{\dim-1},
$$
where $H=\log|\det Df|$.

Then $\mu_c$ is a well-defined positive measure because 
$T$ has continuous local potentials, and the mass of
$\mu_c$ is finite (= the degree of the critical locus).

\proclaim Lemma~3.1. Let $f=f_h$ be any homogeneous regular polynomial 
endomorphism of $\cx\dim$, and let
$|\det(Df)|$ and $|\det(Df_\Pi)|$ be the Jacobians of $f$ and 
$f_\Pi$ in the Euclidean metric on
$\cx\dim$ and the Fubini-Study metric on $\Pi$, respectively.
Then
$$
|\det(Df)(z)|=\deg\cdot\left({{|f(z)|}\over{|z|}}\right)^\dim
|\det(Df_\Pi)[z]|.
$$

\give Proof. Pick any $z_0\in\cx\dim_*$. After pre- and
post-composing with dilations and unitary maps, we may assume that
$f(z_0)=z_0=(0,\dots,0,1)$.
Since $z_0$ and $[z_0]$ are now fixed points, the
choices of metrics are irrelevant when computing the Jacobians. We use
local coordinates $(\xi,s)$ on $\px\dim$ and $\xi$ on $\Pi$,
where $\xi_i=z_i/z_\dim$ for $1\le i\le\dim-1$ and $s=t/z_\dim$.
In these coordinates, the homogeneity of $f$ allows us to write
$$
\eqalign{
f(\xi,s)
&=(f_1(\xi,1)/f_\dim(\xi,1),\dots,f_{\dim-1}(\xi,1)/f_\dim(\xi,1),
s^\deg/f_\dim(\xi,1)),\cr
f_\Pi(\xi)
&=(f_1(\xi,1)/f_\dim(\xi,1),\dots,f_{\dim-1}(\xi,1)/f_\dim(\xi,1)).}
$$
Since the first $\dim-1$ coordinates in $f(\xi,s)$ do not depend
on $s$, we see that 
$$
\det Df(\xi,s)|_{(\xi,s)=(0,1)}=\deg\cdot\det Df_\Pi(\xi)|_{\xi=0}.
$$
We introduce the factors $|f(z)|^\dim$ and $|z|^\dim$ because of the
pre- and post-compositions with dilations. This completes the proof.
\qed

\proclaim Theorem~3.2. If $f$ is any regular polynomial endomorphism
of $\cx\dim$, then
$$
\sumlyap(f)=\log{\deg}+\sumlyap(f_\Pi)+\int G\,\mu_c.\eqno(3.1)
$$

\give Proof. From the Ergodic Theorem and the definition of $\mu$ we have
$$
\eqalign{
\sumlyap(f)
&=\int_{\cx\dim} H\,\mu\cr
&={1\over2\pi}\int_{\cx\dim} Hdd^cG\wedge T^{\dim-1}.}\eqno(3.2)
$$
Extend $G$ to a function $\hat G$ on $\px\dim$ by declaring
$\hat G=\infty$ on $\Pi$. Then $\hat G$ is locally integrable, and
$dd^c\hat G=2\pi(T-[\Pi])$.  Note that the current of integration over $\Pi$
appears because $\hat G$ has a $+\infty$ logarithmic singularity along $\Pi$. 
We have $H=\dim(\deg-1)\log|z|+O(1)$ near $\Pi-\cC_\Pi$, so
we define $\hat H$ on
$\cx\dim$ by $\hat H(z)=H(z)+\dim(\log^+|z|-\log^+|f(z)|)$. It follows
that $\hat H$
extends to a locally integrable function on $\px\dim$ which is
continuous outside $\cC\cup\cC_\Pi$.
Note that if $x\in\Pi$,
then $\hat H(x)$ depends only on the homogeneous part of $f$ of
degree $\deg$. Thus $\hat H=\log\deg+\log|\det Df_\Pi|$ on $\Pi$
by Lemma~3.1.

By the invariance of $\mu={1\over2\pi}dd^cG\wedge T^{\dim-1}$, we have
$$
\int_{\cx\dim}(\log^+|z|-\log^+|f(z)|)\,dd^cG\wedge T^{\dim-1}=0.
$$
Thus the last integral in~(3.2) equals
$$
{1\over 2\pi}\int_{\cx\dim}\hat H dd^cG\wedge T^{\dim-1}.
$$
Using the formula $dd^c\hat G=2\pi(T-[\Pi])$ and the fact that
$dd^cG\wedge T^{\dim-1}$ is supported on $\cx\dim$ we see that this equals
$$
{1\over2\pi}\int_{\px\dim}\hat Hdd^c\hat G\wedge T^{\dim-1}
+\int_\Pi\hat H\,\mu_\Pi.\eqno(3.3)
$$
By Stokes' Theorem the first term in~(3.3) is
$$
\eqalign{
{1\over2\pi}\int_{\px\dim}\hat G\,dd^c\hat H\wedge T^{\dim-1}
&={1\over2\pi}\int_{\cx\dim}G\,dd^c\hat H\wedge T^{\dim-1}\cr
&={1\over2\pi}\int_{\cx\dim}G\,
dd^c(H+\dim(\log^+|z|-\log^+|f(z)|))\wedge T^{\dim-1}\cr
&={1\over2\pi}\int_{\cx\dim}G\,dd^c H\wedge T^{\dim-1}\cr
&=\int G\,\mu_c.}
$$
The first equality holds because
$dd^c\hat H\wedge T^{\dim-1}$ puts no mass on $\Pi$, and the third
equality follows because $f_*(GT^{\dim-1})=GT^{\dim-1}$ by (2.7),
so 
$$
\eqalign {
\int_{\cx\dim}(dd^c\log^+|f(z)|)\,GT^{\dim-1}
&=\int_{\cx\dim}(f^*dd^c\log^+|z|)\,GT^{\dim-1} \cr
=\int_{\cx\dim}(dd^c\log^+|z|)\,f_*(GT^{\dim-1})
&=\int_{\cx\dim}(dd^c\log^+|z|)\, GT^{\dim-1}.}
$$
Further, by the remark above the second term in~(3.3) is equal to
$$
\int_\Pi(\log\deg+\log|\det Df_\Pi|)\,\mu_\Pi
=\log\deg+\sumlyap(f_\Pi).
$$
This completes the proof.
\qed

\proclaim Corollary~3.3. $\sumlyap(f)\ge{\dim+1\over 2}\log\deg$.

\give Proof. It is a result of Briend [Bri] that the Lyapunov exponents
of $f_\Pi$ with respect to $\mu_\Pi$ are bounded below by
${1\over 2}\log\deg$ (see~\S6). In particular,
$\sumlyap(f_\Pi)\ge{\dim-1\over 2}\log\deg$, so Theorem~3.2 gives
$\sumlyap(f)\ge\log\deg+{\dim-1\over 2}\log\deg={\dim+1\over 2}\log\deg$.
\qed

%
%

\section 4 Stable manifolds and a local model near $\Pi$.

In~\S5--\S7 we will show how the current $T^{\dim-1}\contract A$
has a laminar structure, and how this allows us to describe $\mu$
in terms of external rays. The laminar structure is easier to handle
---and visualize---in
the case when $f_\Pi$ is uniformly expanding on $J_\Pi$.
The expansion enables us to invoke the (uniform) stable manifold
theorem, and thus provides us with a continuous family of
local stable manifolds. These manifolds define currents of integration,
which are responsible for the laminar structure of
$T^{\dim-1}$ in a neighborhood of $\Pi$.

In this section we prove a version of the stable manifold theorem.
We also show that $f$, restricted to the union of the local stable
manifolds, is conjugate to the homogeneous map $f_h$. The homogeneous
map $f_h$ may be seen as the canonical model for $f$ near $\Pi$ and
the conjugacy as a generalization of
the B\"ottcher coordinate at infinity for one-dimensional polynomials.

Let us therefore assume that $f_\Pi$ is 
uniformly expanding on $J_\Pi$. This means
that there exist constants $c>0$ and $\lambda>1$ such that
$$ 
|Df^n_xv|\ge c\lambda^n|v|
\quad x\in J_\Pi,\ v\in T_x\Pi,\ n\ge 1.\eqno(4.1)
$$
If $f$ is expanding on $J_\Pi$ and $a\in J_\Pi$, then the tangent space
$T_a{\px\dim}$ splits into a direct sum $E^u(a)\oplus E^s(a)$, 
where $E^u(a)=T_a\Pi$ and $E^s(a)$ is the eigenspace of $Df_a$ 
associated with the zero eigenvalue. We clearly have
$Df_a(E^{u/s}(a))\subset E^{u/s}(f_\Pi a)$, and $E^{u/s}(a)$
depends continuously on $a$. Therefore, with the definition
given in Appendix B,  $f_\Pi$ is hyperbolic on $J_\Pi$.

The stable manifold theorem (see [Ru, p.~96] or [PS, Theorem~5.2]) 
asserts that there is a {\it local stable manifold} $W^s_{\rm loc}(a)$
at each point of $a$ in $J_\Pi$. This is defined by
$$
W^s_{\rm loc}(a):=\{x\in\px\dim: d(f^jx,f^ja)<\delta
{\rm\ for\ all\ }j\ge 0\}\eqno(4.2)
$$
for small $\delta>0$, and is an embedded real 2-dimensional
disk. In fact, since $f$ is holomorphic, $W^s_{\rm loc}(a)$
is a complex disk, i.e.\ the image of an injective holomorphic
immersion of $\D$. Moreover, the local stable manifolds depend
continuously on $a$ in the $C^1$ topology.

It will be convenient to work with 
neighborhoods of $\Pi$ defined in terms of the Green function,
so let $A_0:=\{G>R_0\}$ and $A_n=f^{-n}A_0=\{G>\deg^{-n}R_0\}$, where
$R_0>0$ and $n\in\Z$.
Thus $A_n\subset A_{n+1}$, $\bigcap_nA_n=\Pi$ and $\bigcup_nA_n=A$.

Since $W^s_{\rm
loc}(a)\cap\Pi=\{a\}$, it follows from Lemma~1.2 that $G|W^s_{\rm loc}(a)$
is harmonic  on the complement of $a$ and equal to $+\infty$ at $a$. 
If we choose
$R_0$ greater than the maximum of $G$ on $\partial W^s_{\rm loc}(a)$,  it
follows from the maximum principle that
$$
W^s_0(a):=W^s_{\rm loc}(a)\cap A_0
$$
is a properly embedded disk in $A_0$ for all $a\in J_\Pi$. We
call $W^s_0(a)$ the {\it local stable disk} at $a$.

We also define global stable manifolds by
$$
W^s(a)=\{x\in\px\dim: d(f^jx,f^ja)\to 0{\rm\ as\ }j\to\infty\}.
$$
In contrast to the diffeomorphism case, the global
stable manifolds may have singular points.
Notice also that $W^s(a)$ contains all the local stable manifolds
$W^s_{\rm loc}(b)$ for $b\in J_\Pi$ with $f_\Pi^nb=f_\Pi^na$,
$n\ge 0$. We will in fact prove that $W^s(a)$ is dense
in the support of $T^{\dim-1}\contract A$ (see Corollary~8.5).
The global stable manifolds may have infinitely many
components or be connected (see the example following Proposition~4.2).

In Theorem~4.3 we will show that $f$ is conjugate to $f_h$. In the
proof of this theorem we will use a holomorphic homotopy between
$f$ and $f_h$, defined by $f_\tau=f_h+\tau(f-f_h)$ for $\tau\in\C$.
Note that $\tau=0$ and $\tau=1$ correspond to $f_h$ and $f$,
respectively, and 
that the restriction of $f_\tau$ to $\Pi$ is $f_\Pi$ for all
$\tau$. Hence there are local stable manifolds for $f_\tau$ for
all $\tau$. We will need to control the dependence of these
manifolds on $\tau$. To get this control, we prove 
a version of the Stable Manifold Theorem adapted to our
situation.

It will be natural to consider the {\it stable set of $J_\Pi$}, i.e.\
$$
W^s(J_\Pi)=W^s(J_\Pi,f)=
\{x\in\px\dim: d(f^nx,J_\Pi)\to 0{\rm\ as\ }j\to\infty\}.
$$
Note that the expansion of $f_\Pi$ on $J_\Pi$ and the superattracting
nature of $\Pi$ implies that $W^s(J_\Pi)$ is closed in $A$ and that
$W^s(J_\Pi)\cap\Pi=J_\Pi$.

Let $G_\tau$ be the Green function for $f_\tau$, 
$A_{0,\tau}=\{G_\tau>R_0\}$, etc.

\proclaim Theorem~4.1. If $\delta$ is small enough and 
$R_0$ is large enough, then for all
$\tau$ with $|\tau|<2$ the following holds
\item{(1)} $W^s_{\rm loc}(a;f_\tau)$ is proper in $A_{0,\tau}$
  for $a\in J_\Pi$ and
  $W^s_{0,\tau}(a):=W^s_{\rm loc}(a;f_\tau)\cap A_{0,\tau}$
  is a properly embedded disk in $A_{0,\tau}$.
\item{(2)} $W^s_{0,\tau}(a)$ is the connected component of
  $W^s(a,f_\tau)\cap A_{0,\tau}$ containing $a$. In particular,
  $W^s_{0,\tau}(a)$ does not depend on the choice of $\delta$.
\item{(3)} $W^s_{0,\tau}(a)$ depends continuously on $a$ and
  holomorphically on $\tau$.
\item{(4)} $G_\tau$ is harmonic and has no critical points on 
  $W^s_{0,\tau}(a)$.

\give Proof. To avoid cumbersome notation we will write $f$
instead of $f_\tau$. However, it is important that the constructions
below hold uniformly in $\tau$ (for $|\tau|<2$).

Our first task is to define good coordinate charts.
Pick $a=[a_1:\dots:a_\dim]\in J_\Pi$. After a unitary change
of coordinates we may assume that $a=[0:\dots:0:1]$.
Let $\zeta=(\zeta_1,\dots,\zeta_{\dim-1})$, where
$\zeta_j=z_j/z_\dim$ and let $t=1/z_\dim$.
We denote the ball 
$|\zeta|<\epsilon_1$ by $U_a=U_a(\epsilon_1)$, the disk 
$|t|<\epsilon_2$ by $V_a=V_a(\epsilon_2)$ and the box
$U_a\times V_a$ by $B_a=B_a(\epsilon)=B_a(\epsilon_1,\epsilon_2)$
for $\epsilon_1,\epsilon_2>0$.
Note that $\Pi$ corresponds to $\{t=0\}$ and the line $L_a$
to $\{\zeta=0\}$. Also,
the Euclidean metric on $B_a$ and the Fubini-Study metric on
$\px\dim$ differ by at most a multiplicative constant $C\ge1$
and $C$ is close to one if $\epsilon_1$ and 
$\epsilon_2$ are small.

We may find an iterate $f^N$, such that~(4.1)
holds with $n=N$, $c=1$ and $\lambda=3C$. Thus,
if $\epsilon_1$ and $\epsilon_2$ small enough, then we have
\item{(i)} If $a,b\in J_\Pi$ and $a\ne b$, then there is an $n\ge 0$
  such that $d(f_\Pi^na,f_\Pi^nb)>3C\epsilon_1$.
\item{(ii)} If $a$,$b\in J_\Pi$, $a\ne b$ and
$f_\Pi^Na=f_\Pi^Nb$, then $\overline{B_a}\cap\overline{B_b}=\emptyset$.
\item{(iii)} If $a\in J_\Pi$, then $f^N$ has no critical points
  in $\overline{B_a}-\Pi$.

We define a vertical disk in $B_a$ to be a disk of the form
$\{\zeta={\rm const}\}$ and a vertical-like disk to be the graph
of a holomorphic map $U_a\to V_a$. Similarly we define
horizontal and horizontal-like disks (these have codimension 1).

By choosing $1\gg\epsilon_1\gg\epsilon_2>0$,
we get that for all $a\in J_\Pi$ and for all $f=f_\tau$ with $|\tau|<2$:
\item{(iv)} $f^N(B_a)\cap B_{f_\Pi^Na}
  \subset U_{f_\Pi^Na}\times{1\over 2}V_{f_\Pi^Na}$.
\item{(v)} $f^{-N}(B_{f_\Pi^Na})\cap B_a
  \subset{1\over 2}U_a\times V_a$.
\item{(vi)} If $\Sigma$ is a horizontal disk in $B_a$, then
  $f^N(\Sigma)\cap B_{f_\Pi^Na}$ is a horizontal-like disk
  in $B_{f_\Pi^Na}$ and the restriction of $f^N$ to
  $\Sigma\cap f^{-N}(B_{f_\Pi^Na})$ is a biholomorphism.

Conditions (iv)--(vi) are illustrated in Figure 1 (with N=1).

\centerline{\psfig{figure=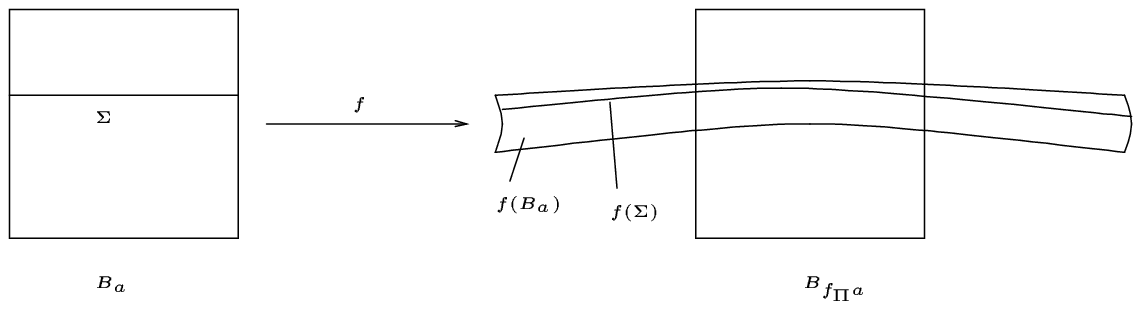,width=0.9\hsize}}
\centerline{Figure 1.}
\medskip

To produce stable manifolds we have to iterate backwards. We claim that
\item{(vii)} If $\Sigma'$ is a vertical-like disk in $B_{f_\Pi^Na}$, then
  $f^{-N}(\Sigma')\cap B_a$ is a vertical-like disk in $B_a$.

To see this, note that $f^{-N}(\Sigma')\cap B_a$ is an 
analytic variety in
$B_a$. Let $\Sigma$ be a horizontal disk in $B_a$. We claim that
$f^N(\Sigma)$ intersects $\Sigma'$ in exactly one point. Indeed, by
(v) we may write $f^N(\Sigma)\cap B_{f_\Pi^Na}=\{t=g(\zeta)\}$
and $\Sigma'=\{\zeta=h(t)\}$, where $g$ and $h$ are holomorphic.
Hence the intersection between these two sets is the
unique fixed point of the holomorphic map 
$g\circ h:U_{f_\Pi^Na}\to{1\over2}U_{f_\Pi^Na}$.
By (v) it follows that $\Sigma$ intersects 
$f^{-N}(\Sigma')\cap B_a$ in exactly one point. This proves that 
the latter set is a vertical-like disk.

Now define $B_a^n=B_a\cap f^{-N}(B_{f_\Pi^Na})\cap\dots
\cap f^{-nN}(B_{f_\Pi^{nN}a})$
for $n\ge 0$ and $B_a^\infty=\cap_{n\ge 0}B_a^n$
Using the Kobayashi metric on $U_a$, 
it follows from (v) and (vi) that there is a constant $\kappa>0$,
independent of $a$ and $n$, such that the diameter of 
$\Sigma\cap B_a^n$ is less that $\kappa2^{-n}$ for every horizontal
disk $\Sigma$ in $B_a$.
We claim that $B_a^\infty$ is a vertical-like disk in $B_a$.
Indeed, the estimate above 
implies that $\Sigma\cap B_a^\infty$ consists of at most
one point for every $\Sigma$. 
On the other hand, a repeated application of
(vii) shows that the set $\gamma_n(a)$, defined inductively
by $\gamma_0(a)=\{0\}\times V_a$ and
$\gamma_n(a)=f^{-N}(\gamma_{n-1}(f_\Pi^Na))\cap B_a$, is a
vertical-like disk in $B_a$, contained in $B_a^n$.
Hence any limit of a subsequence $\gamma_{n_j}(a)$ is a vertical-like
disk in $B_a$. By the remark above, this disk must be exactly
$B_a^\infty$. We also see that $\gamma_n(a)$, and hence $B_a^\infty$,
depends holomorphically on $\tau$.

Clearly $f(B_a^\infty)\subset B_{f_\Pi a}^\infty$ for all $a\in J_\Pi$.
We claim that the sets $B_a^\infty$ are
pairwise disjoint. Suppose $a\ne b$. By~(i)
there exists an $n\ge 0$ such that
$f_\Pi^nb\notin3B_{f_\Pi^na}$. Hence 
$B_{f_\Pi^na}^\infty\cap B_{f_\Pi^nb}^\infty=\emptyset$,
so $B_a^\infty$ and $B_b^\infty$ are disjoint.

We next show that the disks $B_a^\infty$ depend continuously on $a$.
Note that if $\epsilon_2'<\epsilon_2$ and
$\epsilon'=(\epsilon_1,\epsilon_2')$, then
$B_a^\infty(\epsilon')$ is the 
restriction to $V_a(\epsilon_2')$ of the vertical-like disk defining
$B_a^\infty(\epsilon)$.
Pick any $\epsilon_2'<\epsilon_2$ and $M$ be larger than
the Lipschitz constant for all $f_\tau$ on $\px\dim$. Assume that
$b$ is close to $a$ and choose $n$ maximal so that
$M^{nN}C^3d(a,b)<
(\epsilon_1^2+\epsilon_2^2)^{1/2}-(\epsilon_1^2+\epsilon_2'^2)^{1/2}$.
Then a simple calculation shows that 
$B_b^\infty$ is contained in $B_a^n$, and the latter set intersects
every horizontal disk in a set of diameter at most $\kappa2^{-n}$. Hence
$B_a^\infty$ depends continuously on $a$.

It follows from the superattractive nature of $\Pi$ that if
$\epsilon_2$ is small enough, then
$d(fx,fa)<d(x,a)$ whenever $a\in J_\Pi$ and $x\in B_a^\infty$.
Hence, if $\delta>0$ is small enough, then
$W^s_{\rm loc}(a)=\{x\in B_a^\infty: d(x,a)<\delta\}$.
Thus $W^s_{\rm loc}(a)$ is a
complex disk, compactly contained in $B_a$, depending
continuously on $a$ and holomorphically on $\tau$.  Thus $W^s_0(a)$ depends
continuously on $a$ and holomorphically on $\tau$.

For $a\in J_\Pi$ in a neighborhood of $a_0\in J_\Pi$, let
$\zeta$ denote a local holomorphic coordinate for $W^s_{\rm loc}(a)$ such
that 
$\zeta=0$ corresponds to $a$.  By Lemma~1.1 and Lemma~1.2,
the restriction $G|W^s_{\rm loc}(a)$ has
the form $\log|\zeta|+g_a(\zeta)$, where $g_a$ is bounded and harmonic in
a neighborhood of $\zeta=0$.  It follows that for $\zeta$ sufficiently
small, 
$G|W^s_{\rm loc}(a)$ has no critical points.  Thus for $R_0$ sufficiently
large, $G|W^s_0(a)$ has no critical points. 

We claim that if $a\in J_\Pi$, then
$$
f^{-1}W^s_0(a)\cap A_0=\bigcup_{f_\Pi b=a}W^s_0(b).\eqno(4.3)
$$
Indeed, $X:=f^{-1}W^s_0(a)\cap A_0$ is a subvariety of $A_0$. Any component
of $X$ must meet $\Pi$, and this must happen at a point
in $f_\Pi^{-1}a$. But then a neighborhood of $\Pi$ in $X$
is contained in the union of $B_b^\infty$, where $b\in f_\Pi^{-1}a$.
Thus~(4.3) holds. Repeating the same argument, we see that
$$
f^{-j}W^s_0(a)\cap A_0=\bigcup_{f_\Pi^jb=a}W^s_0(b).\eqno(4.4)
$$

It remains to be seen that $W^s_0(a)$ is the connected component
of $W^s(a)\cap A_0$ containing $a$.
We may assume that $R_0$ is so large that
$W^s(J_\Pi)\cap A_0$ is contained in the union of the boxes $B_a$.
In particular $W^s(a)\cap A_0$ is contained in the union these boxes
for all $a\in J_\Pi$.

Thus, by the definition of $B_a^\infty$, we see that
if $x\in W^s(a)\cap A_0$, 
then $f^{nN}(x)\in B_{f_\Pi^{nN}a}^\infty$ for large $n$.
Thus~(4.4) implies that $x\in\bigcup_{f_\Pi^{nN}b=a}W^s_0(b)$.
Hence we have shown that
$$
W^s(a)\cap A_0=\bigcup_{b\in W^s(a,f_\Pi)}W^s_0(b),\eqno(4.5)
$$
where $W^s(a,f_\Pi)=\cup_{j\ge 0}f_\Pi^{-j}f_\Pi^na$.
Since the disks $W^s_0(a)$ are disjoint, it follows
that $W^s_0(a)$ is the connected component of $W^s(a)\cap A_0$
containing $a$.
\qed

\proclaim Proposition~4.2. For $R_0$ large enough we have
$$
W^s(J_\Pi)\cap A_0=\bigcup_{a\in J_\Pi}W^s_0(a).
$$

We let $\cW^s(J_\Pi)$ denote the partition of $W^s(J_\Pi)$ by
global stable manifolds.  Proposition~4.2 implies that 
$\cW^s(J_\Pi)\cap A_0$ is a Riemann surface
lamination (see [C] for the definition).
Now the iterates of $f$ are local biholomorphisms 
outside the set $\cC_\infty:=\bigcup_{n\ge 0}f^{-n}(\cC)$.
The expansion of $f_\Pi$ on $J_\Pi$ implies that 
$\cC_\infty\cap W^s(J_\Pi)$ is closed and nowhere dense 
in $W^s(J_\Pi)$. Thus $\cW^s(J_\Pi)-\cC_\infty$ is also a
lamination. In Figure 2 we see how the local stable disks in $A_0$ can join
at higher levels and create a lamination whose leaves are not simply
connected.

\centerline{\psfig{figure=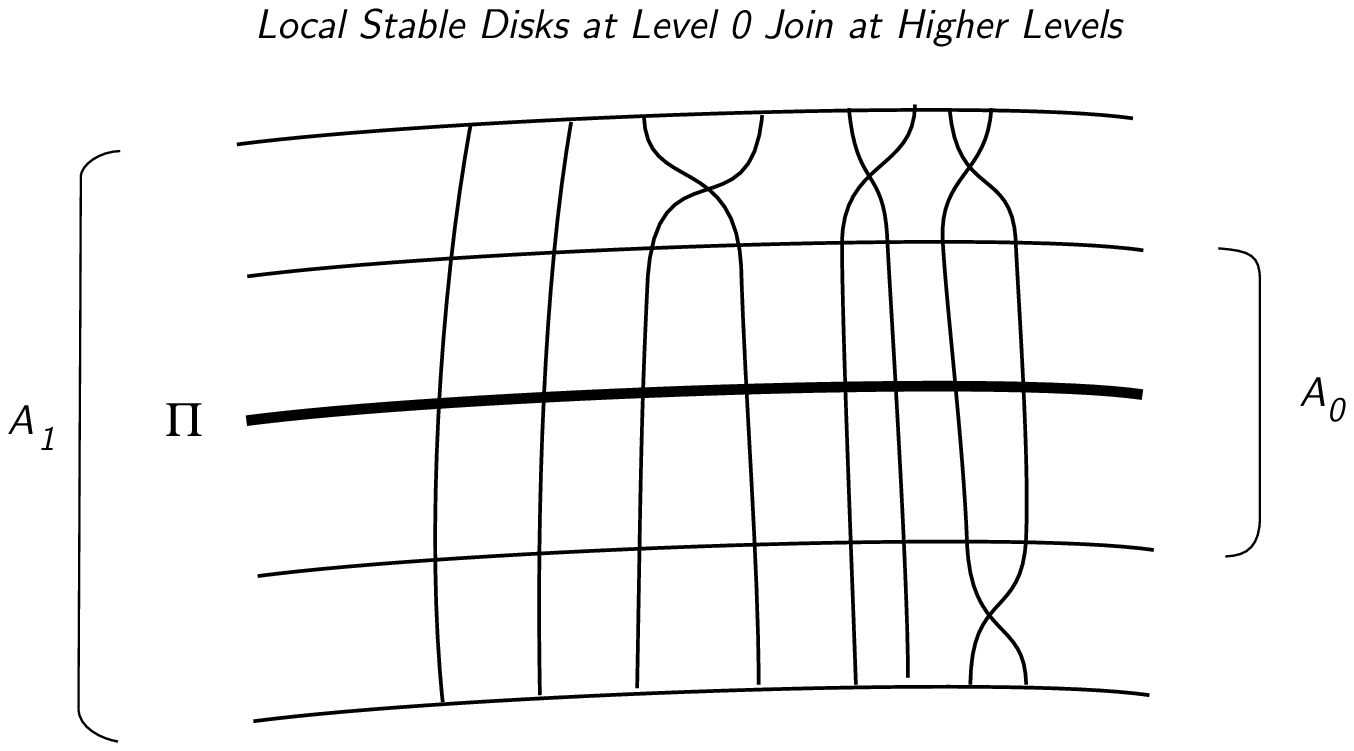,width=0.7\hsize}}
\centerline{Figure 2.}
\medskip

\give Proof of Proposition~4.2.  This can be proved by observing that the
natural extension
$\widehat{J_\Pi}$ has local product structure (see Proposition~B.6),
but we will give a direct proof. The inclusion ``$\supset$'' is trivial.
After replacing $f$ by an iterate we may assume that~(4.1) holds
with $n=c=1$ and $\lambda=3$.
Let $M\ge 1$ be larger than the 
Lipschitz constant for $f$ on $\px\dim$. Let $\eta>0$ be so small
that if $a\in J_\Pi$, then all branches of 
$f_\Pi^{-1}$ are single-valued on the ball
$B(f_\Pi a,4M\eta)$ in $\Pi$ and the branch mapping $f_\Pi a$ 
to $a$ maps $B(f_\Pi a,4M\eta)$ into the ball $B(a,2M\eta)$.
Now let $x\in W^s(J_\Pi)\cap A_0$. Let $n$ be so large that 
$d(f^{n+j}x,J_\Pi)<\eta$ for $j\ge 0$ and pick points 
$a_j\in J_\Pi$ such that $d(f^{n+j}x, a_j)<\eta$ for $j\ge 0$.
Then $(a_j)_{j\ge 0}$ is an $2M\eta$-pseudoorbit in $J_\Pi$, i.e.\
$d(f_\Pi a_j,a_{j+1})<2M\eta$.
Let $g_j$ be the branch of $f_\Pi^{-1}$ on 
$B(f_\Pi a_j,4M\eta)$
mapping $f_\Pi a_j$ to $a_j$. 
Then $g_j(a_{j+1})\in B(f_\Pi a_{j-1},4M\eta)$
so the point $b^{(j)}:=g_0\circ\dots\circ g_j(a_{j+1})$ is well-defined.
Moreover $d(f_\Pi^i(b^{(j)}),a_i)<2M\eta$ for $0\le i\le j$. Letting 
$j\to\infty$ and using the compactness of $J_\Pi$ we find a
point $b\in J_\Pi$ such that 
$d(f_\Pi^ib,a_i)<3M\eta$ for all $i\ge 0$. Hence 
$d(f^{n+i}x,f_\Pi^ib)<4M\eta$ for all $i\ge 0$.
Assume that $4CM\eta<\epsilon$, with $C$ and $\epsilon$ from the proof
of Theorem~4.1. It follows that
$f^nx\in W^s_0(b)$, so by~(4.4) we have $x\in W^s_0(c)$ for some
$c\in f_\Pi^{-n}b$. This completes the proof.
\qed

\give Remark. Proposition~4.2 holds for $f_\tau$ for $|\tau|<2$
(with a uniform $R_0$).
\bigskip
The following family of examples shows that the behavior of the lamination $\cW^s(J_\Pi)$ can be
simple when the leaves are simply connected, and complicated otherwise.

\give Example.  Let $f(z,w)=(z^2+c,w^2)$.  We use the affine coordinate
$\zeta=w/z$ on $\Pi$, so
that
$f_\Pi(\zeta)=\zeta^2$, and $J_\Pi=\{|\zeta|=1\}$.  Let $K_c$ denote the
(1-dimensional) filled
Julia set of $p_c(z)=z^2+c$, and let $G_c(z)=\log|z|+o(1)$ denote the
Green function for $K_c$. 
It follows that $G(z,w)=\max(G_c(z),\log^+|w|)$, and
$W^s(J_\Pi)=\{G_c(z)=\log^+|w|>0\}$.  We may
choose a harmonic conjugate function $G_c^*(z)$
such that $\phi_c(z):=\exp(G_c(z)+iG_c^*(z))\approx z$
is single-valued and analytic for $z$ large.  For $\zeta\in
J_\Pi$, the local stable manifold $W^s_0(\zeta)$ is given
(for $z$ large) as the graph $w=\zeta\phi_c(z)$. If
$c$ belongs to the Mandelbrot set, then $K_c$ is connected, and $\phi_c$
extends analytically to
$\C-K_c$.  The stable manifolds are then countable unions of closed
disks in $A$, each of which is a graph of the form
$\{w=\zeta\phi_c(z):z\in\C-K_c\}$, $\zeta\in J_\Pi$.

In case $K_c$ is not connected, we let
$\Phi:\C\times\cx*\to\C\times\cx*$ be the biholomorphic mapping given by
$(u,v)=\Phi(z,w)=(z/w,1/w)$.  Thus $h:=\Phi\circ f\circ \Phi^{-1}$ is
given by
$h(u,v)=(u^2+cv^2,v^2)$, which is a homogeneous mapping.
Let $G_h(u,v)$ denote the
logarithmically homogeneous Green function associated to $h$. Thus
$$
G_h(u,v)=G_h(uw,vw)-\log|w|=G_h(z,1)-\log|w|=G_c(z)-\log|w|.
$$
The set $\{G_h<0\}$ is the basin of attraction of $(0,0)$ for the mapping
$h$, and the boundary of the basin is given by $\{G_h=0\}$.

The set $\{G_c(z)=G_c(u/v)>0,G_h(u,v)=0\}\subset\{G_h=0\}$ is 
the image of
$W^s(J_\Pi)=\{G_c(z)>0, G_c(z)-\log|w|=0\}$ under the mapping $\Phi$.
Thus $\Phi$ transfers the Riemann surface lamination of the one set to the
lamination of the other. The set
$\{G_c(z)>0,G_h(u,v)=0\}$ is exactly that part of the boundary of the
basin where $G_h$ is pluriharmonic, and thus it part of the boundary
that lies inside the Fatou set.
The leaves of the Riemann surface lamination of this set have been
studied by Hubbard and Papadopol [HP] and Barrett [Ba] and are shown to be
dense and to have infinite topological type as well as other complicated
behaviors.

\bigskip
\vbox{
  \hbox to \hsize{
    \hfil
    \psfig{figure=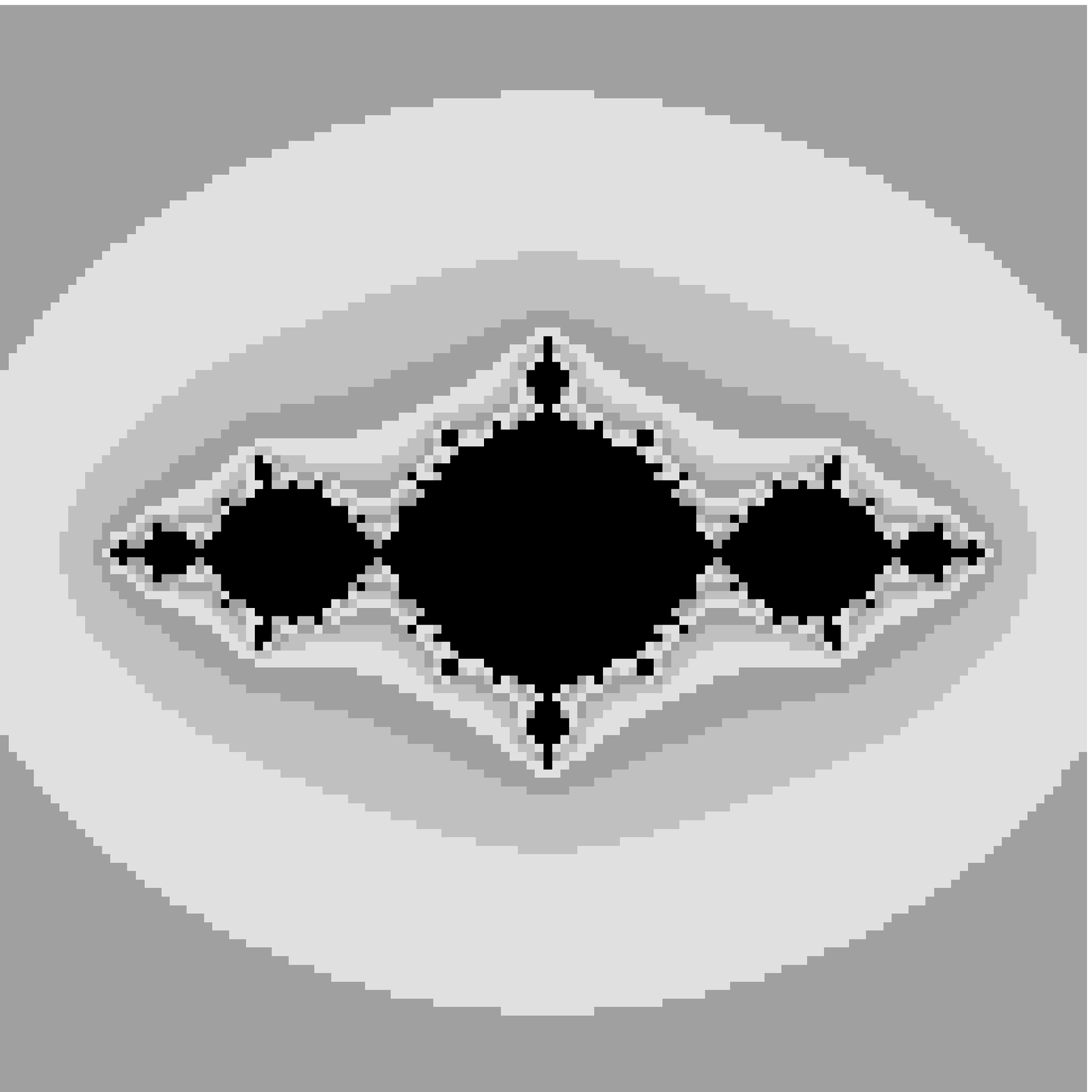,width=0.32\hsize}
    \hfil
    \psfig{figure=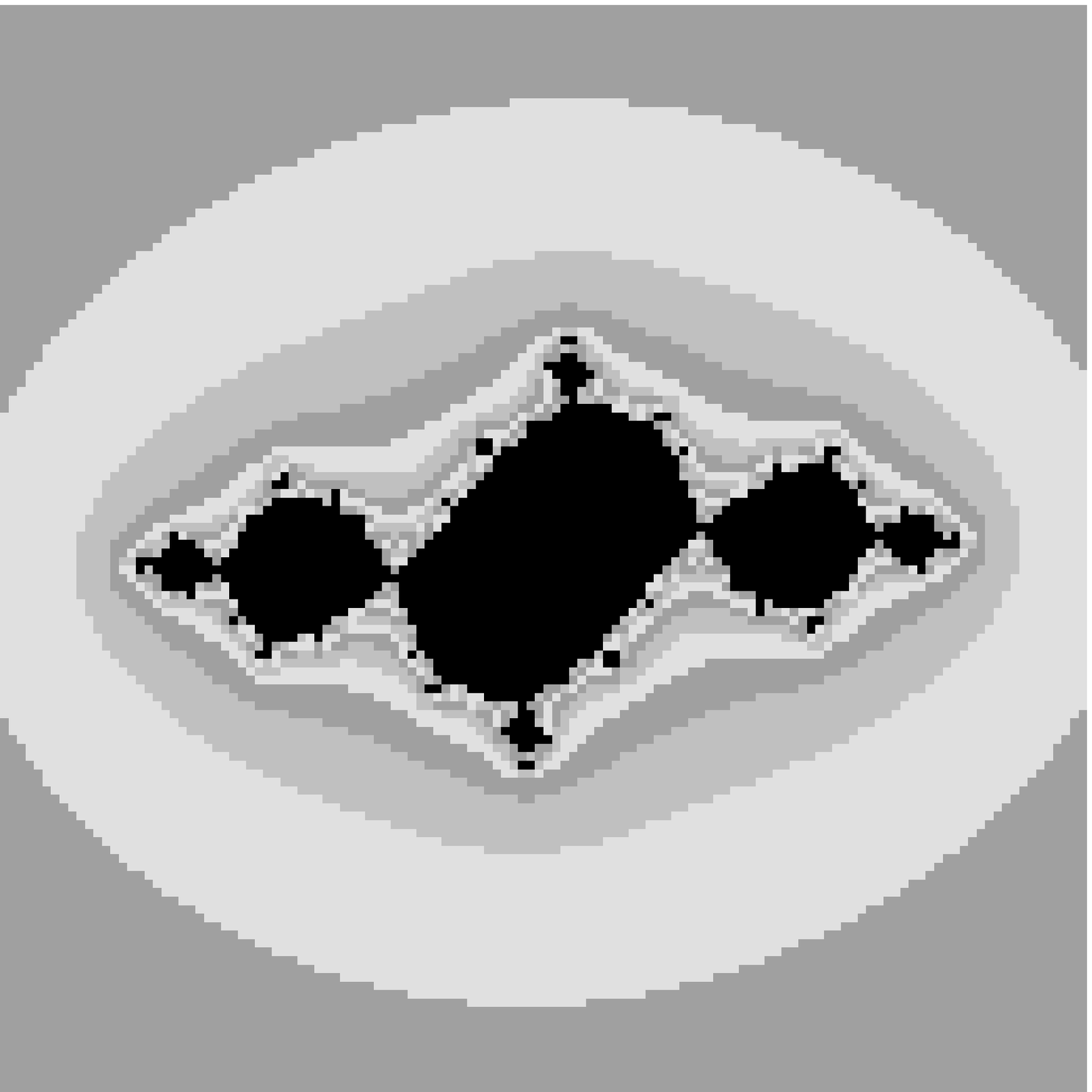,width=0.32\hsize}
    \hfil
    \psfig{figure=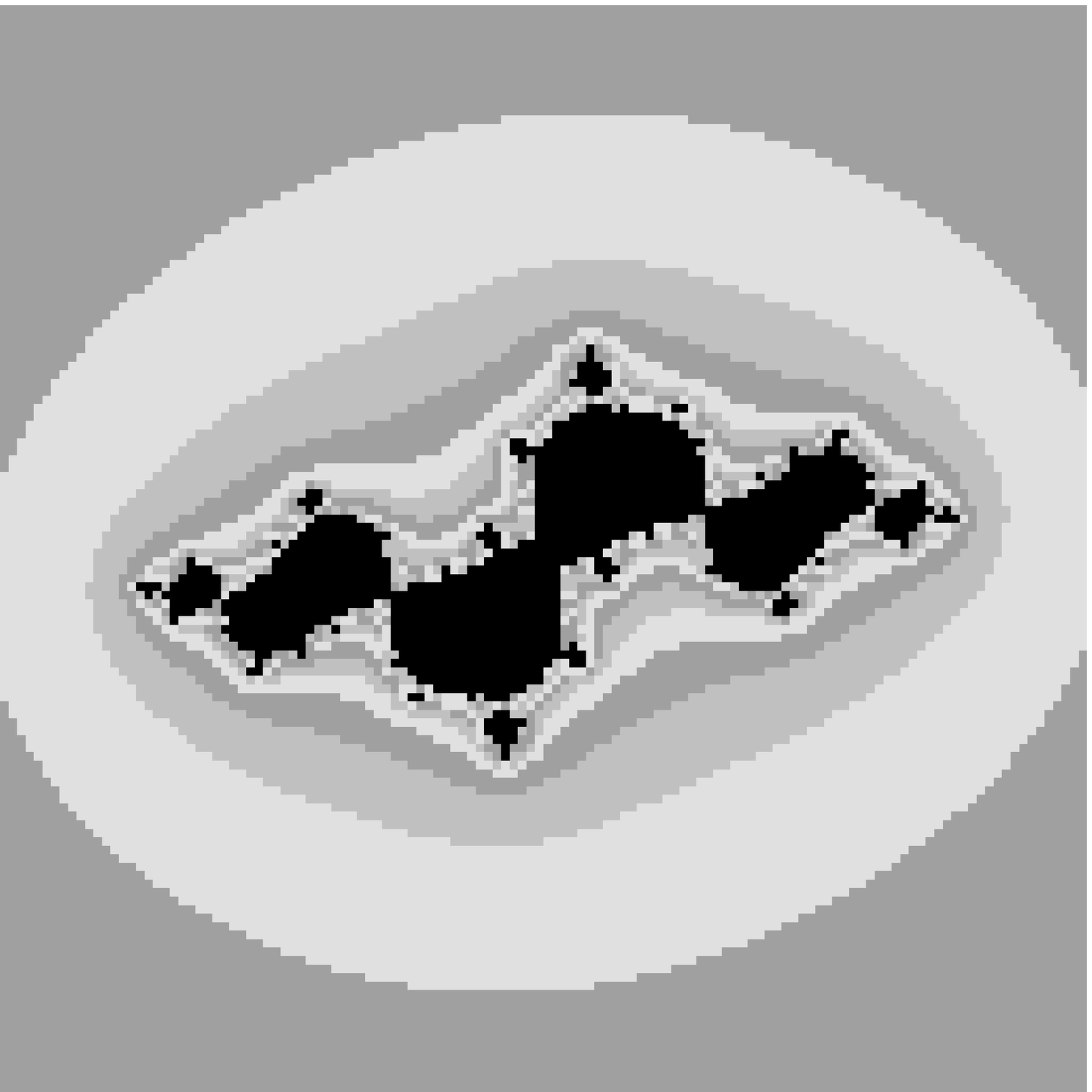,width=0.32\hsize}
    \hfil
  }
  \hbox to \hsize{
    \hfil
    \psfig{figure=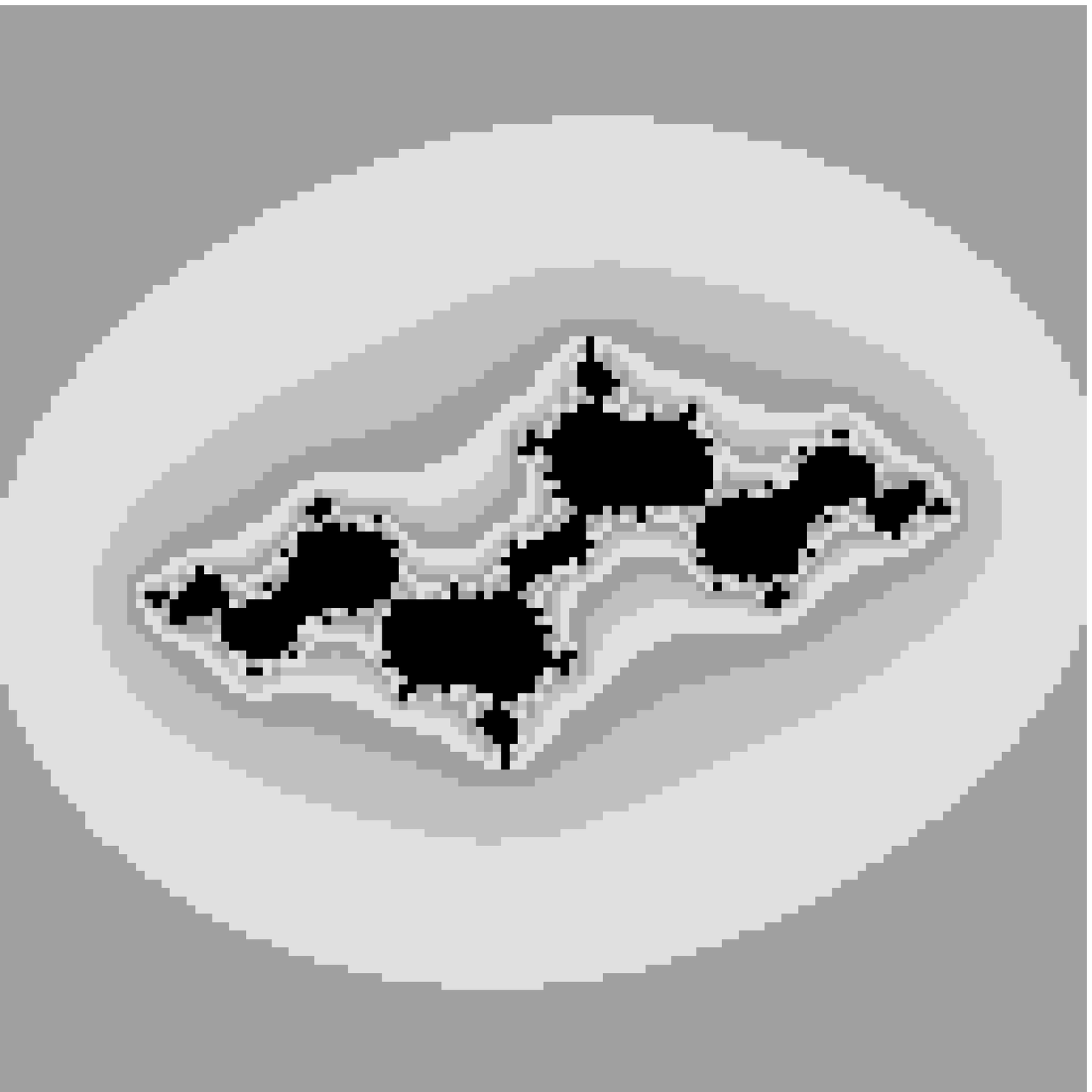,width=0.32\hsize}
    \hfil
    \psfig{figure=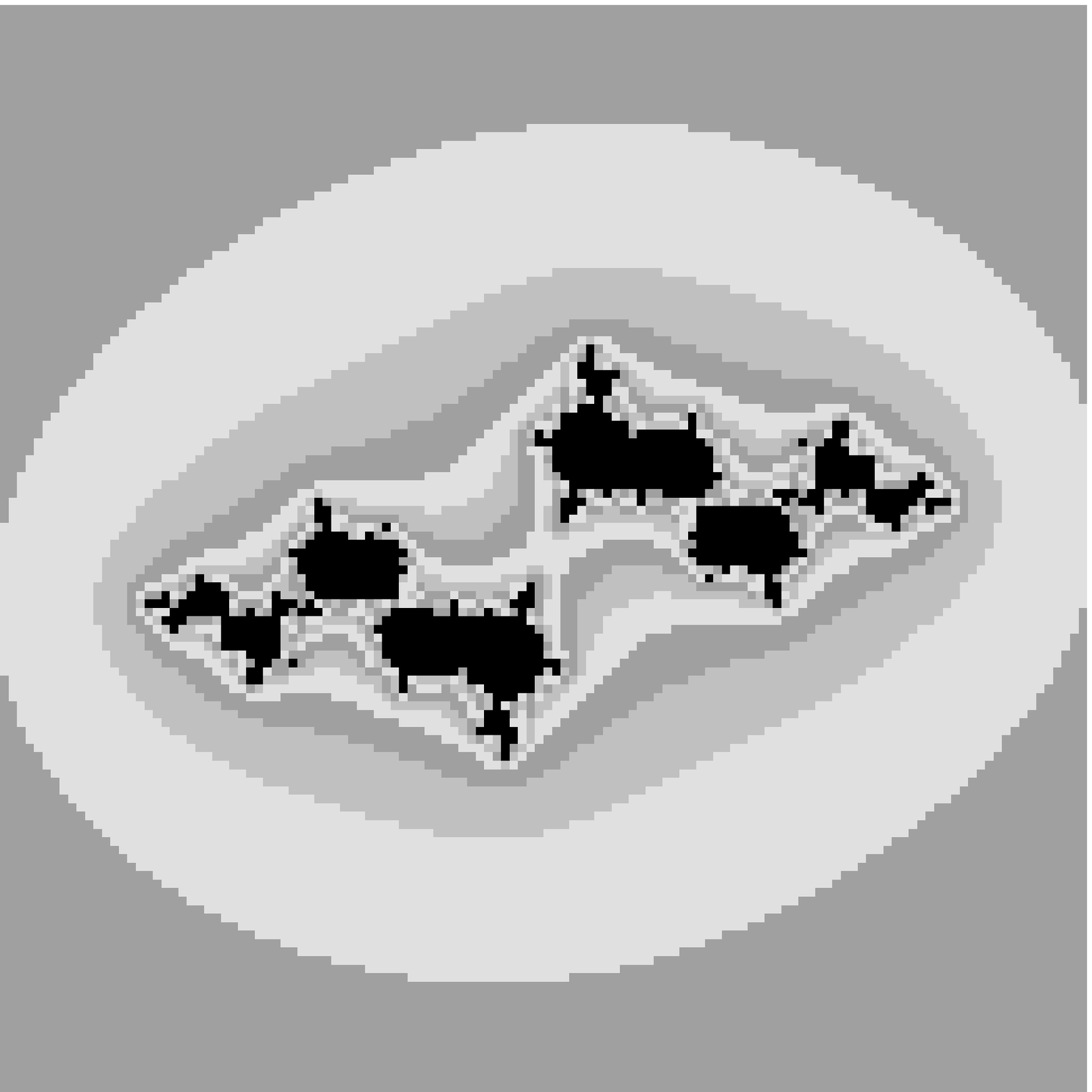,width=0.32\hsize}
    \hfil
    \psfig{figure=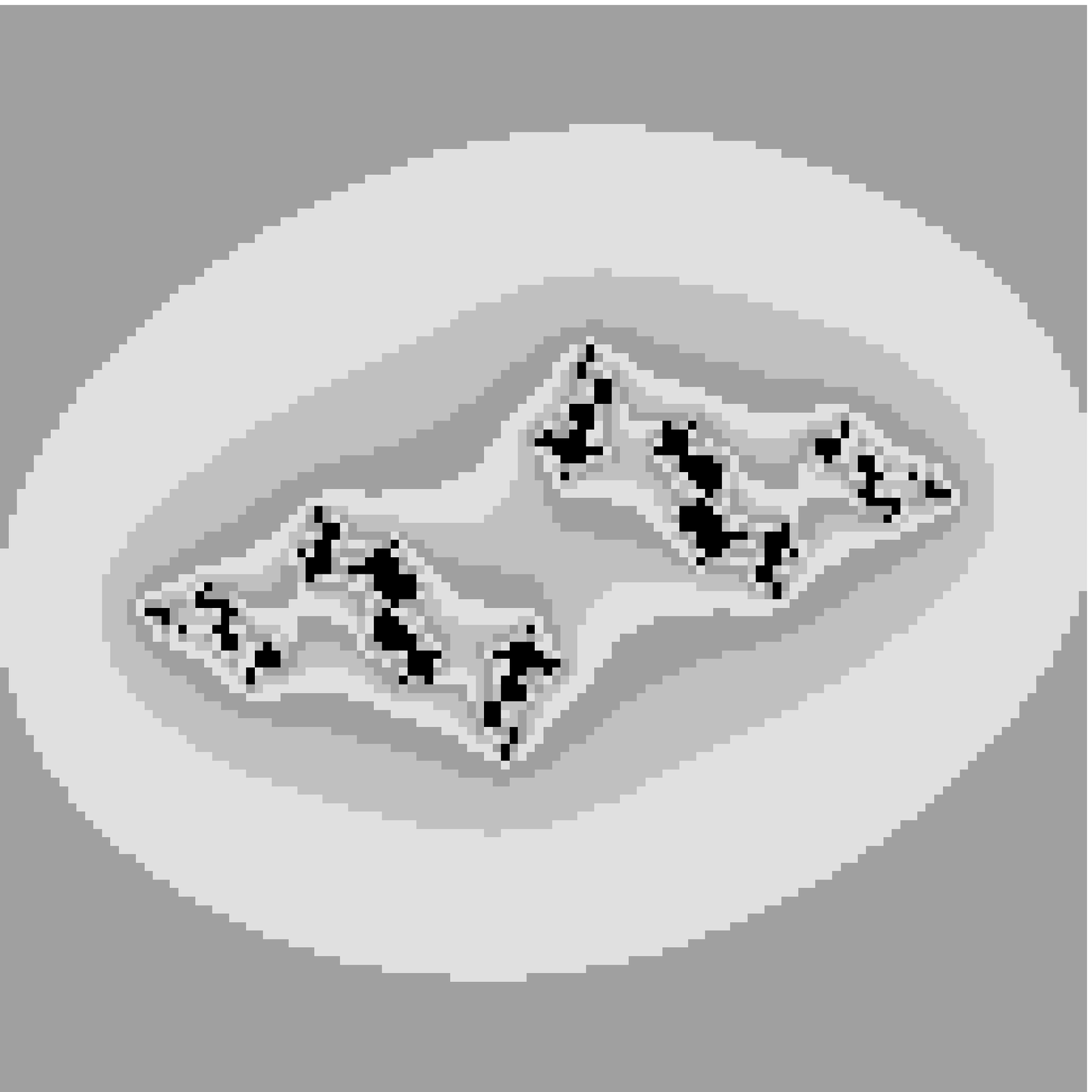,width=0.32\hsize}
    \hfil
  }
}
\medskip
\centerline{Figure 3.}
\medskip

Figure 3 shows slices of
$W^s(J_\Pi)$ by complex lines $\{z=c\}$ for the map
$f(z,w)=(z^2-0.1,w^2-z^2+0.2z-0.5i)$.
In the coordinate $\zeta=w/z$, we have $f_\Pi(\zeta)=\zeta^2-1$.
The first picture is the Julia set of the map $\zeta\to\zeta^2-1$.
By Proposition~4.2, the slices above
converge (suitably scaled) to this picture as $c\to\infty$.
The remaining five pictures show the slices by the lines 
$\{z=2\}$, $\{z=1.3\}$, $\{z=1.2\}$, $\{z=1.15\}$, $\{z=1.1\}$.

We will now show that $f$, restricted to 
$W^s(J_\Pi)\cap A_0$, is conjugate to the canonical model
$f_h$. Note that if $f=f_h$ is homogeneous, then
$W^s(J_\Pi)=C(J_\Pi)\cap A$, where $C(J_\Pi)$ is the complex cone of
lines $L_a$, for $a\in J_\Pi$. Let $A_{0,h}:=\{G_h>R_0\}$.

\proclaim Theorem~4.3. Suppose that $f_\Pi$ is uniformly expanding
on $J_\Pi$.  If $R_0$ is large enough, then there is a homeomorphism
$\Psi:W^s(J_\Pi,f_h)\cap A_{0,h}\to W^s(J_\Pi,f)\cap A_0$
conjugating $f_h$ and $f$. Further, $G\circ\Psi=G_h$ and 
the restriction of $\Psi$ to the local stable disk
$W^s_0(a,f_h)$ is a biholomorphism onto
$W^s_0(a,f)$ for all $a\in J_\Pi$.

\give Remark.  Using the Green functions $G_h$ and $G$, it is easy to
construct a biholomorphism of
$W^s_0(a,f_h)$ onto $W^s_0(a,f)$ taking $G_h$ to $G$. Such a biholomorphism is
unique up to a rotation of the disk $W^s_0(a,f_h)$. The difficulty
in constructing $\Psi$ is to chose these rotations in a continuous way.
Note that the disk $W^s_0(a,f)$ will not, in general, be tangent to 
$W^s_0(a,f_h)$ at
$a$. On the other hand, if $f=f_h$, then we may use
$\Psi=\id$. Our approach will be to use a holomorphically varying homotopy
between $f$ and $f_h$. 

\give Proof. The idea is to define the conjugacy $\Psi$ as
$\lim_{n\to\infty}f^{-n}\circ f_h^n$, the difficulty being to
define $f^{-n}$. We will use the notation from the proof of Theorem~4.1.

Fix $a\in J_\Pi$, $n\ge 1$ and $\tau$ with $|\tau|<2$. Write
$f_\tau=f_h+\tau(f-f_h)$.  Let $\Delta_n$ be the disk defined by 
$$
\Delta_n=f_h^{nN}W^s_0(a,f_h)=W^s_0(f_\Pi^{nN}a,f_h)\cap A_{-nN,h}.
$$
Then $f_\tau^{-nN}(\Delta_n)$ is a variety in
$f_\tau^{-nN}A_{-nN,0}$, all of whose components must meet $\Pi$
at some point in $f_\Pi^{-nN}f_\Pi^{nN}a$. Hence, we may write
$$
f_\tau^{-nN}(\Delta_n)=\bigcup_{b\in f_\Pi^{-nN}f_\Pi^{nN}a}\beta_{b,n},
$$
where $\beta_{b,n}$
is contained in a vertical-like disk in $B_b$. Thus
$f_\tau^{nN}$ maps $\beta_{b,n}$ onto $\Delta_n$ as a branched covering
of degree $\deg$, branched only at $a$.

Hence there are $\deg^n$
locally defined branches of $f_\tau^{-nN}\circ f_h^{nN}$ mapping
$W^s_0(a,f_h)-\{a\}$ into $\beta_{a,n}\subset B_a$.
These branches depend holomorphically on $\tau$.
Let $\psi_{a,\tau,n}$ be the branch obtained by analytic
continuation of $\psi_{a,0,n}=\id$. Then
$\psi_{a,\tau,n}$ is well-defined on $W^s_0(a,f_h)$, 
depends continuously on $a$ and holomorphically on $\tau$.

The mappings $\psi_{a,\tau,n}$ map 
$W^s_0(a,f_h)\times\{|\tau|<2\}$ into
$B_a$, hence they form a normal family.
We claim that in fact $\psi_{a,\tau,n}$ converges as $n\to\infty$.
To see this, we first note that
$\deg^{-n}G_h\circ f_\tau^n\to G_\tau$ uniformly on
compact subsets of $(A_{0,h}-\Pi)\times\{|\tau|<2\}$ by Lemma~1.1.
Hence any limit point $\psi_{a,\tau}$ of
$\psi_{a,\tau,n}$ must have the following properties:
\item{(i)} $\psi_{a,\tau}$ depends holomorphically on $\tau$.
\item{(ii)} $G_\tau\circ\psi_{a,\tau}=G_h$.
\item{(iii)} $\psi_{a,\tau}$ maps $W^s_0(a,f_h)$ 
  biholomorphically onto $W^s_0(a,f_\tau)$.
\item{(iv)} $\psi_{a,0}=\id$.

Now suppose that $\psi_{a,\tau}'$ and $\psi_{a,\tau}''$ are two such
limits. By~(iii) the mapping
$\nu_{a,\tau}:=(\psi_{a,\tau}'')^{-1}\circ\psi_{a,\tau}'$
is a biholomorphism of the disk $W^s_0(a,f_h)$. By~(i)
$\nu_{a,\tau}$ depends holomorphically on $\tau$ and
by~(ii) $\nu_{a,\tau}$ is a rotation for all $\tau$.
Hence $\nu_{a,\tau}$ is a constant (not depending on $\tau$)
times the identity, so by~(iv) $\nu_{a,\tau}=\id$,
i.e.\ $\psi_{a,\tau}'=\psi_{a,\tau}''$ for all $\tau$.

Thus $\psi_{a,\tau,n}$ converges to a map $\psi_{a,\tau}$ having
the properties~(i)--(iv).
Define $\Psi_\tau$ by
$\Psi_\tau=\psi_{a,\tau}$ on $W^s_0(a,f_h)$. Then
$\Psi_\tau$ is a bijection of 
$W^s(J_\Pi,f_h)\cap A_{0,h}$ onto $W^s(J_\Pi,f_\tau)\cap A_{0,\tau}$.

We claim that $\Psi_\tau$ conjugates $f_h$ to $f_\tau$.
This amounts to showing that, for all $a\in J_\Pi$,
$\psi_{f_\Pi a,\tau}\circ f_h=f_\tau\circ\psi_{a,\tau}$
on $W^s_0(a,f_h)$.
Now these two mappings are both branched coverings of 
$W^s_0(a,f_h)$ onto $W^s_0(f_\Pi a,f_\tau)\cap A_{-1,\tau}$
of degree $\deg$, branched only at $a$. Moreover, we have
$$
G_\tau\circ\psi_{f_\Pi a,\tau}\circ f_h
=d\cdot G_h
=G_\tau\circ f_\tau\circ\psi_{a,\tau}.
$$
Hence there exists a complex number $\nu_{a,\tau}$ of unit
modulus such that
$$
f_\tau\circ\psi_{a,\tau}\circ\nu_{a,\tau}=\psi_{f_\Pi a,\tau}\circ f_h.
$$
Since $\nu_{a,\tau}$ depends holomorpically on $\tau$ and
$\nu_{a,0}=1$ we see that $\nu_{a,\tau}=1$ for all $\tau$.

We complete the proof by showing that $\Psi_\tau$ is continuous
for all $\tau$. Fix $a\in J_\Pi$ and pick parametrizations
$\chi_b:\hat\C-\bar\D_{R_0}\to W^s_0(b,f_h)$ for $b\in J_\Pi$ 
close to $a$, such that
$\chi_b$ depends continuously on $b$ and $G_h\circ\chi_b=\log|\cdot|$.
It suffices to prove that $\psi_{b,\tau}\circ\chi_b$ converges to
$\psi_{a,\tau}\circ\chi_a$ as $b\to a$. Again this follows,
using the Green functions and the holomorphic dependence on $\tau$.
\qed

\give Remark. The above result is similar to Theorem~9.3 in [HP].
\bigskip

%
%

\section 5 Uniform laminar structure of $T^{\dim-1}$ near $\Pi$.

The next two sections will be devoted to the laminarity of the current
$T^{\dim-1}$ on $A$. In~\S6, we cover some rather general situations.  It is
worthwhile, however, to start with the case where
$f_\Pi$ is uniformly expanding on $J_\Pi$, and we obtain~(5.1), which is our
strongest laminarity property.

\proclaim Theorem~5.1. If $f_\Pi$ is expanding on $J_\Pi$ and 
$R_0$ is large enough, then
$$
T^{\dim-1}\contract A_0
=\int\left[W^s_0(a)\right]\,\mu_\Pi(a).\eqno(5.1)
$$

\give Proof. It follows from Proposition~2.1 that
$$
\eqalign{
  {1\over{\deg^{j(\dim-1)}}}
  \left(f^j\right)^*
  &\left(T_h^{\dim-1}\contract A_0\right)
  \contract A_0\cr
  &={1\over{\deg^{j(\dim-1)}}}
  \left(\left(f^j\right)^*\int
    \left[L_a\cap A_0\right]\mu_\Pi(a)\right)
  \contract A_0,}
$$
for all $j\ge 0$.
We will show: (1) the left hand side tends to $T^{\dim-1}\contract A_0$ as
$j\to\infty$, and (2) the right hand side tends to
$\int[W^s_0(a)]\,\mu_\Pi(a)$ as 
$j\to\infty$.

To prove (1) it suffices to show that
$\deg^{-j}G_h\circ f^j\to G$ uniformly on $A_0-\Pi$. But from
Lemma~1.1 we know that $G-G_h$ is bounded on $A_0$. Thus
$$
\eqalign{
{1\over \deg^j}G_h\circ f^j-G
&={1\over \deg^j}\left(G_h-G\right)\circ f^j\cr
&=O({1\over \deg^j}).}
$$
To show (2), we use Lemma~2.3 and calculate
$$
\eqalign{
{1\over{\deg^{(\dim-1)j}}}
&\left(\left(f^j\right)^*\int
  \left[L_a\cap A_0\right]\,\mu_\Pi(a)\right)
\contract A_0\cr
&=\int {1\over{\deg^{(\dim-1)j}}}
\left[f^{-j}\left(L_a\cap A_0\right)\cap A_0\right]
\,\mu_\Pi(a)\cr
&=\int {1\over{\deg^{(\dim-1)j}}}
\left[f^{-j}\left(L_{f_\Pi^ja}\right)\cap A_0\right]
\,\mu_\Pi(a),}
$$
where we have used the invariance of $\mu_\Pi$.
From the proof of Theorem~4.1 we know that 
$f^{-j}L_{f_\Pi^ja}\cap A_0$ is a union
of $\deg^{(\dim-1)j}$ disjoint complex disks $\gamma_j(b)$,
over $b\in f_\Pi^{-j}f_\Pi^ja$ (at least if $j$ is a multiple of $N$,
with $N$ from the same proof). Hence we get
$$
\eqalign{
\int {1\over{\deg^{(\dim-1)j}}}
\left[f^{-j}\left(L_{f_\Pi^ja}\right)\cap A_0\right]
\,\mu_\Pi(a)
&=\int {1\over{\deg^{(\dim-1)j}}}
\sum_{b\in f_\Pi^{-j}f_\Pi^ja}\left[\gamma_j(b)\right]\,\mu_\Pi(a)\cr
&=\int\left[\gamma_j(a)\right]\,\mu_\Pi(a),}
$$
since $f_\Pi^*\mu_\Pi=\deg^{\dim-1}\mu_\Pi$.
Moreover, from the same proof it follows that $\gamma_j(a)$ converges
to the local stable disk $W^s_0(a)$ in $C^1$-topology,
uniformly in $a$. Hence the last line above converges to 
$\int[W^s_0(a)]\,\mu_\Pi(a)$ as $j\to\infty$,
completing the proof.
\qed

Theorem~5.1 allows us to describe the support of 
$T^{\dim-1}\contract A$ in dynamical terms.

\proclaim Corollary~5.2. If $f_\Pi$ is expanding on $J_\Pi$, then
$\supp(T^{\dim-1})\cap A=W^s(J_\Pi)$.

\give Proof. It follows from Theorem~5.1, from the continuity of
$a\to W^s_0(a)$ and from Proposition~4.2 that
the support of $T^{\dim-1}\contract A_0$ is equal 
to $W^s(J_\Pi)\cap A_0$.
This proves the corollary, because 
the sets $\supp(T^{\dim-1})\cap A$ 
and $W^s(J_\Pi)$ are both completely invariant and any compact subset of 
either of them is mapped by some iterate of $f$ into $A_0$.
\qed

Another consequence of Theorem~5.1 is that $T^{\dim-1}$ has a 
uniform laminar structure on $A_n=f^{-n}A_0$ for every $n\ge 0$,
hence on every relatively compact subset of $A$.

\proclaim Corollary~5.3. For every $n\ge 0$ we have
$$
T^{\dim-1}\contract A_n
=\int{1\over\deg^{(\dim-1)n}}
\left[f^{-n}W^s_0(f_\Pi^na)\right]\,\mu_\Pi(a).\eqno(5.2)
$$

\give Proof. This is an easy consequence of Theorem~5.1.
Indeed,
$$
\eqalign{
T^{\dim-1}\contract A_n
&={1\over \deg^{(\dim-1)n}}
\left(f^n\right)^*\left(T^{\dim-1}\contract A_0\right)\cr
&={1\over \deg^{(\dim-1)n}}
\left(f^n\right)^*\left(\int\left[W^s_0(a)\right]\,\mu_\Pi(a)\right)\cr
&={1\over \deg^{(\dim-1)n}}
\int\left[f^{-n}\left(W^s_0(a)\right)\right]\,\mu_\Pi(a)\cr
&=\int{1\over \deg^{(\dim-1)n}}
\left[f^{-n}\left(W^s_0\left(f^n(a)\right)\right)\right]\,\mu_\Pi(a).}
$$
\qed

%
%

\section 6 Laminar Structure of $T^{\dim-1}$ on $A$.

The main goal of this section is to show that the current
$T^{\dim-1}\contract A$ has a laminar structure. We have seen in the
previous section that if $f_\Pi$ is uniformly expanding on $J_\Pi$,
then $T^{\dim-1}\contract A$ has a uniformly laminar structure in a
neighborhood of $\Pi$ with respect to the Riemann surface lamination
given by the local stable disks $W^s_0(a)$.
Also, $T^{\dim-1}$ is uniformly laminar on $A_n$ for each $n\ge 0$,
hence on each compact subset of $A$.

Here we show that there is a nonuniform laminar
structure in general. When the expansion of $f_\Pi$ is not uniform,
we still have stable manifolds by Pesin theory. Without uniformity,
however, we are not able to bound the topological type of the stable
manifolds in a neighborhood of $\Pi$. Despite this, there is a
formulation of the laminarity of $T^{\dim-1}\contract A_n$
(Theorem~6.4) in terms of currents of integration over subvarieties of
$A_n$.  This allows us to express the restriction of the critical
measure $\mu_c$ to $A_n$ as an intersection product with the critical
locus (Corollary~6.5). A more global laminar formulation for
$T^{\dim-1}\contract A$ (Theorem~6.10) is obtained by subdividing the
manifolds in the lamination into disks, which is done by cutting along
the gradient lines of $G$.

Our starting point is the result by Briend [Bri] that the Lyapunov
exponents of $f_\Pi$ with respect to $\mu_\Pi$ are strictly
positive. More precisely, for $\mu_\Pi$-almost every $a$ and all
$v\in T_a\Pi$, $v\ne 0$, we have
$$
\liminf_{j\to\infty}{1\over j}
\log|Df_\Pi^j(a)v|\ge {1\over 2}|v|\log\deg .
$$

Thus, if we set $E^u_a=T_a\Pi$, then $f$ is (nonuniformly)
expanding on the subspace $E^u_a$. When we consider the mapping $f$
at a point $a$ of $\Pi-\cC_\Pi$, there is a unique one-dimensional
subspace $E_a^s$ of $T_a\px\dim$ such that the restriction of $Df$
to $E_a^s$ is zero. We also have $Df(T_a\Pi)\subset T_{f_\Pi a}\Pi$.
Note that $\mu_\Pi(\cC_\Pi)=0$, since $\cC_\Pi$
is pluripolar in $\Pi$ and $\mu_\Pi$ has continuous local potentials.
Hence, for $\mu_\Pi$-a.e.\ $a$ the tangent space
of $\px\dim$ at $a$ is the direct sum of
two subspaces on which $Df^j$ is asymptotically expanding and
contracting, respectively. In other words,
$\mu_\Pi$ is a hyperbolic measure for $f$.

By Pesin theory there exists a local stable
manifold through almost every point of $J_\Pi$.
In general, we write
$$
W^s_{\rm loc}(a)=\{x\in\px\dim: d(f^jx,f^ja)<\delta\ \forall j\ge 0\}
$$
for small $\delta>0$.
Since $\mu_\Pi$ is a hyperbolic
measure, it is a consequence of
Pesin Theory that for $\mu_\Pi$-almost every $a\in J_\Pi$
there exists a $\delta=\delta(a)>0$ such that $W^s_{\rm loc}(a)$
is an embedded real 2-dimensional disk in $\px\dim$, tangent to
$E^s_a$ at $a$. Since $f$ is
holomorphic, $W^s_{\rm loc}(a)$ is in fact a complex
disk in $\px\dim$. We may choose $m=m(a)\ge 0$ such that
$W^s_{\rm loc}(a)$ is proper in the neighborhood $A_{-m}$ of $\Pi$.


The precise statement from Pesin Theory that we will need
is the following, which is an adaptation of Corollary~5.3 of [PS].

{\sl For every $\eta>0$ there exists a compact subset $F=F_\eta$ of
$J_\Pi$ with $\mu_\Pi(F)\ge 1-\eta$ and an integer $m=m(\eta)\ge 0$
such that the following holds:
\item{(a)} $F$ has no isolated points and does not intersect the set
  $\bigcup_{j\in\Z}f^j(\cC_\Pi)$.
\item{(b)} For each $a\in F$, the local stable manifold
  $W^s_{\rm loc}(a)$ is proper in $A_{-m}$ and
  $W^s_{-m}(a):=W^s_{\rm loc}(a)\cap A_{-m}$
  is a properly embedded disk in $A_{-m}$. The Green function
  $G$ is harmonic on $W^s_{-m}(a)-\{a\}$ and has no critical point there.
\item{(c)} The map $a\to W^s_{-m}(a)$ is continuous in the $C^1$
  topology and the set of disks $\{W^s_{-m}(a):a\in F\}$ defines a
  Riemann surface lamination in $A_{-m}$.
\item{(d)} For each $a\in F$, we let $L_a$ denote the complex line
  in $\px\dim$ defined by $a$, and we let $D_j(a)$ denote the
  component of $f^{-j}(L_{f^j_\Pi a})\cap A_{-m}$ containing $a$.
  Then $D_j(a)$ is a complex disk in $A_{-m}$ for $j\ge 0$ and
  $D_j(a)$ converges in the $C^1$ topology to $W^s_{-m}(a)$ as
  $j\to\infty$. This convergence is uniform in $a$ for $a\in F$.
\item{}}

Property~(a) is a consequence of the fact that $\mu_\Pi$ does
not give mass to pluripolar sets. For~(b), we know from Lemma~1.2
that $G$ is harmonic on the local stable manifolds.
Since the local stable manifolds vary continuously
over $F$, we see from Lemma~1.1 that $G$ has no critical points
on $W^s_{\rm loc}(a)$ if $\delta$ is small enough. Thus we may
choose $m\ge 0$ so that $W^s_{\rm loc}(a)\cap A_{-m}$ is a
properly embedded complex disk in $A_{-m}$. Properties~(c) and~(d)
follow from the construction of the local stable manifolds as
in~[PS].

We will call a set $F$ satisfying~(a)--(d) a {\it Pesin box\/}
and the associated disks $W^s_{-m}(a)$ {\it Pesin disks\/}.

In general we must let $m(\eta)\to\infty$ as $\eta\to 0$ to insure that
$W^s_{\rm loc}(a)$ is a proper disk in $A_{-m}$ for all $a\in F$. 
Figure~4 illustrates this phenomenon.

\centerline{\psfig{figure=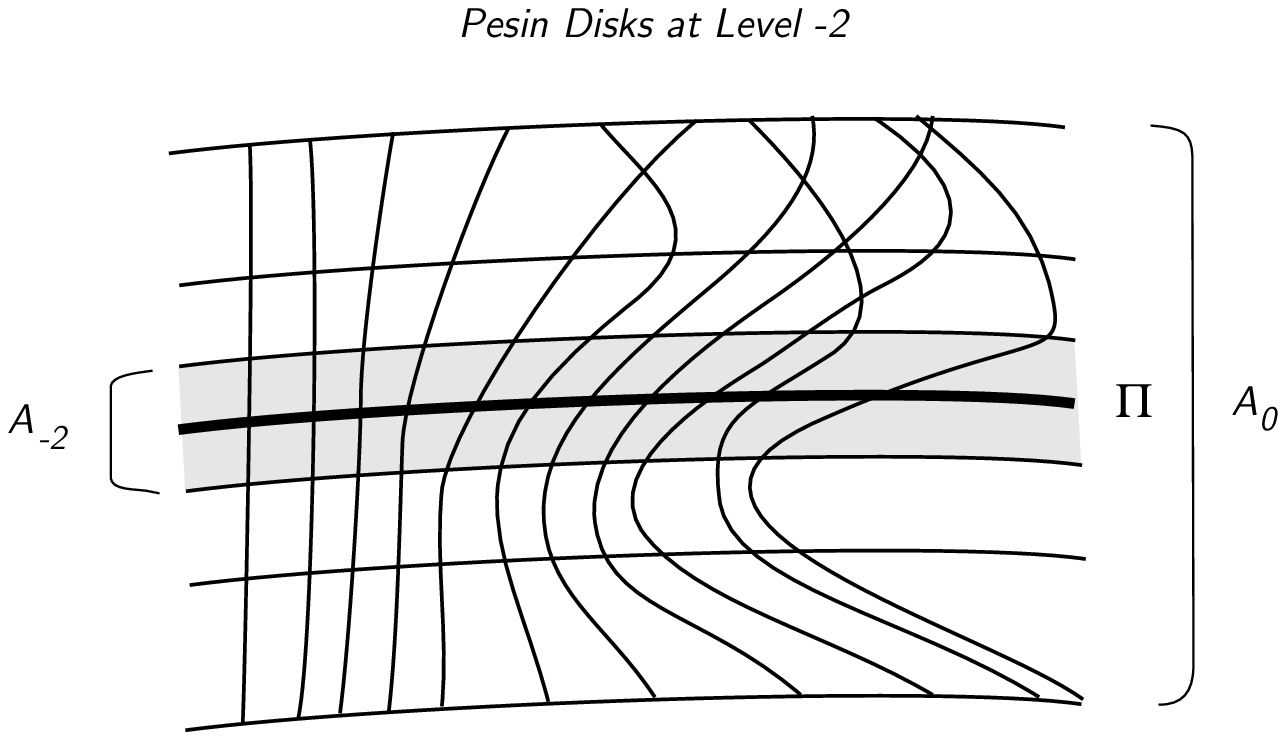,width=0.7\hsize}}
\centerline{Figure 4.}
\medskip

In~\S5 we saw that the (uniform) hyperbolicity of $f$ on $J_\Pi$
implies that $T^{\dim-1}$ has a uniform laminar structure in a
neighborhood of $\Pi$. The following lemma is a corresponding
result in the nonuniformly hyperbolic case. By restricting to a small
neighborhood $A_{-m}$ of $\Pi$ we get many Pesin disks
and these account for a large part of $T^{\dim-1}$ in this neighborhood.

\proclaim Lemma~6.1. Given $\eta>0$ there exists $m\ge 0$
and a Pesin box $F=F_\eta$ having properties~(a)--(d) above
such that
$$
T^{\dim-1}\contract A_{-m}
=\int_{a\in F}[W^s_{-m}(a)]\,\mu_\Pi(a)+S\contract A_{-m},\eqno(6.1)
$$
where $S$ is a positive closed current on
$\px\dim$ with $\Vert S\Vert\le\eta$.

\give Proof. We let $F$ satisfy properties (a)--(d) above.
From Proposition~2.1 we have
$$
\eqalign{
T_h^{\dim-1}\contract A_{-m}
&=\int_{J_\Pi}[L_a\cap A_{-m}]\,\mu_\Pi(a)\cr 
&=\int_{J_\Pi}[L_{f^j a}\cap A_{-m}]\,\mu_\Pi(a)
= \int_F +\int_{{}^cF},\cr}
$$
where we have used the invariance of $\mu_\Pi$. Now we apply $f^{*j}$
to this equation, divide by $\deg^{(\dim-1)j}$,
restrict to $A_{-m}$ and let $j\to\infty$. 
As in Theorem~5.1 we see that the left hand side then 
converges to $T^{\dim-1}\contract A_{-m}$ and, using (b),
${1\over\deg^{(\dim-1)j}}f^{*j}(\int_F)\contract A_{-m}$ converges to
$\int_F[W^s_{-m}(a)]\,\mu_\Pi(a)$.

Finally, we consider
${1\over\deg^{(\dim-1)j}}f^{*j}(\int_{{}^cF})\contract A_{-m}$.
The mass norm of the current of integration over 
a curve of degree $r$ in $\px\dim$ is $r$ (up to a constant only
depending on the volume of $\px\dim$). In particular,
$\Vert L_a\Vert=1$. It follows that the currents
${1\over\deg^{(\dim-1)j}}f^{*j}(\int_{{}^cF}[L_{f^ja}])$
have mass norms that are bounded by $\eta$.
Passing to a subsequence, we obtain a limit $S$ which has the desired
properties in~(6.1). \qed

\proclaim Corollary~6.2. Given the assumptions in Lemma~6.1 and $n\ge 0$,
we have
$$
T^{\dim-1}\contract A_n
=\int_F{1 \over \deg^{(\dim-1)(n+m)}}[f^{-(n+m)}W^s_{-m}(a)]\,\mu_\Pi(a)
+S\contract A_n,\eqno(6.2)
$$
where $S$ is a positive, closed current on $\px\dim$ with mass norm
bounded by $\eta$.

\give Proof. We pull~(6.1) back by $f^{n+m}$ and divide by
$\deg^{(\dim-1)(n+m)}$. Note that the proof of Lemma~6.1 shows that
the current $S$ in~(6.1) satisfies
$\Vert (f^{n+m})^*S\Vert\le\deg^{({\dim-1})(n+m)}$. \qed

Corollary~6.2 shows that we can approximate $T^{\dim-1}$ on $A_n$
by laminar currents over closed varieties in $A_n$. We want to pass
to the limit and obtain an exact formula, just as in~(5.2).

Let us fix a Pesin box $E_0\subset J_\Pi$ satisfying (a)--(d)
above with some $\eta>0$.
We may choose $R_0$ so that $m=m(\eta)=0$.
By Poincar\'e recurrence, the set
$$
E=\{a\in J_\Pi: f^n_\Pi a\in E_0
{\rm \ for\ infinitely\ many\ }n\ge0\}\eqno(6.3)
$$
has full measure. It is clearly invariant under $f_\Pi$.
Recall the notion of the global stable manifold $W^s(a)$ of a point
$a\in J_\Pi$. In general, $W^s(a)$ will be a very complicated object,
but we do have the following.

\proclaim Lemma~6.3. If $a\in E$, and $n\ge 0$,
then $W^s(a)\cap A_n$ is
a disjoint countable union of connected varieties in
$A_n$, each of which intersects $\Pi$ transversely at
finitely many points.

\give Proof. Let $n\le n_1<n_2<\dots$ be return times for $a$ to
$E_0$, i.e.\ $f_\Pi^{n_j}a\in E_0$.
Let $Z_a^{(j)}=f^{-n_j}W^s_0(f_\Pi^{n_j}a)\cap A_n$ for $j\ge 1$.
Then $Z_a^{(j)}$ is a variety in $A_n$ which intersects
$\Pi$ transversely at finitely many points because of property~(a).
Further, $f^{n_{j+1}-n_j}$ maps $W^s_0(f_\Pi^{n_j}a)$ into
$W^s_0(f_\Pi^{n_{j+1}}a)$, so $Z_a^{(j)}\subset Z_a^{(j+1)}$ and
by definition of $W^s(a)$ we have
$W^s(a)\cap A_n=\bigcup_{j\ge 1}Z_a^{(j)}$. We complete the 
proof by showing that
$Z_a^{(j)}$ is open in $Z_a^{(j+1)}$. Take
$x\in Z_a^{(j)}$ and let $y\in Z_a^{(j+1)}$ be close to $x$.
We have to show that $f^{n_j}y\in W^s_0(f_\Pi^{n_j}a)$.
If $y$ is close to $x$, then the orbit $(f^lf^{n_j}y)_{l\ge 0}$
stays close to $(f^lf^{n_j}a)_{l\ge 0}$ for $0\le l\le n_{j+1}-n_j$.
But $f^{n_{j+1}}y\in W^s_0(f_\Pi^{n_{j+1}}a)$, so the same is true
for $l\ge n_{j+1}-n_j$. Hence $f^{n_j}y\in W^s_0(f_\Pi^{n_j}a)$
and we are done.
\qed

\proclaim Theorem~6.4. Given any regular polynomial endomorphism
of $\cx\dim$ there is a set $E\subset J_\Pi$
with $\mu_\Pi(E)=1$ such that the following holds. If $a\in E$,
$n\ge 0$ and $Z_{a,n}$ is the connected component of $W^s(a)\cap A_n$ 
containing $a$, then $Z_{a,n}$ is a one-dimensional
subvariety of $A_n$ intersecting $\Pi$ in $N_n(a)<\infty$ points, and
$$
T^{\dim-1}\contract A_n
=\int_E{1\over N_n(a)}[Z_{a,n}]\,\mu_\Pi(a).\eqno(6.4)
$$
Further, the restriction of $G$ to $Z_{a,n}$ is harmonic outside
the singular locus of $Z_{a,n}$ for all $a\in E$.

\give Proof. The first part follows from Lemma~6.3 and the harmonicity
of $G$ from Lemma~1.2. It remains to show~(6.4).
Given $\eta>0$, let $F=F_\eta$ be a Pesin box satisfying
(a)--(d) above with $m=m(\eta)\ge 0$.
If $a\in f_\Pi^{-(n+m)}F_\eta\cap E$, then
$f^{-(n+m)}(W^s_{-m}(f_\Pi^{n+m}a))$ is  a subvariety of $A_n$, contained
in $W^s(a)$. By Lemma~6.3 we may therefore write
$$
f^{-(n+m)}(W^s_{-m}(f_\Pi^{n+m}a))
=\bigcup_{b\in f_\Pi^{-(n+m)}f_\Pi^{n+m}a}Z_{b,n},
$$
and, by definition of $N_n(b)$,
$$
\sum_{b\in f_\Pi^{-(n+m)}f_\Pi^{n+m}a}{1\over\deg^{(\dim-1)(n+m)}}
[f^{-(n+m)}(W^s_{-m}(f_\Pi^{n+m}b))]
=\sum_{b\in f_\Pi^{-(n+m)}f_\Pi^{n+m}a}{1\over N_n(b)}[Z_{b,n}].
$$
This, Corollary~6.2, and the invariance of $\mu_\Pi$ yield
$$
T^{\dim-1}\contract A_n=
\int_{f_\Pi^{-(n+m)}F\cap E}{1\over N_n(a)}[Z_{a,n}]\,\mu_\Pi(a)
+S\contract A_n.
$$
Theorem~6.4 follows by letting $\eta\to 0$.
\qed

\proclaim Corollary~6.5. With the assumptions and notation of 
Theorem~6.4, we have
$$
\mu_c\contract A_n
=\int_E{1\over N_n(a)}[\cC\cap Z_{a,n}]\,\mu_\Pi(a).\eqno(6.5)
$$
Further, $\mu_c\contract A_n=0$ holds if and only if
$\cC\cap Z_{a,n}=\emptyset$ for almost every $a$ and
$\mu_c\contract A=0$ if and only if
$\cC\cap W^s(a)=\emptyset$ for almost every $a$.

\give Proof. The formula~(6.5) follows from Lemma~2.2 and Corollary~6.4. 
It follows directly from~(6.5) that $\mu_c\contract A_n=0$ if and
only if $\cC\cap Z_{a,n}=\emptyset$ for almost every $a$. By Lemma~6.3
this happens if and only if $\cC\cap W^s(a)\cap A_n=\emptyset$.
Letting $n\to\infty$ we get $\mu_c\contract A=0$ if and only if
$\cC\cap W^s(a)=\emptyset$ for almost every $a$.
\qed

Formula~(6.4) exhibits $T^{\dim-1}\contract A_n$ as a laminar current
using currents of integration over (closed) subvarieties $Z_{a,n}$ in
$A_n$. These varieties $Z_{a,n}$ are subsets of the global stable manifolds
$W^s(a)$. We would like to have a laminar structure for $T^{\dim-1}$
in the larger set $A$. We could try to do this by attempting to 
extend the varieties $Z_{a,n}$ analytically 
to subvarieties of $A$. However,
in~\S4 we gave an example of a map $f$
(with $f_\Pi$ uniformly expanding on $J_\Pi$)
where the global stable manifolds $W^s(a)$ are connected and have locally
infinite area in $A$. The analytic continuations of $Z_{a,n}$ would be
$W^s(a)$ in this case, hence would not define currents of integration.

Nevertheless, we will show that $T^{\dim-1}\contract A$ does have a laminar
structure. We will accomplish this by dividing the global
stable manifolds $W^s(a)$ into disks $W_a$. These disks are not
in general closed in $A$. The construction of $W_a$ consists of cutting
$W^s(a)$ along gradient lines of $G$. These gradient lines will
be used in~\S7 to exhibit $\mu$ as a quotient of the product of $\mu_\Pi$
and Lebesgue measure on the circle.

Let $E_0$ and $E$ be the sets defined above.
Recall that $G|_{W^s_0(a)}$ has no critical
points for $a\in E_0$.
Let $\cC_n=\cC\cup f^{-1}\cC\cdots\cup f^{-(n-1)}\cC$ be the
critical set of $f^n$.
Let us fix $a\in E$ and $n\ge 0$ 
and recall the definition of $Z_{a,n}$ above.
Suppose that $Z_{a,n}\subset\cC_n$. Since
the restriction of $G$ to $Z_a$ is harmonic outside the singular locus,
$G|_{Z_a}$ cannot be
bounded above, so $Z_{a,n}\cap\Pi$ is nonempty. It follows that
$Z_{a,n}\cap\cC_n\cap\Pi\ne\emptyset$, which contradicts property~(a).
Thus $\cC_n\cap Z_{a,n}$ is a discrete set. In fact it is finite, because
$\cC_n\cap Z_{a,n+1}$ is also discrete. Now $f^n$ is a
local biholomorphism on $\cx\dim-\cC_n$, so $f^{-n}$ serves to
transfer certain properties from the Pesin disks of $E_0$.
Specifically, let $n$ be such that $f_\Pi^na\in E_0$. Such an $n$
will be called a {\it return time}. If $n$ is a return time, 
then $Z_{a,n}-\cC_n$ is a manifold,
and the restriction of $G$ to this manifold is a harmonic
function without critical points.

We wish to define gradient lines of $G$. There is a unique tangent
line to the level sets of $G|_{Z_{a,n}}$ at points off
of $\cC_n$. By the conformal structure, we may define the gradient
vector of $G$ to be the tangent to $Z_{a,n}-\cC_n$
which is orthogonal to the level line, and which points in the
direction of increasing $G$. A {\it gradient line\/} is then an integral
curve of the gradient vector field. Equivalently, a gradient line
is, locally, a level set for a harmonic conjugate to $G$.
We note that these so-called
gradient lines are an artifact of the conformal structure. We could
have as well defined $\tau$-gradient lines, which make an angle of
$\tau$ with respect to the gradient inside the tangent space to
$Z_{a,n}-\cC_n$.

We say that a gradient line $\gamma$ is {\it complete\/} if it
is a complete orbit of the gradient vector field, and
if $\sup_{\gamma}G=+\infty$.
We let $S_{a,n}$ denote the set $Z_{a,n}\cap \cC_n$,
together with points of $Z_{a,n}-\cC_n$ which are not
contained in complete gradient lines.

\proclaim Lemma~6.6. The set $S_{a,n}-\cC_n$ consists of a finite
number of (incomplete) gradient lines.

\give Proof. If $\gamma$ is an incomplete gradient line, it follows
that the closure of $\gamma$ must contain a singular point of the
gradient field. Let us fix a point $p\in Z_{a,n}\cap \cC_n$.
It suffices to show that there are only finitely many gradient lines
whose closures contain $p$. For this, we let $h:\Delta\to
Z_{a,n}$ be a holomorphic mapping with $h(0)=p$, and
which is a homeomorphism with its image. It follows that
$G\circ h(\zeta)=\Re (a_l\zeta^l)+\cdots$\ . By further composing with a
conformal map fixing the origin, we may achieve that
$G\circ h(\zeta)=\Re(\zeta^l)$. Since the gradient lines are conformally
invariant, it follows that they are mapped into the gradient lines of
the function $\Re(\zeta^l)$, which are a finite family of straight
lines through the origin.
\qed

It follows that $S_{a,n}$ is a closed subset of $Z_{a,n}$,
and thus $Z_{a,n}-S_{a,n}$ is a manifold and an
open subset of the subvariety $Z_{a,n}\subset A_n$. An
alternative definition of $Z_{a,n}-S_{a,n}$ is that it
is the largest open subset of $Z_{a,n}-\cC_n$ which is
invariant under the (positive) gradient flow. For every point
$x\in Z_{a,n}-S_{a,n}$, the gradient line starting at $x$
approaches a unique point $a'\in\Pi$.

Let $W_{a',n}$ denote the connected component of
$Z_{a,n}-S_{a,n}$ containing $a'$. Then the gradient lines of
all points of $W_{a',n}$ must approach $a'$. Thus $W_{a',n}\cap
\Pi=\{a'\}$. From the paragraph above, we see that $W_{a,n}$ consists
of all of the complete
gradient lines in $Z_{a,n}-S_{a,n}$ that emanate from $a$. It is evident,
then, that the gradient lines serve as a sort of exponential map from
the tangent space of $W_{a,n}$ at $a$ to $W_{a,n}$.

If $m$ is the number given in conditions (a)--(d) from Pesin Theory given above, then $Z_{a,-m}$
contains no critical points, and so for $a\in E$,
$$W^s_{-m}(a)=W_a\cap A_{-m}=Z_{a,-m}.$$

\proclaim Lemma~6.7. If $n$ is a return time for $a$, then
$Z_{a,n}-S_{a,n}$ consists of $N_n(a)$
components, each of which is of the form $W_{a',n}$ for some
$a'\in f_\Pi^{-n}f_\Pi^na$. The manifolds $W_{a',n}$ are simply connected.
If $n'>n$ is another return time, then $W_{a,n}\subset W_{a,n'}$, and
$$
W_{a,n'}-W_{a,n}\subset A_{n'}-A_n.
$$

\give Proof. The first statement is a consequence from the discussion
above. For the simple connectivity of $W_{a,n}$, we let $\sigma$
denote a simple, closed curve in $W_{a,n}$. Moving $\sigma$ along the
gradient flow of $G$ is a homotopy, which takes $\sigma$ to a
neighborhood of $a$, where it is contractible. The last statements
follow because $W_{a,n}$ is the union of complete gradient lines
emanating from $a$. \qed

We define $W_a=\bigcup_n W_{a,n}$, where the union is taken over any
sequence of return times $n\to\infty$.

\proclaim Corollary~6.8.
$W_a$ is a simply connected Riemann surface, and the topology of $W_a$
as a manifold coincides with the topology of $W_a$ as a topological
subspace of $A$. If $n$ is a return time for $a$, then
$W_{a,n}=W_a\cap A_n$.

\proclaim Lemma~6.9. For any $n\ge 0$ and any $a\in E$ we have
$$
[Z_{a,n}]=\sum_{b\in Z_{a,n}\cap\Pi}[W_b\cap A_n].\eqno(6.6)
$$

\give Proof. First suppose $n$ is a return time for $a$. Then we know
from Lemma~6.7 that
$$
Z_{a,n}=\bigcup_{b\in Z_{a,n}\cap\Pi}W_{b,n}\cup S_{a,n}.
$$
Thus~(6.6) follows from the fact that $S_{a,n}$
has zero area (Hausdorff 2-dimensional measure) in $\px\dim$. This in
turn is a consequence of the fact that each incomplete gradient line
is the countable union of real analytic arcs of finite length (and
thus zero area).

If $n$ is not a return time for $a$, then let $n'>n$ be a return time.
Thus 
$$
Z_{a,n'}=\bigcup_{b\in Z_{a,n'}\cap\Pi}W_{b,n'}\cup S_{a,n'},
$$
and we get~(6.6) by intersecting this formula with
$Z_{a,n}$ and using the fact that $S_{a,n'}$ has zero area.
\qed

From~(6.4),~(6.6), and the definition of $N_n(a)$ we get
$$
T^{\dim-1}\contract A_n=\int_E [W_a\cap A_n]\,\mu_\Pi(a).\eqno(6.7)
$$
If we let $n\to\infty$ in~(6.7), then we
have a sequence which is eventually stationary in the sense that it is
constant on the open set $A_j$ if $n\ge j$. This yields the following.

\proclaim Theorem~6.10. If $E$ is the set defined by~(6.3), then for
$a\in E$, there exists a complex disk $W_a$ in $A$ such that
$$
T^{\dim-1}\contract A=\int_E[W_a]\,\mu_\Pi(a).
$$
In particular, the disk $W_a$ has finite area for $\mu_\Pi$ almost every
$a$.

%
%

\section 7  External Rays, Global Model

In this Section, we do two things. First we define the set $\cE$ of
external rays and show that $\cE$ carries a natural measure
$\nu$. Further, there is an endpoint mapping $e:\cE\to J$
defined almost everywhere, and this mapping
satisfies $e_*\nu=\mu$ (Theorem~7.3). 
Second, we prove Theorem~7.4, which is 
a more global version of Theorem~4.3: the conjugacy is given between
$f|\bigcup_{a\in J_\Pi}W_a$ and the the restriction of $f_h$ to a union of
hedgehogs.  

Let $a$ be a point of $E$, as defined in~(6.3). The set
$W_a$, defined in~\S6, is a simply connected Riemann surface, and
$G|_{W_a}$ is harmonic.  Let $G^*_a$ be a harmonic conjugate for
$G|_{W_a}$, well-defined modulo $2\pi i$ and is unique up to an additive
constant. Since $G|_{W_a}$ has a logarithmic pole at $a$, the function 
$$
\varphi_a:=e^{G+iG^*_a}:W_a\to\C
$$
is analytic on $W_a$, and $\varphi_a$ is locally injective on $W_a$ near
$a$. Note that $\log|\varphi_a|=G_a$, and this condition determines
$\varphi_a$ uniquely up to multiplication with a constant of unit modulus.
Since $G^*_a$ is constant on the gradient lines of $G$, these gradient
lines are taken to radial lines. By the construction of $W_a$, there are at
most finitely many gradient lines $\gamma$ with the property that
$\inf_{\gamma}G\ge\deg^{-n}R_0>0$ for any $n\ge 0$. Thus the range
$H_a:=\varphi_a(W_a)$ is a hedgehog domain of the form
$$
H_a=\hat\C-\left(\bar\D\cup\bigcup_{j=1}^N
(e^{i\theta_j},e^{r_j+i\theta_j}]\right)
$$
for some $0\le N\le\infty$, and $r_j>0$ is a sequence of
points with $r\to 0$. The case
$N=0$ is interpreted as 
$H_a=\hat\C-\bar\D$. Since $W_a$ is invariant under the
gradient
flow, it follows that $\varphi_a$ is injective. We let
$$
\psi_a:=\varphi_a^{-1}:H_a\to W_a\eqno(7.1)
$$
denote the inverse.
Thus the gradient lines in $W_a$ correspond to the images under $\psi_a$
of
the rays in
$H_a$.

Let $\cE_a$ denote the set of all gradient lines in $W_a$ for $a\in E$,
and let $\cE$ be the union of all $\cE_a$.
For $a\in E$, the gradient lines are naturally parametrized
by a choice of argument $\theta$. The function $\psi_a$ represents one
assignment of argument. Although $\theta$ is not uniquely defined, the
induced measure ${d\theta\over2\pi}$ is on $\cE_a$. Thus the measure
$\nu=\mu_\Pi\otimes{d\theta\over2\pi}$ is well defined on $\cE$.
Note that $f$ maps gradient lines to gradient lines. Thus $f$
induces a measurable mapping $\sigma:\cE\to\cE$.

We would like to assign an
endpoint to $\nu$ almost every ray. For $\gamma\in\cE$ and $r>0$ we set
$e_r(\gamma)=\gamma\cap\{G=r\}.$
For each $a\in E$, $e_r$ is defined for all but possibly finitely many
rays lying in $W_a$ and we will write $e_{a,r}$ for the
restriction of $e_r$ to the rays in $W_a$.
The mapping $\psi_a$ represents $e_r$ in the sense
that if $\gamma_\theta$ is the ray in $W_a$ corresponding to argument
$\theta$, then
$$
e_r(\gamma_\theta)=e_{a,r}(\gamma_\theta)=\psi_a(e^{r+i\theta}).
$$

\proclaim Lemma~7.1. If $a\in E$, and if $W_a$ has finite area as a
subset
of $\px\dim$, then
$$
\lim_{r\to0^+}e_r(\gamma_\theta)=\lim_{r\to0^+}\psi_a(e^{r+i\theta})
$$
exists for almost every $\theta$.

\give Proof. We work in affine coordinates in $\cx\dim\subset\px\dim$.
Let $\tilde a\in\cx\dim$ denote a point 
with $|\tilde a|=1$ such that
$\pi(\tilde a)=a\in\Pi$.
Thus we may write $\psi_a(\zeta)=\zeta^{-1}\tilde a +h_a(\zeta)$,
where $h_a$ is analytic on $H_a$.
Away from the hyperplane at infinity, the Euclidean
metric on $\cx\dim$ is equivalent to the
Fubini-Study metric on $\px\dim$. The condition that $W_a$ has
finite area in $\px\dim$ is equivalent to
$\int_{H_a}|\nabla h_a|^2<\infty$. It follows
that
$$
\int_0^1|\nabla h_a(re^{i\theta})|^2rdr = \int_0^1\left|{\partial
h_a(re^{i\theta})\over\partial r}\right|^2rdr <\infty
$$
for almost every $\theta$. Thus radial limits exist for these values of
$\theta$. \qed

It follows that there is a measurable mapping $e:\cE\to\partial K$
such that
$$
e(\gamma)=\lim_{r\to0^+}e_r(\gamma)
$$
for $\nu$ a.e.\ $\gamma$. Our next
step will be to show that the mapping $e$ pushes $\nu$ forward to $\mu$. 
We will write $G_r:=\max(G,r)$.

\proclaim Lemma~7.2. For $a\in E$,
$$
(e_{a,r})_*{d\theta }=dd^c_{W_a}G_r|_{W_a}.
$$

\give Proof.  Let us first note that the statment of the Lemma is conformally invariant.  Under a
conformal transformation, the Green function transforms by composition; the gradient lines and
level sets are preserved, and so the map $e_{a,r}$ transforms by composition.  Similarly, the
operator $dd^c$ on the right hand side is invariant, so the right hand side transforms correctly. 
Now we transform under the map
$\psi_a$.  Since $\psi_a$ is nonsingular at infinity, the measure $d\theta$ is preserved.  The
image of $W_a$ is  $H_a$, and  $G$ is taken to $\log|\zeta|$. The mapping $e_{a,r}$ then
takes the angle $\theta$ to the point $re^{i\theta}$ and
the Lemma is reduced to an elementary calculation involving $\log|z|$.
\qed

\proclaim Theorem~7.3. The endpoint mapping $e$ satisfies
$e_*\nu=\mu$. Further, $f\circ e=e\circ\sigma$, where
$\sigma:\cE\to\cE$ is the map induced by $f$.

\give Proof. Integrating the formula of Lemma~7.2 with respect to
$\mu_\Pi$, and dividing by $2\pi$, we obtain
$$
(e_{r})_*\nu=\int
{1\over2\pi}dd^c_{W_a}G_r|_{W_a}\,\mu_\Pi(a).\eqno(7.2)
$$
Let us choose $n$ sufficiently large that $\{G\ge r\}\subset A_n$.
By Theorem~6.4 we have
$$
T^{\dim-1}\contract A_n=\int_E{1\over N_n(a)}[Z_{a,n}]\,\mu_\Pi(a),
$$
and by Lemma~2.2 with $X={1\over2\pi}dd^cG_r$, we have
$$
{1\over 2\pi}dd^cG_r\wedge T^{\dim-1}
=\int_E{1\over2\pi N_n(a)}dd^cG_r\wedge[Z_{a,n}]\,\mu_\Pi(a).\eqno(7.3)
$$
Here we have dropped $\contract A_n$ because the support of $dd^cG_r$ is
contained in $A_n$. By Lemma~6.9, we have 
$$
[Z_{a,n}]=\sum_{b\in Z_{a,n}\cap\Pi}[W_b\cap A_n],
$$
where the sum has $N_n(a)$ terms.
Since $G_r$ is continuous, the wedge product of $dd^cG_r$ with the
current of integration over $Z=Z_{a,n}$ is equal to the slice
measure $dd^c_ZG_r|_Z$. Since
$G_r$ is bounded, $dd^cG_r|_Z$ can put no mass on a polar set. In
particular there can be no mass on the intersection of a 
gradient line of $G$ with $\{G=r\}$, which is an isolated point.
Thus we have 
$$
\eqalign{
dd^cG_r\wedge[Z_{a,n}]
&=\sum_{b\in Z_{a,n}\cap\Pi}dd^cG_r\wedge[W_b\cap A_n]\cr
&=\sum_{b\in Z_{a,n}\cap\Pi}dd^cG_r\wedge[W_b].\cr}
$$
Again the $\contract A_n$ is dropped, since $dd^cG_r$ is supported in
$A_n$. Thus~(7.3) becomes
$$
{1\over 2\pi}dd^cG_r\wedge T^{\dim-1}
=\int_E{1\over2\pi}dd^cG_r\wedge[W_a]\,\mu_\Pi(a).
$$
By the continuity of $G_r$, we have
$dd^cG_r\wedge[W_a]=dd^c_{W_a}G_r|_{W_a}$, so that by~(7.2) we have
$$
{1\over 2\pi}dd^cG_r\wedge T^{\dim-1}=(e_r)_*\nu.
$$
Finally, as we let $r\to0$, the left hand side converges to
$\mu={1\over2\pi}dd^cG\wedge T^{\dim-1}$,
and the right hand side converges to $e_*\nu$, which shows that
$e_*\nu=\mu$.

It is evident from the definition of $\sigma$ that
$f\circ e_{\deg r}=e_r\circ\sigma$. Thus $f\circ e=e\circ\sigma$.
\qed

\give Example. Let $f(z_1,z_2)=(z_1^2+c_1,z_2^2+c_2)$.
Since $f$ is a product, $J=J_{c_1}\times J_{c_2}$, where
$J_c$ is the one-dimensional Julia set of
$p_c(z)=z^2+c$. In the coordinate $\zeta=w/z$, $f_\Pi(\zeta)=\zeta^2$,
so $J_\Pi=\{|\zeta|=1\}$ and $f_\Pi$ is uniformly expanding on $J_\Pi$.
Using the local canonical model (Theorem~4.3) and Proposition~A.2, we
see that $\cE$ is homeomorphic to $J_\Pi\times S^1=S^1\times S^1$.
In the example at hand, the homeomorphism is given as follows. 
For $\zeta=e^{2\pi i\phi}\in J_\Pi$ and
$e^{2\pi i\theta}\in S^1$, we associate the external ray
$\gamma=\gamma(\phi,\theta)\in\cE$ with  the
properties: the projection of $\gamma$ to the
$z_1$ axis is the external ray $\gamma(\theta)$ for $p_{c_1}$,
which is asymptotic to $re^{2\pi i\theta}$ at infinity,
and the projection of $\gamma$ to the
$z_2$ axis is the external ray $\gamma(\phi+\theta)$ for $p_{c_2}$.
If the endpoint map $e$ is defined,
then we may write it as 
$e(\phi,\theta)=(e_{c_1}(\phi),e_{c_2}(\phi+\theta))$.
Any identifications in the quotient by $e$ (i.e.~pairs
$(\phi',\theta')$, $(\phi'',\theta'')$ with
$e(\phi',\theta')=e(\phi'',\theta'')$)
are thus determined by the identifications in the quotients by
$e_{c_1}$ and $e_{c_2}$.  

\medskip
We conclude this section by proving a global version of the conjugacy
given by Theorem~4.3.
First note that if $f_\Pi$ is uniformly expanding on $J_\Pi$, then we
can assume that $E_0=J_\Pi$ and that
$W^s(J_\Pi)\cap A_0\cap\cC_\infty=\emptyset$,
where $\cC_\infty=\bigcup_{j\ge 0}f^{-j}\cC$. Thus
$W^s(J_\Pi)\cap\cC_\infty$ is closed and nowhere dense in $W^s(J_\Pi)$,
and $W^s(J_\Pi)-\cC_\infty$ is a Riemann surface lamination on
which the gradient flow induced by $G$ is defined.
For $x\in W^s(J_\Pi)-\cC_\infty$,
let $\gamma(x)$ be the gradient line containing $x$. Define a numbers
$r(x)\in(0,\infty]$, $s(x)\in[0,\infty)$ for $x\in W^s(J_\Pi)$ be
declaring $r(x)=s(x)=G(x)$ if $x\in\cC_\infty$ and
$r(x)=\sup\{G(y): y\in\gamma(x)\}$,
$s(x)=\inf\{G(y): y\in\gamma(x)\}$ for $x\in W^s(J_\Pi)-\cC_\infty$.
Let $S:=\{x\in W^s(J_\Pi): r(x)<\infty\}$. Then $S$ is the union
of incomplete gradient lines.

\proclaim Theorem~7.4. If $f_\Pi$ is uniformly expanding on $J_\Pi$,
then the following hold:
\item{(1)} The set $S$ of points contained in incomplete gradient
  lines is closed and nowhere dense in $W^s(J_\Pi)$.
  Further, $W^s(J_\Pi)-S$ is a Riemann surface lamination in $A-S$.
  The leaves of this lamination are exactly the disks $W_a$ and these
  are properly embedded in $A-S$.
\item{(2)} There is a closed and nowhere dense subset $S_h$ of
  $W^s(J_\Pi,f_h)$ such that $\Psi$, defined in Theorem~4.3 extends
  to a homeomorphism
  $$
  \Psi:W^s(J_\Pi,f_h)-S_h\to W^s(J_\Pi,f)-S,
  $$
  conjugating $f_h$ to $f$.
  The set $S_h$ is a union of real rays through the origin in $\cx\dim$
  and if $D_a$ is the disk in $W^s(J_\Pi,f_h)$ passing through
  $a\in J_\Pi$, then $\Psi$ maps $D_a-S_h$ biholomorphically onto $W_a$.
\item{(3)} We may identify the set $\cE$ of external rays with the
  boundary of $W^s(J_\Pi,f_h)$, i.e.\ the union of the circles
  $\partial D_a$, $a\in J_\Pi$. With this identification, the maps
  $e_r$ are defined by $e_r(\gamma)=\Psi(e^r\cdot\gamma)$ and
  $e(\gamma)=\lim_{r\to 0}\Psi(e^r\gamma)$.
\item{(4)} We have $\mu_c(A)=0$ if and only if
  $W^s(J_\Pi)\cap\cC=\emptyset$. In this case, $S=S_h=\emptyset$
  and hence there is a conjugacy
  $\Psi:W^s(J_\Pi,f_h)\to W^s(J_\Pi,f)$ conjugating $f_h$ to $f$.

\give Proof.
 {(1)} We first prove that $S$ is closed. Note that
    $S\cap A_0=\emptyset$. Take
    $x\in W^s(J_\Pi)-S$. By definition there is a gradient line
    $\gamma=\gamma(x)$ containing $x$ and ending at a point
    $a\in J_\Pi$. We have to
    show that the same is true for all points in a neighborhood of $x$
    in $W^s(J_\Pi)$. We may assume that $x\notin A_0$, because
    otherwise there is nothing to prove. Pick $t>0$ and $n\ge 0$ with 
    $s(x)<t<G(x)$ and
    $R_0<\deg^nt<R_0\deg$. Thus
    $f^n\gamma$ is a complete gradient line in $A_{-1}$.
    In fact, all gradient lines are complete in $A_0$,
    so there is a complete
    gradient line $\gamma'$ in $A_0$ containing $f^n\gamma$.
    We have $\gamma\cap C_{-n}=\emptyset$. Thus there exists
    a branch $g$ of $f^{-n}$, defined in a neighborhood of
    $\gamma'\cap\{G>\deg^nt\}$, such that
    $g\circ f^n=\id$ on $\gamma\cap\{G>t\}$.
    The map $\Psi$ is defined on $W^s(J_\Pi,f_h)$ and
    $\gamma'':=\Psi^{-1}(\gamma')$ is a gradient line
    (i.e.\ a real line segment) in $A_{0,h}$.
    Let $U$ be a simply connected neighborhood of
    $\gamma''\cap\{\deg^nrt<G_h<\infty\}$ in
    $W^s(J_\Pi,f_h)\cap\{\deg^nt<G_h<\infty\}$
    consisting of gradient lines for $G_h$. Then $\Psi(U)$ is a simply
    connected neighborhood of $\gamma'\cap\{\deg^nt<G<\infty\}$
    consisting of gradient lines for $G$. We may assume that
    $g$ is defined in $\Psi(U)$. Thus $g\circ\Psi(U)$ is an
    open set in $W^s(J_\Pi)\cap\{t<G<\infty\}$ consisting of
    gradient lines in $\{t<G<\infty\}$. Since $x\in g\circ\Psi(U)$
    we have proved that $S$ is closed in $W^s(J_\Pi)$.
 
Since $S$ is closed, it is also nowhere dense if it contains no relative interior of
$W^s(J_\Pi)$.  But if $S$ were to contain a relative interior point of $W^s(J_\Pi)$, then it
contains relative interior of a global stable manifold
$W^s(a)$.  However, this is not possible, since $S\cap W^s(a)$
consists of at most a countable number of curves.

{(2)} Let $D_a$ be the disk in $W^s(J_\Pi,f_h)$ associated 
    with $a\in J_\Pi$. The discussion in the beginning of~\S7
    shows that there is a uniquely defined
    closed subset $S_{a,h}$ of $D_a$ such that
    $G_h>R_0$ on $S_{a,h}$ and
    such that $\Psi|_{D_a\cap\{G_h>R_0\}}$ extends to a
    biholomorphism of $D_a-S_{a,h}$ onto $W_a$.
    In fact, if we identify $D_a$ with $\hat\C-\bar\D$, then 
    $D_a-S_{a,h}$ equals the set $H_a$ and
    $\Psi|_{D_a-S_{a,h}}=\psi_a$ defined in~(7.1).
    Let $S_h=\bigcup_{a\in J_\Pi}S_{a,h}$. Clearly
    $S_{a,h}\cap A_{0,h}=\emptyset$.
    By analytic continuation in $D_a$ we get
    $f \circ\Psi=\Psi\circ f_h$ on $W^s(J_\Pi,f_h)-S_h$.
    We show that $S_h$ is closed in $W^s(J_\Pi,f_h)$.
    Take $x\in W^s(J_\Pi)-S_h$. Thus $x\in D_a-S_{a,h}$ 
    for some $a\in J_\Pi$. Pick $0<s<t<1$ such that
    the point $sx\in D_a-S_{a,h}$ and let $\gamma$ be the gradient
    line containing $x$. Choose $n\ge 0$ such that
    $R_0<\deg^nsx<\deg R_0$ and let $\gamma_1$ be the gradient line
    in $A_{0,h}$ containing $f_h^n(x)$. Then
    $\gamma'=\Psi(\gamma)$ is a gradient line for $G$
    containing $\Psi(sx)$ and ending at $a$ and
    $\gamma_1'=\Psi(\gamma_1)$ is a gradient line in $A_0$
    containing $f^n(\Psi(x))$. There is a neighborhood $U$
    of $\Psi(tx)$ consisting of gradient lines $\gamma''$
    for $G$ such that $r(\gamma'')=\infty$ and
    $s(\gamma'')<t$. Let $h$ be a branch of $f_h^{-n}$
    defined near $\gamma_1$ such that $g(f_h^n(x))=x$. If
    $U$ is small enough, then $g\circ\Psi^{-1}(U)$ is
    a neighborhood of $x$ in $W^s(J_\Pi,f_h)$ disjoint from
    $S_h$. Thus $S_h$ is closed in $W^s(J_\Pi,f_h)$.
   Since $S_h$ intersects each disk of $W^s(J_\Pi,f_h)$ in a nowhere dense
set, it is nowhere dense.

To finish the proof, we note that {(3)} is a consequence of~(1) and~(2), and 
 {(4)} follows from Corollary~6.5, from Propostion~4.2, and from~(2).
\qed

\give Remark.  Both the external rays and the global model are not
``canonical'' in the sense that they were not defined by dynamical
behaviors.  In the definitions of $\cE$ and the hedgehog sets $H_a$, it
is equally natural to replace the ``gradient'' lines by the family of
curves in $W^s(a)$ which cross the level sets of $G|W^s(a)$
at a constant angle $\tau$.  The flexibility of considering values of
$\tau$ different from ${\pi\over2}$ has proved useful in the case $\dim=1$
(see [Le]).  

\section 8 Properties of the Support of $T^{\dim-1}$

In this Section we show that under rather general conditions the (local)
stable disks given by Pesin theory actually determine the set
$A\cap\supp(T^{\dim-1})$.  First we show that
$W^s(a)$ is dense in the support of
$T^{\dim-1}\contract A$ for $\mu_\Pi$-a.e.~$a$. As in the case of
polynomial automorphisms of $\cx2$ ([BLS, Prop~2.9])
we do this by proving a convergence result for currents. See
also~[FS4, Cor~5.13].  The other result of this Section is that when the
critical measure $\mu_c$ vanishes on $A_n$, the closure of the Pesin family of
disks (and thus support of
$T^{\dim-1}\contract A_n$) has a uniformly laminar structure.

First we prove some convergence results for measures on $J_\Pi$.
The main tool is Lemma~8.1 below, due to Forn{\ae}ss and Sibony.
Given $a\in J_\Pi$ and $j\ge 0$, define the measures
$\nu_{a,j}$ and $\nu_{a,j}'$ by
$$
\nu_{a,j}':={1\over\deg^{(\dim-1)j}}(f_\Pi^j)^*\delta_a
={1\over\deg^{(\dim-1)j}}\sum_{b\in f_\Pi^{-j}a}\delta_b
$$
and
$$
\nu_{a,j}:=\nu_{f_\Pi^ja,j}'
={1\over\deg^{(\dim-1)j}}\sum_{b\in f_\Pi^{-j}f_\Pi^ja}\delta_b.
$$

\proclaim Lemma~8.1 ([FS3, Lemma~8.3]). There is a constant $C>0$
such that if $\phi$ is a $C^2$ test function on $\Pi$ and $s>0$, and 
$$
E'(\phi,s,j):=\{a\in J_\Pi: 
|\langle\nu_{a,j}',\phi\rangle-\langle\mu_\Pi,\phi\rangle|>s\},
$$
then
$$
\mu_\Pi(E'(\phi,s,j))\le{C|\phi|_{C^2}\over{\deg^j s}}.
$$

\proclaim Lemma~8.2. As $j\to\infty$, we have
$\nu_{a,j}\to\mu_\Pi$ for $\mu_\Pi$-a.e.\ $a$. If $f_\Pi$
is expanding on $J_\Pi$, then $\nu_{a,j}\to\mu_\Pi$
for every $a\in J_\Pi$.

\give Proof. Fix $\phi\in C^2$. It is sufficient to prove that
$\langle\nu_{a,j},\phi\rangle\to\langle\mu_\Pi,\phi\rangle$
for almost every $a$. Let
$$
E_j:=\{a\in\Pi: 
|\langle\nu_{a,j},\phi\rangle-\langle\mu_\Pi,\phi\rangle|>\deg^{-j/2}\}.
$$
By applying Lemma~8.1 with $s=\deg^{-j/2}$ and using 
the invariance of $\mu_\Pi$ we get
$\mu_\Pi(E_j)\le C\deg^{-j/2}$ (the industrious reader may check
that Lemma~8.3 in~[FS3] remains valid if $s$ depends on $j$).
It follows that the set of $a$ such that $a\notin E_j$ for sufficiently
large $j$ has full measure.
Clearly $\langle\nu_{a,j},\phi\rangle\to\langle\mu_\Pi,\phi\rangle$
for these $a$.

If $f_\Pi$ is expanding on $J_\Pi$, then there is an $\epsilon>0$
such that all branches of $f_\Pi^{-j}$ are single-valued on
balls in $\Pi$ of radius $\epsilon$ centered at points in $J_\Pi$. 
Further, the diameters of the preimages under $f^j$
of these balls tend to zero uniformly as $j\to\infty$. It follows
that if $d(a,a')<\epsilon$, then 
$\langle\nu_{a,j}',\phi\rangle-\langle\nu_{a',j}',\phi\rangle\to 0$.

Now let $a\in J_\Pi$ be given. If $j$ is large enough, then
there is a point $b_j\notin E_j$ close to $f_\Pi^ja$. Thus
$$
\langle\nu_{a,j},\phi\rangle=
\langle\nu_{f_\Pi^ja,j}',\phi\rangle\to 
\langle\mu_\Pi,\phi\rangle
{\rm\ as\ }j\to\infty.
$$
\qed

\give Remark. The statements of Lemma~8.2 also hold with the measures
$\nu_{a,j}$ replaced by $\nu_{a,j}'$.

\proclaim Corollary~8.3. $W^s(a,f_\Pi)$ is dense in $J_\Pi$
for $\mu_\Pi$-a.e.\ $a\in J_\Pi$. If $f_\Pi$ is expanding
on $J_\Pi$, then this holds for every $a\in J_\Pi$.

Now we consider convergence of currents.

\proclaim Proposition~8.4. For almost every $a\in J_\Pi$ we have
$$
\liminf_{j\to\infty}{1\over\deg^{(\dim-1)j}}f^{j*}[W_{f_\Pi^ja}]
\ge T^{\dim-1}\contract A.\eqno(8.1)
$$
If $f_\Pi$ is expanding on $J_\Pi$, then
$$
\lim_{j\to\infty}{1\over\deg^{(\dim-1)j}}f^{j*}[W_{f_\Pi^ja}]
=T^{\dim-1}\contract A\eqno(8.2)
$$
for every $a\in J_\Pi$.

\give Proof.
Fix $n\ge 0$. It suffices to show~(8.1) and~(8.2)
on $A_n$.

Fix $\eta>0$ and let $F=F_\eta\subset E$ be a Pesin box with
$\mu_\Pi(F)\ge 1-\eta$ satisfying (a)--(d) above.
Since the Pesin disks $W^s_{-m}(b)=W_b\cap A_{-m}$ depend
continuously on $b$ on $F$, it follows from
Lemma~8.2 that
$$
{1\over\deg^{(\dim-1)j}}\sum_{b\in f_\Pi^{-j}f_\Pi^ja\cap F}
[W_b\cap A_{-m}]\to
\int_F[W_b\cap A_{-m}]\,\mu_\Pi(b)
$$
as $j\to\infty$ for almost every $a$.
Hence, if we define
$$
X_j(a)={1\over\deg^{(\dim-1)j}}f^{j*}[W_{f_\Pi^j(a)}],
$$
then 
$$
\liminf_{j\to\infty} X_j(a)\contract A_{-m}\ge
\int_F[W_b\cap A_{-m}]\,\mu_\Pi(b)
$$
for almost every $a$. If we pull this back by $f^{n+m}$, then we get
$$
\liminf_{j\to\infty} X_j(a)\contract A_n\ge
\int h[W_b\cap A_n]\,\mu_\Pi(b),\eqno(8.3)
$$
where $h=h_\eta=\deg^{-(\dim-1)(m+n)}f^{n+m}_*\chi_F$.
In particular we have $0\le h\le 1$ and $\int h=\mu_\Pi(F)$.
By letting $\eta\to 0$ we get $h_\eta\to 1$ so by dominated
convergence we find that the right hand side of~(8.3)
converges to
$\int [W_b\cap A_n]\,\mu_\Pi(b)=T^{\dim-1}\contract A_n$
for almost every $a$.

If $f_\Pi$ is expanding on $J_\Pi$, then we may take $F=J_\Pi$
and use the second part of Lemma~8.2. Thus
$X_j(a)\contract A_n\to T^{\dim-1}\contract A_n$ for every
$a\in J_\Pi$.
\qed

\proclaim Corollary~8.5. For almost every $a\in J_\Pi$ we have
$\overline{W^s(a)}=\supp(T^{\dim-1}\contract A)$. If $f_\Pi$
is uniformly expanding on $J_\Pi$, then this holds for every
$a\in J_\Pi$.

\give Proof. It is clear that $W_a\subset W^s(a)$ for all $a\in E$.
We claim that $W_a\subset\supp(T^{\dim-1}\contract A)$. By the 
construction of $W_a$ it suffices to show that $W^s_0(a)$
is contained in the support
of $T^{\dim-1}\contract A$ for $a\in E_0$. But this follows
from the continuity of $W^s_0(a)$ on $E_0$, from the fact that
$E_0$ has no isolated points and from~(6.4).

Thus $\overline{W^s(a)}\subset\supp(T^{\dim-1}\contract A)$.
The reverse inclusion is a consequence of Proposition~8.4.
\qed

If $f_\Pi$ is uniformly expanding on $J_\Pi$, then by increasing
$R_0$, we have $\mu_c(A_0)=0$. We now consider the property
$$
\mu_c(A_n)=[\cC]\wedge(T^{\dim-1}\contract A_n)=0
$$
for some $n\in\Z$, without assuming uniform expansion on $J_\Pi$.
By Corollary~6.5, this implies that
$\cC\cap Z_{a,n}=\emptyset$ for $\mu_\Pi$ almost every $a\in E$.
Using the invariance of $\mu_\Pi$ and the fact that
$fZ_{a,n}\subset Z_{f_\Pi a,n}$ we also get that
$Z_{a,n}\cap\cC_\infty=\emptyset$ for almost every $a\in E$,
where $\cC_\infty=\bigcup_{j\ge 0}f^{-j}\cC$.

For such $a$, the construction of $W_a$ involves the removal
of no gradient lines in $A_n$. Thus $W_a\cap A_n=Z_{a,n}$, and
$W_a\cap A_n$ is a properly embedded disk in $A_n$.
Further, for these $a$, the mapping $\psi_a$, defined
by~(7.1) maps the disk
$\Delta_n=\{|\zeta|>\exp(\deg^{-n}R_0)\}$ biholomorphically
onto $W_a\cap A_n$. Let us summarize this.

\proclaim Proposition~8.6.
\item{(1)} If $\mu_c(A_n)=0$, then for almost every $a$,
  $W_{a,n}:=W_a\cap A_n$ is a properly embedded disk in $A_n$, and
  the restriction of $\psi_a$ to
  $\Delta_n$ is a biholomorphism onto $W_{a,n}$.
\item{(2)} If $\mu_c(A)=0$, then for almost every $a$
  $W_a$ is a properly embedded disk in $A$, and
  $\psi_a$ maps $\Delta=\{|\zeta|>1\}$
  biholomorphically onto $W_a$. Further,
  $\psi_a:\Delta\to\px\dim-J$ is proper.

\give Proof. Everything except the last statement follows from the
discussion above. By Theorem~7.3, it follows that for $\mu_\Pi$-almost
every $a$ the boundary values of the disk
$\psi_a:\Delta\to W_a$ lie inside $J$ for almost every $\theta$.
Thus $\psi_a:\Delta\to\px\dim-J$ is proper by
the theorem of Alexander (see [A] and [Ro]).
\qed 

Suppose $\mu_c(A_n)=0$. We want to show that the family of
disks $W_{a,n}$ extends to a Riemann surface lamination in $A_n$.
Let $\Gamma_n$ denote all the uniformizing
mappings $\psi_a:\Delta_n\to W_{a,n}$.

\proclaim Lemma~8.7.  Either $\Gamma_n$ is a normal family, or there
is a nonconstant holomorphic mapping $h:\C\to\Pi$ such that
$h(\C)\subset\Pi-E$.

\give Proof.  If $\Gamma_n$ is not a normal family, there is a sequence
$\{\psi_j\}\subset\Gamma_n$ without a convergent subsequence.
By the renormalization technique of Brody [La, pp.\ 68--71]
there is a sequence $r_j\to\infty$ and a
sequence of M\"obius transformations $\rho_j:\{|\zeta|<r_j\}\to\Delta_n$
such that $\psi_j\circ\rho_j$ converges to a nonconstant mapping
$h:\C\to\px\dim$. Further, $G\circ \psi_j\circ\rho_j$ is a
sequence of positive, superharmonic functions which
converge normally on $\C$ to $G\circ h$. Since every positive,
superharmonic function on $\C$ is constant (or ${\equiv}+\infty$),
it follows that either $h(\C)\subset\{G=c\}$ or
$h(\C)\subset\Pi$. The set $\{G=c\}$ is a compact subset of $\cx\dim$, so
we cannot have $h(\C)\subset\{G=c\}$. Thus we have $h(\C)\subset\Pi$.
Finally, the disks $\{W_a\}$ are pairwise disjoint, so an application
of the Hurwitz Theorem shows that either
$h(\C)\subset W_a$ or $h(\C)\cap W_a=\emptyset$. The first case is not
possible since $h$ is nonconstant, so we have $h(\C)\subset\Pi-\{a\}$.
Thus $h(\C)\subset\Pi-E$.
\qed

\give Remark. If $\dim=2$, there can be no nonconstant mapping of $\C$
into $\Pi-E$, so $\Gamma_n$ is necessarily a normal family.
For $\dim>2$ we let $X_1,X_2,\dots$ denote the distinct irreducible
components in $\bigcup_{j\in\Z}f_\Pi^j\cC_\Pi$.
Since $\mu_c(A_n)=0$ we conclude, using the Hurwitz Theorem, that for any
$X_j$ we have either $h(\C)\subset X_j$ or $h(\C)\subset\Pi-X_j$.
Thus $\Gamma_n$ is a normal family if for any finite family
$X_{j_1},\dots,X_{j_p}$, there is no nonconstant mapping of $\C$
into $X_{j_1}\cup\dots\cup X_{j_p}-\bigcup'X_i$, where $\bigcup'$
denotes the union of the remaining components.  Since normality of
$\Gamma_n$ is the essential ingredient of the
following theorem, it would appear to apply in many cases with $\dim>2$.

\proclaim Theorem~8.8. If $\mu_c(A_n)=0$ and $\dim=2$, then for each
$a\in J_\Pi$ there is a complex disk $W_{a,n}$ which is properly embedded
in $A_n$, such that $a\mapsto W_{a,n}$ is continuous, and such that
$W_{a,n}\cap W_{b,n}=\emptyset$ if $a\ne b$. Thus the
family $\{W_{a,n}:a\in J_\Pi\}$ is a Riemann surface lamination in $A_n$.
Further, we have
$$
T\contract A_n=\int [W_{a,n}]\,\mu_\Pi(a)\eqno(8.4)
$$
and for each $a$ we have $W_{a,n}\cap\cC=\emptyset$ or
$W_{a,n}\subset\cC\cup\{a\}$.
If $\mu_c(A)=0$, then all of the above conclusions hold on $A$.

\give Proof. By Lemma~8.7 and the remark above,
$\Gamma_n$ is a normal family.
Let $\tilde\Gamma_n$ denote the mappings $\psi:\Delta_n\to A_n$
which are normal limits of $\Gamma_n$. Note that
$G\circ\psi(\zeta)=\log|\zeta|$
for any $\psi\in\tilde\Gamma_n$. Since the disks $W_a$ are
pairwise disjoint, it follows from the Hurwitz Theorem that if  
$\psi',\psi''\in\tilde\Gamma_n$, then either
$\psi'(\Delta_n)=\psi''(\Delta_n)$ or
$\psi'(\Delta_n)\cap\psi''(\Delta_n)=\emptyset$.
Thus for each $a\in J_\Pi$ there is a
unique image $W_{a,n}:=\psi(\Delta_n)$,
$\psi\in\tilde\Gamma_n$, which contains $a$. The continuity of
$a\mapsto W_{a,n}$ follows from the normality of $\Gamma_n$.
We have $Z_{a,n}=W_{a,n}$ and $N_n(a)=1$ for almost every $a$
and hence~(8.4) follows from~(6.4). Applying Lemma~2.2 to~(8.4)
we obtain
$$
\mu_c\contract A_n=\int [W_{a,n}\cap\cC]\,\mu_\Pi(a)=0.\eqno(8.5)
$$
By continuity of $a\to W_{a,n}$, the property that
$W_{a,n}\cap\cC\ne\emptyset$ but $W_{a,n}-(\cC\cup\{a\})\ne\emptyset$
is open in $J_\Pi$. Thus this property never holds by~(8.5).
\qed

\give Example.  Let $f(z,w)=(z^2,{1\over 4}z^2+w^2)$.  In the coordinate
$\zeta=w/z$ on $\Pi$, we have $f_\Pi(\zeta)=\zeta^2+{1\over 4}$.  Let
$K_\Pi\subset\Pi$ denote the filled Julia set for $f_\Pi$.  The point
$\zeta={1\over2}$ is a parabolic fixed point, and all points of
$\{{1\over2}\}\cup{\rm int}(K_\Pi)\subset\Pi$ approach $\{\zeta={1\over 2}\}$
in forward time.  Thus the stable set $W^s({1\over2})$ for $f$ contains a
neighborhood of ${1\over 2}$ inside the cone of complex lines
$C(\{{1\over2}\}\cup{\rm int}(K_\Pi))$, which contains an open set in
$\px\dim$.
 
Since $f=f_h$ is homogeneous, each (local) Pesin disk $W^s_{-m}(a)$ is the
complement of a closed disk (centered at the origin) inside the complex line
$L_a$.  Thus the family of Pesin disks has an extension to the lamination
inside the complex cone of lines
$C(J_\Pi)$.  In the example at hand, the critical locus is 
$\cC=\{z=0\}\cup\{w=0\}$, and
$T_h$ is supported on
$C(J_\Pi)$.  Since
$\{\zeta=0,\infty\}$ is disjoint from $J_\Pi$, and 
$A$ is a neighborhood of
$\Pi$ disjoint from $0\in\cx2$, it follows that $\mu_c\contract A=[\cC]\wedge
T_h\contract A=0$.  Thus this example also satisfies the hypotheses of
Theorem~8.8.
\bigskip

%
%

\section 9 Axiom A in $\cx 2$.

In the next section we will impose certain hyperbolicity assumptions
(see Definition~10.1) on the dynamics on $f$ in
order to prove that all of the external rays land
(and land continuously) on $J$. Most of these assumptions are 
related to Axiom
A, which was introduced by Smale as a property of a smooth dynamical
system which enables the understanding of its global dynamics.
In this Section, we discuss Axiom A in the setting of polynomial
endomorphisms of $\cx2$, chiefly to clarify our assumptions in
Definition~10.1.

The literature on hyperbolic dynamics is
vast, but most expositions consider only diffeomorphisms. A regular
polynomial endomorphism of $\cx 2$ of degree $\deg\ge 2$ is not
invertible, and the hyperbolic theory is slightly different.
There seems to be no general, detailed treatment of exactly the
results we need, so we will give
further definitions and results in Appendix B. More
details can be found in [J2]. We also refer to [FS4], where the authors
study hyperbolic endomorphisms of $\px 2$.

Suppose that $f$ is a regular polynomial endomorphism of $\cx 2$; as
usual we regard $f$ as a holomorphic map of $\px 2$.
Since $f$ is not injective, we will often have to work with histories
of points instead of the points themselves. Precisely, a
{\it history\/}
of a point $x\in\px 2$ is a sequence $(x_i)_{i\le 0}$ of points
in $\px 2$ such that $x_0=x$ and $fx_i=x_{i+1}$ for all $i<0$.
We will use the notation $\hat x$ for a history $(x_i)$.

Let $L$ be a compact subset of $\px 2$ with $fL=L$. We refer to 
Appendix B for a definition of what it means for
$f$ to be (uniformly)
hyperbolic on $L$. Let us only recall that the definition
involves the compact set $\hat L$ of histories in $L$. The pair
$(\hat L,\hat f)$, where $\hat f$ is the left shift on $\hat L$, 
is often called the natural extension of $f|_L$. There is a natural
projection $\pi:\hat L\to L$ such that $\pi(\hat x)=x_0$.
We say that $L$ has unstable index
$i$ if the stable bundle $E^s$ has constant dimension $2-i$ on $L$.
If $L$ has unstable index 2, then $f$ is said to be (uniformly)
expanding on
$L$ (see Appendix B for an alternative definition).
If $f$ is hyperbolic on $L$, then to every point in $x\in L$
and every history $\hat x\in\hat L$ there is an associated 
local stable and unstable manifold respectively, defined by
$$
\eqalign{
W^s_{\rm loc}(x)
&=\{y\in\px 2: d(f^iy,f^ix)<\delta\ \forall i\ge 0\}\cr
W^u_{\rm loc}(\hat x)
&=\{y\in\px 2: \exists\hat y\in\widehat{\px 2}, \pi(\hat y)=y,
d(y_i,x_i)<\delta\ \forall i\le 0\},}
$$
for small $\delta>0$.
Then $W^s_{\rm loc}(x)$ and $W^u_{\rm loc}(\hat x)$ are complex disks
of $\px 2$. If $f$ is uniformly expanding on $L$,
then the local stable manifolds 
are empty and the local unstable manifold at $\hat x$ is a 
neighborhood of $x_0$ in $\px2$.

We also define global stable and unstable manifolds by declaring
$$
\eqalign{
W^s(x)
&=\{y\in\px 2: d(f^iy,f^ix)\to 0\ {\rm as}\ i\to\infty\}\cr
W^u(\hat x)
&=\{y\in\px 2: \exists\hat y\in\widehat{\px 2}, \pi(\hat y)=y,
d(y_i,x_i)\to 0\ {\rm as}\ i\to-\infty\}.}
$$
Note that if $n\ge 0$, $y\in L$ and $f^ny=f^nx$, then 
$W^s(x)$ contains $W^s_{\rm loc}(y)$. Hence the global
stable manifolds are in general large and quite complicated objects
(compare with Corollary~8.5).
Both the stable and unstable manifolds may have singularities; 
this is in contrast to the case of polynomial automorphisms of $\cx 2$,
where they are immersed copies of $\C$ [BS1], 

We now turn to Axiom A regular polynomial endomorphisms of $\cx 2$.
A point $x\in\px 2$ is wandering 
if for every neighborhood $V$ of $x$ there exists an $n\ge 1$ such that
$f^n(V)\cap V\ne\emptyset$. The {\it non-wandering set\/} $\Omega$ of $f$ 
is the set of all non-wandering points; it is a compact set.
A regular polynomial endomorphism $f$ of $\cx 2$ is
{\it Axiom A\/} if the periodic points of $f$ are dense in $\Omega$ and
$f$ is hyperbolic on $\Omega$. If $f$ is Axiom A, then
Smale's spectral decomposition theorem (Theorem~B.9) asserts that
$\Omega$ can be written in a unique way as a finite union of 
disjoint compact invariant sets $\Omega_j$, called basic sets, such that
$f|_{\Omega_j}$ is transitive, i.e.\ has a dense orbit. Thus each basic
set has a well-defined unstable index.

Let us investigate what the possible basic sets are for a Axiom A 
regular polynomial endomorphism $f$ of $\cx 2$. To do this, we
first observe that the four sets $\Pi$, $\cx 2-K$, ${\rm int}(K)$
and $\partial K$ are all completely invariant and see what basic sets
each one of them may contain.

To begin with, it is clear that $\Omega(f)\cap\Pi=\Omega(f_\Pi)$.
Now $f_\Pi$ is a rational map and from one-dimensional
dynamics we know that $f_\Pi$ is Axiom A if and only if $f_\Pi$
is uniformly expanding on $J_\Pi$ (see [M]). Hence, if 
$f$ is Axiom A, then the basic sets in $\Pi$ are $J_\Pi$, which
is of unstable index 1, and a finite union of attracting periodic
points, all of whose unstable index is zero.

All the points in the open set $\cx 2-K$ are attracted to $\Pi$ so
$(\cx 2-K)\cap\Omega$ is empty. It is clear that
$\{f^n\}$ is normal on the interior of $K$, so if $f$ is Axiom A, then
the only basic sets in ${\rm int}(K)$ are attracting periodic
points, all of whose unstable index are zero.

The boundary of $K$ contains the most complicated dynamics. 
Clearly, no basic sets in $\partial K$ can have unstable index 0. 
Let $S_2$ and ${S_1}$ be the 
union of the basic sets in $\partial K$ of index 2 and 1,
respectively. We note that ${S_1}$ can be empty, as in the
example $f(z,w)=(z^2+c,w^2+c)$, with $c$ outside
the Mandelbrot set.
On the other hand, $J$ is a basic set of unstable
index 2 (see [FS2, Theorem~7.4]), so $J\subset S_2$.
The question arises whether this inclusion is ever strict or,
equivalently, whether $f$ can have repelling periodic points
outside $J$. Hubbard and Papadopol [HP] have in fact given an example
of a regular polynomial endomorphism of $\cx 2$ with a repelling 
periodic point outside $J$ but is seems difficult to check 
whether their map can be made Axiom A. In any case we have the
following.

\proclaim Lemma~9.1. Let $f$ be an Axiom A regular polynomial
endomorphism of $\cx2$. Then $f^{-1}S_2=S_2$ if and
only if $S_2=J$, i.e.\ if all repelling periodic points are 
contained in $J$.

\give Remark. A proof is given in [FS4]. We give it here for the
convenience of the reader.

\give Proof. The ``only if'' part is trivial since $f^{-1}(J)=J$, so 
suppose that $f$ is Axiom A and $f^{-1}(S_2)=S_2$ but $S_2\ne J$.
Let $N$ be an open neighborhood of $J$ such that 
$f^{-1}(N)\subset N$ and $\cap_{n\ge 0}f^{-n}(N)=J$.
Then $N-J$ contains only wandering points, so 
$S_2-J$ is at a positive distance from $J$ and is therefore a
completely invariant compact set. Let $N'$ be an 
open neighborhood of $S_2-J$ disjoint from $J$ with 
$f^{-1}(N')\subset N'$. Then $N'$ has positive capacity and
if $x\in N'$ then $(f^n)^*\delta_x/\deg^{2n}$ cannot converge to
$\mu$ as $n\to\infty$. This contradicts Lemma~8.3 in [FS3].
\qed

Let $f$ be an Axiom A regular polynomial endomorphism
of $\cx 2$ with $f^{-1}S_2=S_2$. It follows from Corollary~B.10
and the above discussion that any history of a point $\cx 2$
which is not an attracting periodic point must converge to either
$J$ or ${S_1}$.
We define the unstable set of $J$ to be the
set of points in $\cx 2$ all of whose histories converge to $J$, i.e.\
$$
W^u(J)=\{x\in\cx 2: (\hat x\in\widehat{\cx 2}, \pi(\hat x)=x)
\Rightarrow x_i\to J\ {\rm as}\ i\to-\infty\}.
$$
We note that this definition differs from the one in [FS4], where
$W^u(J)$ is defined as the set of points having {\sl at least}
one history converging to $J$.
On the other hand we define the unstable set of ${S_1}$
as
$$
W^u(S_1)=\{x\in\cx 2: \exists\hat x\in\widehat{\cx 2}, \pi(\hat x)=x,
x_i\to S_1\ {\rm as}\ i\to-\infty\}.
$$
Let $N$ be a neighborhood of $J$ in $\cx 2$ as in the proof of Lemma~9.1.
Clearly $N\subset W^u(J)$ and every point in
$\cx 2$ which is not an attracting periodic point is 
contained in precisely one of the sets $W^u(J)$ and $W^u({S_1})$.

\proclaim Lemma~9.2. If $x\in W^u(J)$, then there exists an $n\ge 0$ such
that $f^{-n}(x)\subset N$. In particular, $W^u(J)$ is open in $\cx 2$
and $W^u({S_1})$ is closed in $\cx 2$ except possibly at some
of the attracting periodic points.

\give Proof. Let $Z$ be the set of points $y$ in $\cx 2$ such that 
for all $n\ge 0$, there is a point in $f^{-n}(y)$ outside $N$. 
It is clear that if $y\in Z$, then $y$ has at least one preimage in
$Z$, so every point $y\in Z$ has a whole history inside $Z$. Such a 
history cannot converge to $J$ so it follows that $Z\cap W^u(J)=\emptyset$,
which completes the proof.
\qed

For the proof of the main result in~\S10 (Theorem~10.2),
we will work with slightly weaker hyperbolicity hypotheses.

\proclaim Definition 9.3. A regular polynomial 
endomorphism $f$ of $\cx 2$ satisfies condition (\dag) if the
following four properties hold:
\item{(\dag1)} $f_\Pi$ is uniformly expanding on $J_\Pi$.
\item{(\dag2)} $f$ is uniformly expanding on $J$.
\item{(\dag3)} The nonwandering set of $f$ in $\partial K$ consists
 of $J$ and a hyperbolic set $S_1$ of unstable index 1.
\item{(\dag4)} $W^u(S_1)=\bigcup_{\hat{x}\in\hat{S_1}}W^u(\hat{x})$.

\proclaim Proposition~9.4. Let $f$ be an Axiom A regular
polynomial endomorphism of $\cx 2$ with $f^{-1}S_2=S_2$.
Then $f$ satisfies condition (\dag).

\give Proof. From the above discussion we know that $f$
satisfies conditions (\dag1), (\dag2) and (\dag3), and
(\dag4) follows from Corollary~B.10.
\qed
%
%

\section 10 Continuous landing of rays in $\cx2$.

So far we have been able to understand the
dynamics in the set $W^s(J_\Pi)$, or at least on the support of
$T^{\dim-1}\contract A$, for rather general $f$. In this Section, we
approach the dynamics of $f$ on $J$ by proving Theorem~10.2, which shows
that $e:\cE\to J$ is a continuous surjection (under suitable
assumptions). Our approach is restricted
to the case $\dim=2$ because we work with unstable manifolds $W^u(\hat q)$
as Riemann surfaces; if $\dim>2$, the unstable manifolds can have
dimension $>1$.  If $\dim=2$, then $\Pi$ is the Riemann sphere,
and $f_\Pi$ is a rational mapping. Up to~\S8, the only hyperbolicity
assumption that we have been concerned with has been uniform expansion
on $J_\Pi$, which for $\dim=2$ means that $f_\Pi$ is a hyperbolic
rational mapping. The reason for this is that we have dealt with the
dynamics on $A$ and not directly on $K$. For the continuity of $e:\cE\to J$
we need to consider the dynamics on $K$ (or, rather, $\partial K$).
Notice that hyperbolicity of $f_\Pi$ does not exclude complicated
dynamics on $K$. 


To motivate (\ddag5) in Definition~10.1 below,
let us revisit the example presented after Theorem~7.3, namely
$f(z_1,z_2)=(z_1^2+c_1,z_2^2+c_2)$. We have
$J=J_{c_1}\times J_{c_2}$, where $J_c$ is the one-dimensional
Julia set of $p_c(z)=z^2+c$. Further,
$\cE\simeq S^1\times S^1$, so $\cE$ is connected and
locally connected. Thus, if $e$ maps $\cE$ continuously onto
$J$, then $J$ also is connected and locally connected. This,
in turn, is equivalent to $J_{c_j}$ being connected and locally
connected for $j=1,2$. Now $J_c$ is connected if and only if the
critical point $0$ of
$p_c$ is not in the basin of attraction of infinity, i.e.\
the parameter value $c$ is in the Mandelbrot set.
Using this one can see that $J$ is connected if and only if
$W^s(J_\Pi)\cap\cC=\emptyset$.
The question of whether $J_c$ is locally connected is more delicate,
but a sufficient condition is that $J_c$ is connected and
$p_c$ is uniformly expanding on $J_c$. Thus $J$ is locally connected
if $J$ is connected and $f$ is uniformly expanding on $J$. 
 
\proclaim Definition 10.1. We say that a regular polynomial 
endomorphism $f$ of $\cx 2$ satisfies condition (\ddag) if $f$ 
satisfies condition~(\dag) in Definition~9.3, and
$W^s(J_\Pi)\cap\cC=\emptyset$, i.e. if the
following five properties hold:
\item{(\ddag1)} $f_\Pi$ is uniformly expanding on $J_\Pi$.
\item{(\ddag2)} $f$ is uniformly expanding on $J$.
\item{(\ddag3)} The nonwandering set of $f$ in $\partial K$ consists
  of $J$ and a (possibly empty) hyperbolic set $S_1$ of unstable
  index 1.
\item{(\ddag4)} $W^u(S_1)=\bigcup_{\hat{x}\in\hat{S_1}}W^u(\hat{x})$.
\item{(\ddag5)} $W^s(J_\Pi)\cap\cC=\emptyset$
  (or, equivalently, $\mu_c(A)=0$).

Let us comment on these conditions.  It follows from Proposition~9.4 that if $f$ is Axiom A,
$f^{-1}(S_2)=S_2$, and satisfies~(\ddag5), then $f$ satisfies~(\ddag).
Using this, one can show that perturbations of the map
$f(z,w)=(z^\deg,w^\deg)$ satisfy~(\ddag).

Conditions 
(\ddag1) and~(\ddag5) guarantee that
$e_r:\cE\to\{G=r\}$ is well defined and continuous for $r>0$
(in general it is defined almost everywhere on $\cE$). 
We have $e=\lim e_r$ as $r\to 0$,
and $e_r(\cE)=f^{-1}e_{\deg r}(\cE)$. Thus, once we know that
$e_r(\cE)$ is in a small neighborhood of $J$ for all sufficiently
small $r$, then
condition~(\ddag2) helps us to show that $e_r$ converges uniformly
as $r\to 0$. However, we only know that $e_r$ accumulates on
$\partial K$, which in general is a larger set than $J$. In fact,
the main difficult in proving continuity for $e$ is to show that
$e_r(\cE)$ accumulates only at $J$. In order to do this, we use
properties~(\ddag3) and~(\ddag4).

Let us make some remarks about the connection between the endpoint map $e$
and the conjugacy $\Psi:W^s(J_\Pi,f_h)\to W^s(J_\Pi,f)$ between $f_h$ to
$f$ given in Theorem~7.4.
Let $\Delta_a:=L_a\cap A_h$ be the disk in $W^s(J_\Pi,f_h)$
corresponding to $a$. We may identify the set $\cE_a$
of external rays in $W_a$ with $\partial\Delta_a$ and $\cE$
with the boundary of $W^s(J_\Pi,f_h)$, i.e.\ the union of
$\cE_a$ over $a\in J_\Pi$.
This defines the topology on $\cE$. With these identifications
we have $e_r(\gamma)=\Psi(e^r\gamma)$ for $r>0$. It follows
that $e_r$ is well defined and continuous for all $r>0$.
Thus $\Psi$ extends continuously to $\cE$ if and only if
$e$ is continuous, and in this case $e$ coincides with
the restriction of $\Psi$ to $\cE$.
The selfmap $f_h$ on $W^s(J_\Pi,f_h)$ induces the selfmap $\sigma$
of $\cE$.

Our main goal is to prove the following result.

\proclaim Theorem~10.2. If the regular polynomial endomorphism $f$
of $\cx2$ satisfies condition~(\ddag),
then $e$ maps $\cE$ H\"older continuously onto $J$ and
$f\circ e=e\circ\sigma$. In particular, the conclusions
hold if $f$ is Axiom A, if
$W^s(J_\Pi)\cap\cC=\emptyset$, and if the expanding part $S_2$ of the
nonwandering set of $f$ satisfies $f^{-1}(S_2)=S_2$.

As mentioned above, the main difficulty in proving Theorem~10.2
is to show that the external rays accumulate
only at $J$. In particular, there must be
no heteroclinic intersection
between $S_1$ and $J_\Pi$, i.e.\ no complete 
orbit $(x_i)_{i\in\Z}$ such that $x_i\to S_1$ as $i\to-\infty$ and
$x_i\to J_\Pi$ as $i\to\infty$.

\proclaim Lemma~10.3. If $f$ satisfies condition (\ddag), then
$W^s(J_\Pi)\cap W^u({S_1})=\emptyset$.

We postpone the proof of Lemma~10.3 for the moment and head towards
the proof of Theorem~10.2.

We may identify the disk $\Delta_a$ in $W^s(J_\Pi,f_h)$
with $\Delta=\{|\zeta|>1\}$ in such a way that
$G_h(\zeta)=\log|\zeta|$. The restriction of $\Psi$ to $\Delta_a$
then induces a conformal equivalence $\psi_a:\Delta\to W_a$ such that
$G\circ\psi_a(\zeta)=\log|\zeta|$. This last condition
defines $\psi_a$ uniquely up to rotation; the precise choice of
rotation will not be important in what follows. Given a choice of
$\psi_a$, the conjugacy between $f_h$ and $f$ translates into
$$
f\circ\psi_a(\zeta)=\psi_{f_\Pi a}(\nu_a\zeta^\deg),\eqno(10.1)
$$
where $|\nu_a|=1$.

We first show that the maps $\psi_a$ are uniformly H\"older continuous.

\proclaim Lemma~10.4. There exist constants $\alpha>0$ and $C>0$ such that
$$
d(\psi_a(\zeta),\psi_a(\zeta'))\le Cd(\zeta,\zeta')^\alpha,\eqno(10.2)
$$
for all $a\in J_\Pi$ and $\zeta,\zeta'\in\Delta$.

\give Proof.
The expansion of $f$ on $J$ implies that there exists a neighborhood $N$
of $J$ with $f^{-1}(N)\subset N$, $\lambda>1$ and a metric equivalent to
the Euclidean metric such that
$|Df(x)v|\ge\lambda|v|$ for all $x\in N$ and all $v\in T_x\cx 2$
with respect to this metric.
By Lemma~9.2 and Lemma~10.3 we know that the set 
$W^s(J_\Pi)\cap\{1\le G\le \deg\}$ is a compact subset of the open set
$W^u(J)$ so by pulling back by $f$ we see that there exists an
$R>1$ such that $\psi_a\{1<|\zeta|\le R\})\subset N$ for all $a$.
Let $\alpha>0$ be so small that $\deg^\alpha<\lambda$ and assume
that $R$ is so small that $R^{\deg-1}\deg^\alpha<\lambda$.
It is sufficient to prove~(10.2) for $1<|\zeta|,|\zeta'|\le R$.

By differentiating~(10.1) and using the estimates above we get that, for 
$1<|\zeta|<R^{\deg^{-1}}$,
$$
\left|D\psi_a\left(\zeta\right)\right|\le\lambda^{-1}
\left|D\psi_{f_\Pi a}\left(\nu_a\zeta^{\deg}\right)\right|
\deg|\zeta|^{\deg-1}.\eqno(10.3)
$$
Define 
$$
m(r)=\sup_{a\in J_\Pi}\sup_{|\zeta|=r}|D\psi_a(\zeta)|,
$$
for $1<r\le R$. Then there exists a constant $C'<\infty$ such that
$$
m(r)\le C'(r-1)^{\alpha-1},\eqno(10.4)
$$
for $R^{\deg^{-1}}\le r\le R$. Using the estimate~(10.3) we prove
inductively that~(10.4) holds for $1<r\le R$. Thus~(10.2) follows
by integrating~(10.4).
\qed

\give Proof of Theorem~10.2.
We know that $e_r$ is continuous for each $r>0$.
Lemma~10.4 shows that $e_r$ converges uniformly on $\cE_a$ for each $a$.
Thus $e_r$ converges uniformly to $e$, so $e$ is continuous and it
follows from Theorem~7.3 that $e(\cE)=J$.
We have $f\circ e_r=e_{\deg r}\circ\sigma$, so be letting $r\to 0$
we get $f\circ e=e\circ\sigma$.
It remains to be seen that $e$ is H\"older continuous.

Let $N$ and $\lambda$ be as in the proof of Lemma~10.4 and let
$\epsilon>0$ be small. We may assume that $d(fx,fy)\ge\lambda d(x,y)$
for $x,y\in N$, $d(x,y)\le\epsilon$, and that
$d(\sigma\gamma,\sigma\gamma')\ge\lambda d(\gamma,\gamma')$
for $\gamma,\gamma'\in\cE$, $d(\gamma,\gamma')\le\epsilon$.
There is a number $M$ so that
$d(\sigma\gamma,\sigma\gamma')\le M d(\gamma,\gamma')$ for all
$\gamma,\gamma'\in\cE$. Pick $\alpha>0$ such that
$M^\alpha<\lambda$. Let $C>0$ be so large that
$d(e_r(\gamma),e_r(\gamma))\le Cd(\gamma,\gamma')^\alpha$ if
$d(\gamma,\gamma')\ge\epsilon$ and $r>0$.

Since $e(\cE)=J$, there exists an $R>0$ such that $e_r(\cE)\subset N$ 
for $0<r\le R$. Now suppose $r\le R\deg^{-j}$ for some $j\ge 0$ and
that $d(\gamma,\gamma')\ge\epsilon/\lambda^j$.
Then
$$
\eqalign{
d(e_r(\gamma),e_r(\gamma'))
&\le\lambda^{-j}d(e_{\deg^jr}\sigma^j\gamma,e_{\deg^jr}\sigma^j\gamma')\cr
&\le\lambda^{-j}Cd(\sigma^j\gamma,\sigma^j\gamma')^\alpha\cr
&\le\lambda^{-j}M^{\alpha j}Cd(\gamma,\gamma')^\alpha\cr
&\le Cd(\gamma,\gamma')^\alpha}
$$
The theorem thus follows by letting $r\to 0$.
\qed

We now turn to the proof of Lemma~10.3 and proceed in a number of steps.
An intersection between $W_a$ and $W^s(J_\Pi)$ is a
{\it heteroclinic connection\/}.  There
are two types, as shown in Figure 5. That is, the intersection
between $W_a$ and
$W^u(S_1)$ can be either a discrete set or a relatively open set.
The next Lemma says that the discrete intersection does not occur.

\centerline{\psfig{figure=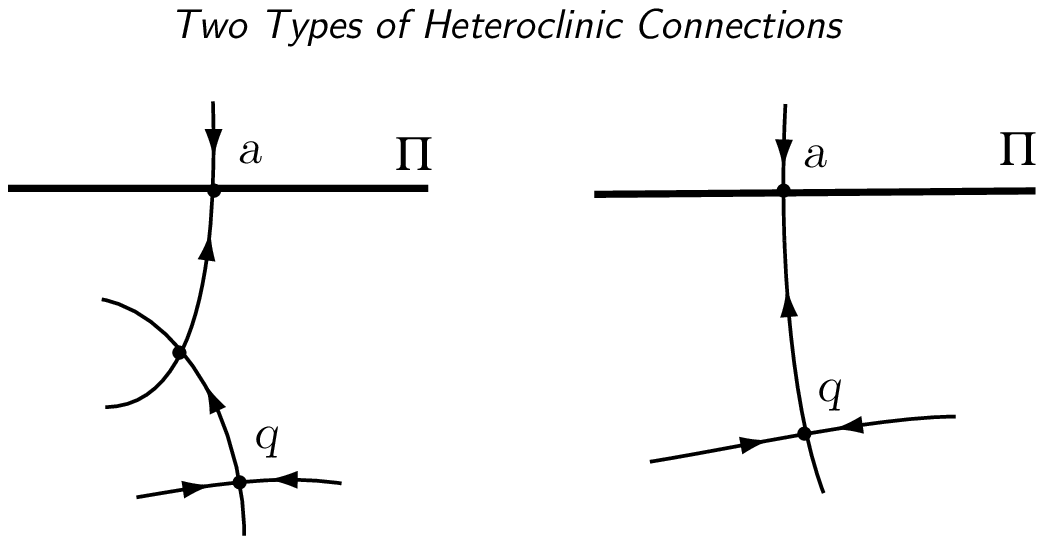,width=0.7\hsize}}
\centerline{Figure 5.}
\medskip

\proclaim Lemma~10.5. Let $W_a$ be the stable disk of a point $a\in J_\Pi$.
Then either $W_a\cap W^u({S_1})=\emptyset$ or there exists a point
$\hat x\in\hat{S_1}$ such that 
$W_a^*\subset W^u(\hat x)$, where $W_a^*=W_a-\{a\}$.

The key observation in proving the dichotomy is the following.

\proclaim Lemma~10.6. If $U$ is a simply connected open subset of 
a punctured stable disk $W_a^*$, then
all branches of $f^{-i}|_U$ for all $i>0$ are well-defined and
holomorphic on $U$ and they form a normal family there.

\give Proof. That the branches are well-defined follows from
condition (\ddag5). If $V$ is relatively compact in $U$
then all branches of $f^{-i}$ on $V$ map $V$ into 
a fixed compact subset of $\cx 2$. Thus they form
a normal family on $U$.
\qed

\give Proof of Lemma~10.5. Suppose that $y\in W_a\cap W^u(S_1)$. 
Then by condition (\ddag4) there exists a point 
$\hat x\in\hat{S_1}$ such that
$y\in W^u(\hat x)$, i.e.\ $y$ has a history $\hat y$ such that
$d(y_i,x_i)\to 0$ as $i\to-\infty$. Let $U$ be any simply connected
open subset of $W_a^*$ containing $y$ and let $g_i$ be the unique
sequence of branches of $f^{-i}|_U$ such that $g_i(y)=y_i$.
Then $\{g_i\}$ is equicontinuous by Lemma~10.6, so there is a 
small neighborhood $V$ of $y$ in $U$ such that the maximal
distance from $g_i(V)$ to $x_i$ is uniformly small
as $i\to\infty$. Hence $V\subset W^u(\hat x)$ and, by normality
of $\{g_i\}$, $U\subset W^u(\hat x)$. Since $U$ 
was arbitrary it follows that $W_a^*\subset W^u(\hat x)$.
\qed

\proclaim Corollary~10.7. Let $J_\Pi'$ be the set $a\in J_\Pi$ such that
$W_a^*\subset W^u({S_1})$. Then $J_\Pi'$ is closed, 
$f_\Pi(J_\Pi')=J_\Pi'$ and $J_\Pi'\ne J_\Pi$.

\give Proof. If $a\notin J_\Pi'$, then 
$W_a^*\cap W^u({S_1})=\emptyset$ by Lemma~10.5. 
Hence $W_a\cap\{G=1\}$ is a compact
subset of the open set $W^u(J)$ so by continuity 
there is an open neighborhood
$X$ of $a$ in $J_\Pi$ such that $W_b\cap\{G=1\}\subset W^u(J)$ for all
$b\in X$. By Lemma~10.5 it follows that $X\cap J_\Pi'=\emptyset$ and
we conclude that $J_\Pi-J_\Pi'$ is open. That $f_\Pi(J_\Pi')=J_\Pi'$
follows from the fact that $f(W^u({S_1}))=W^u({S_1})$.

Finally suppose $J_\Pi'=J_\Pi$. Then 
$W^s(J_\Pi)\subset J_\Pi\cup W^u(S_1)$, so
$W^s(J_\Pi)$ does not intersect $W^u(J)$. This contradicts
Theorem~7.3, because 
$W^u(J)$ contains a neighborhood of $J=\supp(\mu)$.
\qed
We say that a stable disk $W_b$ {\it lands\/} on $J$ if
$\psi_b$ extends continuously to $S^1$ and $\psi_b(S^1)\subset J$.
This does not depend on the specific choice of parametrization
$\psi_b$.

\proclaim Lemma~10.8. There exists a dense set of $b\in J_\Pi$ 
such that $W_b$ lands on $J$.

\give Proof. Since periodic points are dense in $J_\Pi$ and
$J_\Pi-J_\Pi'$ is open and nonempty, we can find a periodic point
$b'\in J_\Pi-J_\Pi'$, say of period $n$. Furthermore, $f$ is expanding
on $J$, so there exists a neighborhood
$N$ of $J$ and $\lambda>1$ with $f^{-1}(N)\subset N$ and
$$
|Df^n(y)v|\ge\lambda^n|v|,\eqno(10.5)
$$
for all $y\in N$ and all tangent vectors $v$ (we may have
to increase $n$).
Now the annulus
$\psi_{b'}\{2^{1/\deg^n}\le|\zeta|\le2\}$ in 
$W_{b'}$ is a compact subset
of $W^u(J)$, so the inverse images under sufficiently high iterates 
of $f$ of points in this annulus will be in $N$. In particular, since
$b'$ is periodic, it follows that there exists an $R>1$ such that
$\psi_{b'}\{1<|\zeta|\le R\}\subset N$. Then, using the 
estimate~(10.5) above, we may prove
that $\psi_{b'}$ extends to a H\"older continuous map of
$\bar\Delta$, mapping $S^1$ into $J$. The proof is very similar
to the proof of Lemma~10.4 and is therefore omitted.

We conclude that $W_{b'}$ lands on $J$ and so does $W_b$ 
for all preimages $b$ of $b'$ under iterates of $f$. 
Such preimages are dense in $J_\Pi$.
\qed

Figure 6 illustrates the effect of a heteroclinic connection. 
Here $W_a^*$ is in the unstable set of $S_1$ whereas $W_b$ 
lands on $J$. The stable disks in the middle are of the form $W_{b_n}$, 
where $b_n$ are preimages of $b$ converging to $a$. Note that
the disks $W_{b_n}$ are very ``bent'' for large $n$.
If the $W_{b_n}$'s were graphs over $W_a$, then this would contradict
the maximum principle. We cannot show directly that
$W_{b_n}$ are graphs over $W_a$, but we will nevertheless prove,
using the maximum principle, that the picture is impossible.

\centerline{\psfig{figure=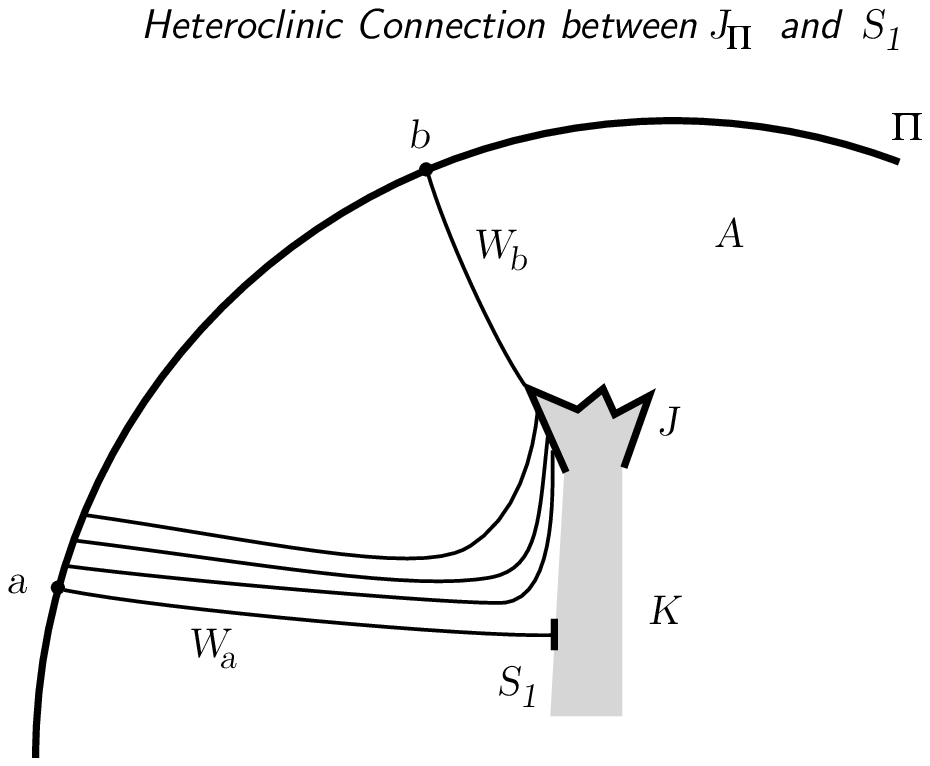,width=0.5\hsize}}
\centerline{Figure 6.}
\medskip

It follows from Lemma~10.5 that for each $a\in J_\Pi'$ there exists
a (not necessarily unique) history $\widehat{p_a}$ in $S_1$ such that 
$W_a^*\subset W^u(\widehat{p_a})$. In general, an unstable manifold
$W^u(\hat{q})$ of a history $\hat{q}$ in $S_1$ is a complicated
object, but, as we will see, 
the information that $W_a^*\subset W^u(\widehat{p_a})$ implies
that $W^u(\widehat{p_a})$ is in fact an algebraic subvariety of $\cx 2$.
Recall that the image of a holomorphic map of a compact Riemann 
surface into $\px 2$ is an algebraic variety. The authors
thank Jeffrey Diller for useful conversations on the proof of the
following result.

\proclaim Lemma~10.9. If $J_\Pi'\ne\emptyset$, then there exists an 
$a\in J_\Pi'$ such that $W^u(\widehat{p_a})$ is an 
algebraic subvariety of $\cx 2$.

\give Proof. Take any point $a\in J_\Pi'$, a complete orbit
$(a_i)_{i\in\Z}$ with $a_0=a$ and a complete
orbit $(p_i)_{i\in\Z}$ be a complete
orbit in $S_1$ such that $W_{a_i}^*\subset W^u((p_{i+j})_{j\le 0})$
for all $i$. We write $\widehat{p_i}$ for the history $(p_{i+j})_{j\le 0}$.
If $\delta>0$ is small enough, then the local
unstable manifolds $W^u_{\rm loc}(\widehat{p_i})$ are complex disks for 
all $i$ and there exist biholomorphisms 
$\eta_i:\D_{\delta_i}\to W^u_{\rm loc}(\widehat{p_i})$ 
with $|D\eta_i(0)|=1$ and complex numbers $\lambda_i\ne 0$ 
such that 
$$
\eta_i(\lambda_{i-1}\zeta)=f(\eta_{i-1}(\zeta)),\eqno(10.6)
$$
for all $i$ and all $|\zeta|<\delta_{i-1}$.
Since $f$ is hyperbolic on $S_1$, the numbers $\delta_i$ are 
uniformly bounded from below and
$\lambda_{i-n}\cdots\lambda_{i-1}\to\infty$ as $n\to\infty$
for all $i$, so~(10.6) allows us to extend $\eta_i$ to
maps of $\C$ into $W^u(\widehat{p_i})$ by defining
$$
\eta_i(\lambda_{i-n}\cdots\lambda_{i-1}\zeta)=f^{n}(\eta_{i-n}(\zeta)),
$$
for $n\ge 0$. 

The maps $\eta_i$ are surjective by the definition of
$W^u(\widehat{p_i})$ but they need not be injective. However, the global 
unstable manifolds $W^u(\widehat{p_i})$ have a natural structure 
as abstract Riemann surfaces given by the maps $\eta_i$.
More precisely, for each $i$ we define a Riemann surface $X_i$ as
the quotient $\C/\sim$, where $z\sim w$
if there are open sets $U\ni z$ and $V\ni w$
such that $\eta_i(U)=\eta_i(V)$. Then the map $\eta_i$ factors
as $\eta_i=\eta_i'\circ\pi_i$, where 
$\pi_i:\C\to X_i$ is surjective, $\eta_i':X_i\to\cx 2$ is locally
injective and the set of points $(z,w)\in X_i\times X_i$ 
with $z\ne w$ and $\eta_i'(z)=\eta_i'(w)$ is discrete.
We will be sloppy and make no distinction between the unstable
manifold $W^u(\widehat{p_i})$ and the Riemann surface $X_i$.
Hence we will sometimes view $W^u(\widehat{p_i})$ as a subset 
of $\cx 2$ and sometimes as
an abstract Riemann surface. The precise meaning should be clear
from the context. 

Now the Riemann surface 
$W^u(\widehat{p_i})$ cannot be hyperbolic, because 
$\eta_i$ maps $\C$ into it so $W^u(\widehat{p_i})$ is biholomorphic
to $\C^*$, $\C$, $\px 1$ or a torus. The last two cases cannot occur, 
because then $W^u(\widehat{p_i})$ would be an algebraic subvariety 
of $\px 2$, which is impossible (because
$W^u(\widehat{p_i})\cap\Pi=\emptyset$).
Hence $W^u(\widehat{p_i})$ is 
biholomorphic to $\C^*$ or $\C$ for all $i$.

Write $W_i$ instead of $W_{a_i}$ and note that
$W^u(\widehat{p_i})$ has an open subset biholomorphic to $W_i^*$. Let
$\Sigma_i$ be the Riemann surface obtained from $W^u(\widehat{p_i})$ 
by filling in the hole at $a_i$. Then $\Sigma_i$ is biholomorphic to
$\C$ or $\px 1$ for all $i$. If $\Sigma_i$ is biholomorphic to $\px 1$
for some $i$, then $\Sigma_i$ is an algebraic subvariety of $\px 2$ 
(in fact a line) and we are done, so assume that $\Sigma_i$ is 
biholomorphic to $\C$ for all $i$.

Suppose that $(\Sigma_i-W_i)\cap W^s(J_\Pi)\ne\emptyset$ for some
$i$. Then $(\Sigma_i-W_i)\cap W_b\ne\emptyset$ for some $b\in J_\Pi$,
$b\ne a_i$.
By the dichotomy given in Lemma~10.5 we then have that 
$W_b^*\subset(\Sigma_i-W_i)$ so by filling in the hole at $b$ we see that
the closure of $\Sigma_i$ in $\px 2$ is an algebraic subvariety of
$\px 2$, which implies that $W^u(\widehat{p_i})$ is algebraic in this case too.

Let us now suppose that $\Sigma_i$ is biholomorphic to $\C$ and that
$(\Sigma_i-W_i)\cap W^s(J_\Pi)=\emptyset$ for all $i$. Pick 
biholomorphisms $\chi_i:\C\to\Sigma_i$ such that $\chi_i(0)=a_i$.
Note that $f$ induces holomorphic maps of $\Sigma_i$ onto $\Sigma_{i+1}$.
Hence we may define entire maps
$h_i$ by $\chi_i\circ h_i=f\circ\chi_{i-1}$ for all $i$.
The restriction of $f$ to $W_{i-1}$ is a branched covering of
$W_i$ of degree
$\deg$, branched only at $a_i$. 
This implies that 
$h_i(\zeta)=\zeta^\deg\exp(u_i(\zeta))$ where $u_i$ is entire.
Moreover, the condition 
$(\Sigma_i-W_i)\cap W^s(J_\Pi)=\emptyset$ implies that the inverse
image of $W_i$ in $\Sigma_{i-1}$ is exactly $W_{i-1}$. Therefore 
$\limsup|h_i(\zeta)|>0$ as $|\zeta|\to\infty$ and
this is only possible if $u_i$ is constant. Hence we may write
$h_i(\zeta)=c_i\zeta^\deg$ for some constants $c_i\ne 0$.

Define $g_i=G\circ\chi_i$.
Then, for each $i$, $g_i\ge 0$ is continuous and 
subharmonic on $\C^*$.
The equation $G\circ f=\deg\,G$ translates into
$g_i\circ h_i=\deg\,g_{i-1}$, i.e.\ $g_i(c_i\zeta^\deg)=\deg\,g_{i-1}(\zeta)$.
Iterating this we see that $g_i(\zeta)$ depends only on $|\zeta|$.
Thus, by the maximum principle, for each $i$ there exists an
$R_i$ such that either $g_i=0$ on $|\zeta|>R_i$ 
or $g_i>0$ for $|\zeta|>R_i$.

If $g_i=0$ for $|\zeta|> R_i$, then $\chi_i$ maps
$|\zeta|>R_i$ into the bounded set $K$ and therefore extends
to a holomorphic map of $\px 1$ into $\px 2$.
Hence $W^u(\widehat{p_i})$ is algebraic.

If $g_i>0$ for $|\zeta|>R_i$, then $\chi_i$ maps $|\zeta|>R_i$
into $\cx 2-K$, and by our previous assumption, the image does
not intersect $W^s(J_\Pi)=\supp(T\contract A)$, so
$g_i$ is harmonic on $|\zeta|>R_i$. Hence there exist constants
$A_i>0$ and $B_i$ such that $g_i(\zeta)=A_i\log|\zeta|+B_i$ for
$|\zeta|>R_i$. Since $G(x)=\log|x|+O(1)$ as $x\to\Pi$, this implies that
$|\chi_i(\zeta)|\le C|\zeta|^{A_i}$ as $\zeta\to\infty$, so
again $\chi_i$ extends to a holomorphic map of $\px 1$ into $\px 2$.
Hence $W^u(\widehat{p_i})$ is algebraic, which completes the proof
of Lemma~10.9.
\qed

We are now in position to prove Lemma~10.3.

\give Proof of Lemma~10.3. Suppose that 
$W^s(J_\Pi)\cap W^u(S_1)\ne\emptyset$.
Then $J_\Pi'\ne\emptyset$ so Lemma~10.9 shows that there exist
$a\in J_\Pi$, a history $\widehat{p}$ in $S_1$ and an irreducible 
polynomial $P(z,w)$ such that $W_a^*\subset W^u(\widehat{p})=\{P=0\}$.
Clearly $W^u(\widehat{p})\cap J=\emptyset$ so there exists an
$\epsilon>0$ such that $|P|\ge 2\epsilon$ on $J$.

By Lemma~10.8 there is a dense set of $b$'s such that $W_b$ lands on
$J$. If we choose $b$ close enough to $a$, then by continuity $W_b$
will intersect the open set $|P|<\epsilon$. Let $U$
be a component of
$\{\zeta\in\Delta^*: |P(\psi_b(\zeta))|<\epsilon\}$. Then
$U$ is relatively compact in $\Delta^*$. Further, $P$ is a nonzero
holomorphic function on $U$, so $-\log|P|$ is harmonic on $U$.
But $|P|<\epsilon$ on $U$ and $|P|=\epsilon$ on $\partial U$, 
contradicting the maximum principle for $-\log|P|$ on $U$.
This completes the proof of Lemma~10.3.
\qed

\proclaim Corollary~10.10. If $f$ satisfies condition (\ddag) and
$J_\Pi$ is connected, then $J$ is connected.
If $J_\Pi$ is also locally connected, then so is $J$.
\qed

\give Proof. If $J_\Pi$ is connected (and locally connected) then
$\cE$ is connected (and locally connected).
\qed

\give Remark. In order to prove that $e$ maps $\cE$ continuously
onto $J$, we could do
without the assumption that $f_\Pi$ is expanding on $J_\Pi$. It would
be sufficient to assume that $\mu_c(A)=0$, and that (\ddag2)--(\ddag4)
hold. Indeed, by Theorem~8.8
we still have a disk lamination $\{W_a: a\in J_\Pi\}$ and $e_r$
is defined and continuous on $\cE$ for all $r>0$. The only place
where the uniform expansion of $f_\Pi$ on $J_\Pi$ is used, is to
get H\"older continuity of $e$. However, we always get
continuity of $e$.

\give Questions. Is $J$ a {\it finite\/} quotient of $\cE$, i.e.\
is $\#e^{-1}(x)$ uniformly bounded on $J$? What are the
possible identifications on $\cE$ introduced by $e$?

%
%
\section Appendix A. The homogeneous model

The model for our study of regular polynomial automorphisms is the
case when $f=f_h$ is a homogeneous mapping of $\cx\dim$. Here we show
that $f_h$ is essentially a skew product over $f_\Pi$. If $g$ is a
homogeneous polynomial of degree $N$, we let $V=\cx\dim\cap\{g=0\}$ and
$V_\Pi=\Pi\cap\{g=0\}$. We let $\Sigma_N$ denote the multiplicative
group of the $N$th roots of unity in $\C$. We let
$\check\C_*=\C_*/\Sigma_N$ denote the quotient. With this
complex structure, the mapping $\C_*\to\C_*$ given by
$\lambda\mapsto\lambda^N$ descends to an isomorphism
$\C_*\to\check\C_*$. Similarly, we define the (finite) quotient
$\check\cx\dim_*:=\cx\dim_*/\Sigma_N$. Thus we have a holomorphic mapping
$$
s:\Pi-V_\Pi\to\check\cx\dim_*
$$ given by $s(z)=g^{-1/N}(z) z$. It
follows that the mapping
$$
\psi:\C_*\times(\Pi-V_\Pi)\to \check\cx\dim_*-V
$$ given by
$\psi(\lambda,[z])=\lambda^{1/N}s(z)$ is biholomorphic. The
homogeneous mapping descends to a finite quotient mapping
$$
\check f_h:\check\cx\dim_*-(V\cup f_h^{-1}V)\to\check\cx\dim_*.
$$
Since $s$ is a section of the bundle $\pi:\check\cx\dim_*\to\Pi$, it follows
that $s\circ f_\Pi(z)$ and $\check f_h\circ s[z]$ define the same line
in $\check\cx\dim_*$. Thus $\chi(z):=f_h(s(z))/s(f_\Pi(z))$ is an $N$
valued holomorphic function on $\Pi-V_\Pi$. A short calculation shows
that
$$
\psi^{-1}\circ \check f_h\circ\psi(\lambda, z) =
(\chi^N(z)\lambda^\deg,f_\Pi(z))\eqno(A.1)
$$
where $\chi^N$ is a (single-valued) analytic function on $\Pi-V_\Pi$.

We note that for $a\in\Pi$, $W^s_{\rm loc}(a,f_h)$ is contained in the
line $L_a$. Let $C(J_\Pi)$ denote the union of $J_\Pi$ and
the complex homogeneous cone in $\cx\dim_*$ over $J_\Pi$. Our canonical
model in~\S4 is given by the restriction of $f_h$ to $C(J_\Pi)$. Let
$\check C(J_\Pi)$ denote the quotient of $C(J_\Pi)$ in
$\check\cx\dim_*\cup\Pi$.

\proclaim Proposition~A.1. Suppose that $\cC_\Pi\cap J_\Pi=\emptyset$.
Then there is a continuous function $\eta$ on $J_\Pi$ with unit
modulus such that the restriction of $\check f_h$ to $\check C(J_\Pi)$
is conjugate to the self-mapping of $\C_*\times J_\Pi$ given by
$$
(\lambda,z)\mapsto(\eta(z)\lambda^\deg,f_\Pi(z)).
$$
In addition, if $\eta$ has a continuous logarithm on $J_\Pi$, then
$\check f_h|_{\check C(J_\Pi)}$ is conjugate to
$(\lambda,[z])\mapsto(\lambda^\deg,f_\Pi[z])$.

\give Proof. If we let $g$ denote the Jacobian determinant of $f_h$
on $\cx\dim$, then $g$ is a homogeneous polynomial of degree
$N=\dim(\deg-1)$.
In particular, $V=\{g=0\}$ is the critical locus of $f_h$, and $V_\Pi$
is disjoint from $J_\Pi$. Thus we may represent $\check f_h$ as in
$(A.1)$.

Let us define $\phi(\lambda,z)=(\alpha(z)\lambda,z)$. Then
$\phi^{-1}\circ\psi^{-1}\circ \check f_h\circ\psi\circ\phi$ is given
by
$$
(\lambda,z)\mapsto(\alpha^\deg(z)\alpha^{-1}(f_\Pi
z)\chi^N(z)\lambda^\deg,f_\Pi z).
$$
Thus we wish to find $\alpha:J_\Pi\to\R$ such that
$$
\alpha^\deg(z)\alpha^{-1}(f_\Pi z)\,|\chi(z)|=1.\eqno(A.2)
$$
Taking logarithms, we have
$$
\log\alpha(z)=-{1\over\deg}\log|\chi^N(z)|+{1\over\deg}\log|\alpha(f_\Pi
z)|.$$ Applying $f^j_\Pi$, dividing by $\deg^j$, and summing
over $j$, we have that
$$
\log\alpha(z)= -\sum_{j=0}^\infty {1\over
\deg^{j+1}}\log|\chi^N(f^j_\Pi z)|\eqno(A.3)
$$
is continuous on $J_\Pi$,
and $\alpha$ solves (A.2). Setting $\eta=\chi^N/|\chi^N|$, we have
the desired form of $\check f_h$.

Finally, if $\eta$ has a continuous logarithm, we may solve (A.2)
without taking absolute value. \qed

\give Examples. We consider two mappings:
$f_1(z_1,z_2)=(z_1^2,z_2^2)$ and $f_2(z_1,z_2)=(z_2^2,z_1^2)$. We
consider the function $g(z)=z_1$, the coordinate $\zeta=z_2/z_1$ on
$\Pi-\{[0:1]\}$, and the mapping
$\psi:\C_*\times(\Pi-\{[0:1]\}\to\cx2_*$ given by
$\psi(\lambda,\zeta)=\lambda(1,\zeta)=(\lambda,\lambda\zeta)$. In
both cases, $J_\Pi=\{|\zeta|=1\}$. The associated mappings on $J_\Pi$
are $f_{1,\Pi}(\zeta)=\zeta^2$ and $f_{2,\Pi}(\zeta)=\zeta^{-2}$. The
normal forms given in the Proposition are
$$
\psi^{-1}f_1\psi(\lambda,\zeta)=(\lambda^2,\zeta^2),\qquad
\psi^{-1}f_2\psi(\lambda,\zeta)=(\zeta^2\lambda^2,\zeta^{-2}).
$$
There is no continuous function $\alpha:J_\Pi\to\C_*$ solving $(A.3)$ for
$f_2$. For in this case we would have
$\alpha^2(\zeta)\alpha^{-1}(\zeta^{-2})\zeta^2=1$. But if $A$ is the
(integer) winding number of $\alpha\{|\zeta|=1\}$ about 0 in
$\C_*$, then the winding numbers of $\alpha^2(\zeta)$ and
$\alpha^{-1}(\zeta^{-2})$ are each $2A$, so we must have $2A+2A+2=0$,
which is a contradiction.

In dimension $\dim=2$ we can often assume that $N=1$.

\proclaim Proposition~A.2. Suppose that $\dim=2$ and that
$J_\Pi\ne\Pi$. Then there is a continuous function $\eta$
on $J_\Pi$ with unit modulus such that the restriction of
$f_h$ to $C(J_\Pi)$
is conjugate to the self-mapping of $\C_*\times J_\Pi$ given by
$$
(\lambda,z)\mapsto(\eta(z)\lambda^\deg,f_\Pi(z)).
$$
In addition, if $\eta$ has a continuous logarithm on $J_\Pi$, then
$f_h|_{C(J_\Pi)}$ is conjugate to
$(\lambda,[z])\mapsto(\lambda^\deg,f_\Pi[z])$.

\give Proof. The proof is the same as for Proposition~A.1, except that
we use $g(z,w)=a_2z-a_1w$, where $a=[a_1:a_2]\in\Pi-J_\Pi$.
\qed

%
%
\section Appendix B. Hyperbolicity for endomorphisms.

In this appendix we present some basic results on hyperbolicity for
smooth endomorphisms. More details can be found in [J2]. 
Our main references are [Ru] and [PS]; see also [FS4]. No proofs
are given in this appendix; they can be found in the above references.

Let $f$ be a $C^\infty$ endomorphism of a finite-dimensional
Riemannian manifold $M$.
Let $L$ be a compact subset of $M$ with $f(L)=L$ and define
$$
\hat L=\{(x_i)_{i\le 0}: x_i\in L, fx_i=x_{i+1}\}.
$$
Then $\hat L$ is a closed subset of $L^{\N}$, hence compact.
We will use the notation $\hat x$
for a point $(x_i)_{i\le 0}$ in $\hat L$.
The restriction
$f|_L$ lifts to a homeomorphism $\hat{f}$ of $\hat{L}$ given by
$\hat{f}((x_i))=(x_{i+1})$. There is a natural projection $\pi$ from
$\hat{L}$ to $L$ sending $\hat x$ to $x_0$ and the pullback under $\pi$
of the restriction to $L$ of the tangent bundle of $M$ is a
bundle on $\hat{L}$ which we call the tangent bundle $T_{\hat{L}}$.
Explicitly, a point in $T_{\hat{L}}$ is of the form
$(\hat x,v)$ where $\hat x\in\hat{L}$ and $v$ is a tangent vector in
$T_{x_0}M$. The derivative $Df$ lifts to a
map $D\hat{f}$ of $T_{\hat{L}}$ in a natural way.

Now $f$ is {\it hyperbolic\/}
on $L$ if there exists a continuous
splitting $T_{\hat{L}}=E^u\oplus E^s$ which is invariant under $D\hat{f}$
and such that $D\hat{f}$ is expanding on $E^u$ and contracting on
$E^s$. More precisely, $D\hat{f}(E^{u/s})\subset E^{u/s}$
and there are constants $c>0$ and $\lambda>1$ such that for all $n\ge 1$
$$
\eqalign{
|D\hat{f}^nv|&\ge c\lambda^n|v|\quad v\in E^u\cr
|D\hat{f}^nv|&\le c^{-1}\lambda^{-n}|v|\quad v\in E^s.}
$$

\give Remark. Such a map is called prehyperbolic in [Ru].

\give Remark. It is possible to make a smooth change of metric
in a neighborhood of $L$ and obtain $c=1$ in the equation above.

Note that whereas the fiber of the unstable bundle $E^u$ at a point
$\hat x\in\hat{L}$ depends on the whole history $\hat x$ of $x_0$,
the fiber of $E^s$ at $\hat x$ depends only on the point $x_0$.
Hence the dimension of the fiber of $E^u$ at a point 
$\hat x$ depends only on $x_0$, so the dimensions of the fibers of 
the bundles $E^u$ and $E^s$ are locally constant.

As a special case of the above we say that $f$ is expanding on $L$
if the bundle $E^s$ is trivial. This means
that there exist constants $c>0$ and $\lambda>1$ such that
$|D\hat{f}^n(x)v|\ge c\lambda^n|v|$ for all $x\in L$, $v\in T_xM$
and all $n\ge 1$.

A basic result in hyperbolic dynamics is the stable manifold
theorem. For each point $p$ in $L$ and each history $\hat q$
in $\hat{L}$, we define local stable and unstable manifolds by
$$
\eqalign{
W^s_{\rm loc}(p)
&=\{y\in M: d(f^iy,f^ip)<\delta\ \forall i\ge 0\}\cr
W^u_{\rm loc}(\hat q)
&=\{y\in M: \exists \hat y, \pi(\hat y)=y,
d(y_i,q_i)<\delta\ \forall i\le 0\},}
$$
for small $\delta>0$.

The following theorem asserts that the local (un)stable manifolds are
indeed nice objects.  [Ru,~\S15] contains an outline of a proof,
whereas~[PS, Theorem~5.2] proves a more general theorem).

\proclaim Theorem~B.1 (Stable Manifold Theorem). If $\delta$ 
is small enough, then
\item{(i)} For all $p\in L$ and
  all $\hat q\in\hat L$, $W^s_{\rm loc}(p)$
  and $W^u_{\rm loc}(\hat q)$ are embedded $C^\infty$ balls of
  $M$ tangent to $E^s(p)$ and $E^u(\hat q)$ at
  $p$ and $q_0$, respectively.
\item{(ii)} $W^s_{\rm loc}(p)$ and $W^u_{\rm loc}(\hat q)$ depend
  continuously on $p$ and $\hat q$, respectively.
\item{(iii)} If $x\in W^s_{\rm loc}(p)$, then $d(f^nx,f^np)\to 0$ 
  exponentially fast as $n\to\infty$. Similarly, every point $x$ in 
  $W^u_{\rm loc}(\hat q)$ has a unique history $\hat x$ such that
  $x_j\in W^u_{\rm loc}({\hat f}^j(\hat q))$ for all 
  $j\le 0$ and $d(x_j,q_j)\to 0$ exponentially fast as $j\to-\infty$.

If $\delta$ is small enough, then by continuity $W^s_{\rm loc}(p)$ and
$W^u_{\rm loc}(\hat q)$ are almost flat, i.e.\ $C^1$ close to the tangents
at $p$ and $q_0$, respectively for all $p\in L$ and all $q\in\hat L$.
Therefore $W^s_{\rm loc}(p)$ and $W^u_{\rm loc}(q)$ intersect in at most 
one point. 

\proclaim Definition B.2. We say that $L$ has 
local product structure 
if $\delta$ can be chosen so that
$W^s_{\rm loc}(p)\cap W^s_{\rm loc}(\hat q)\subset L$ for all
$p$ and $\hat q$.

If $L$ has local product structure, $p\in L$, $\hat q\in\hat L$
and if $p$,$q_0$ are sufficiently
close, then $W^s_{\rm loc}(p)$ and $W^u_{\rm loc}(\hat q)$ intersect
in exactly one point $x\in L$ and $x$ has a history $\hat x$
such that $x_j\in W^u_{\rm loc}({\hat f}^j(\hat q))$ for all $j\le 0$.
It is not a priori clear that $\hat x\in\hat L$, i.e.\ that
$x_j\in L$ for all $j\le 0$. We therefore make another definition.

\proclaim Definition B.3. We say that $\hat L$ has 
local product structure 
if $\delta$ can be chosen so that if the intersection
$W^s_{\rm loc}(p)\cap W^s_{\rm loc}(\hat q)$
is nonempty, then it consists of a unique point $x\in L$ and
the unique history $\hat x$ of $x$ with
$x_j\in W^u_{\rm loc}({\hat f}^j\hat q)$ 
for all $j\le 0$ is contained in $\hat L$.

\proclaim Definition B.4. Let $\eta>0$. An $\eta$-pseudoorbit in $M$
is a sequence $(x_i)_{[t_1,t_2]}$, where 
$-\infty\le t_1<t_2\le\infty$, such that 
$d(fx_i,x_{i+1})<\delta$ for $t_1\le i<t_2$.
An $\eta$-pseudoorbit $(x_i)_{[t_1,t_2]}$ is $\epsilon$-shadowed
by an orbit $(y_i)_{[t_1,t_2]}$ if $d(y_i,x_i)<\epsilon$
for all $i\in[t_1,t_2]$.

For proofs of the remaining results in this appendix see [J2].

\proclaim Theorem~B.5 (Shadowing Lemma). Suppose that $L$ is
hyperbolic and that $\hat L$ has
local product structure. Then for each $\epsilon>0$ there exists
an $\eta>0$ such that every $\eta$-pseudoorbit in $L$ can be
$\epsilon$-shadowed by an orbit in $L$.

Using shadowing we control the orbits of $f$ staying near $L$
in positive or negative time.

\proclaim Proposition~B.6 (Fundamental Neighborhood).
Let $L$ be a hyperbolic set 
for a map $f$. Assume that $\hat L$ has local product structure. 
Then $L$ has a neighborhood $U$ in $M$ such that
\item{(i)} If $x\in U$ and $f^jx\in U$ for all $j\ge 0$, then
  $x\in W^s_{\rm loc}(p)$ for some $p\in L$.
\item{(ii)} If $x\in U$ and $x$ has a history $\hat x$ with
  $x_i\in U$ for all $i\le 0$, then $x\in W^u_{\rm loc}(\hat q)$
  for some $\hat q\in\hat L$. More precisely
  $d(x_i,q_i)<\delta$ for all $i\le 0$.
\item{(iii)} If $(x_i)_{i\in\Z}$ is a complete orbit in $U$ then
  $x_i\in L$ for all $i$.

Next we consider Axiom A endomorphisms. A point $x\in M$ 
is {\sl wandering} if it has a neighborhood $V$ such
that $f^n(V)\cap V=\emptyset$ for all $n\ge 1$; otherwise it is called
{\sl non-wandering}. The {\sl non-wandering set} $\Omega$ of $f$ 
is the set of all non-wandering points; it is a closed set.

\proclaim Definition B.7. A map $f$ is said to be Axiom A if
its non-wandering set satisfies
\item{(i)} $\Omega$ is compact.
\item{(ii)} Periodic points are dense in $\Omega$.
\item{(iii)} $f$ is hyperbolic on $\Omega$.

\give Remark. If $\Omega$ satisfies (ii), then
$f(\Omega)=\Omega$, so (iii) makes sense. Also, if $f$ is Axiom A, then
periodic points (under $\hat f$) are dense in $\hat\Omega$.

The next proposition shows that the preceding results apply to
open Axiom A endomorphisms.

\proclaim Proposition~B.8. If $f$ is Axiom A and open,
then $\hat\Omega$ has local product structure.

\proclaim Theorem~B.9 (Spectral decomposition).
If $f$ is Axiom A, then $\Omega$ can be written 
in a unique way as
a disjoint union $\Omega=\cup_{i=1}^l{\Omega}_i$, 
where each $\Omega_i$
is compact, satisfies $f(\Omega_i)=\Omega_i$ 
and $f$ is transitive
on $\Omega_i$. The sets $\Omega_i$ 
are called the basic sets of $f$. Morover,
each $\Omega_i$ can be further decomposed 
into a finite disjoint union 
$\Omega_i=\cup_{1\le j\le n_i}\Omega_{i,j}$, 
where $\Omega_{i,j}$ is compact,
$f(\Omega_{i,j})=\Omega_{i,j+1}$
($\Omega_{i,n_i+1}=\Omega_{i,1}$) and 
$f^{n_i}$ is topologically mixing on each $\Omega_{i,j}$.

Our final result in this appendix describes forward and backward
orbits for an Axiom A endomorphism.

\proclaim Corollary~B.10. Assume that $f$ is Axiom A and $M$ is compact.
\item{(i)} If $x\in M$, then there is a unique basic set $\Omega_i$
  such that $f^jx\to\Omega_i$ as $j\to\infty$. Moreover, there is
  a (not necessarily unique) $p\in\Omega_i$ such that
  $d(f^jx,f^jp)\to 0$ as $j\to\infty$.
\item{(ii)} If $\hat x\in\hat M$, then there is a unique
  basic set $\Omega_i$ such that $x_j\to\Omega_i$ as $j\to-\infty$.
  Moreover, there is a (not necessarily unique) 
  $\hat q\in\widehat{\Omega_i}$ such that $d(x_j,q_j)\to 0$ 
  as $j\to-\infty$.
%
%

\centerline{\bf References}

\item{[A]} H.\ Alexander, Gromov's method and Bennequin's problem.
  Invent. Math. 125, 135--148 (1996).
\item{[AT]} H.\ Alexander, B.\ A.\ Taylor, Comparison of two capacities
  in $\cx n$, Math.\ Z.\ 186, 407--417 (1984).
\item{[Ba]} D.\ Barrett, Holomorphic extension from boundaries with 
  concentrated Levi form, Indiana U. Math. J., 44 (1995), 1075--1087.
\item{[Be]} B.\ Berndtsson, Personal communication.
\item{[BLS]} E.\ Bedford, M.\ Lyubich, J.\ Smillie, Polynomial 
  diffeomorphims of ${\bf C}^2$ IV: The measure of maximal entropy and 
  laminar currents. Invent. Math. 112, 77--125 (1993).
\item{[BS1]} E.\ Bedford, J.\ Smillie, Polynomial 
  diffeomorphims of ${\bf C}^2$: currents, equilibrium measure and
  hyperbolicity. Invent. Math. 103, 69--99 (1991). 
\item{[BS2]} E.\ Bedford, J.\ Smillie, Polynomial 
  diffeomorphims of ${\bf C}^2$ V: Critical points and Lyapunov exponents.
  To appear in J. Geom. Anal.
\item{[BT1]} E.\ Bedford, B.\ A.\ Taylor, A new capacity for
  plurisubharmonic functions. Acta Math. 149, 1--39 (1982).
\item{[BT2]} E.\ Bedford, B.\ A.\ Taylor, Plurisubharmonic functions with
  logarithmic singularities. Ann. Inst. Fourier (Grenoble) 38, 133-171 (1988).
\item{[Bri]} J.-Y.\ Briend. Exposants de Liapounoff des endomorphismes
  holomorphes de $\C\P^k$. Thesis, Universit\'e Paul Sabatier--Toulouse III,
  1997.
\item{[Bro]} H.\ Brolin, Invariant sets under iteration of rational
  functions. Ark. Mat. 6, 103--144 (1965).
\item{[C]} A.\ Candel, Uniformization of surface laminations.
  Ann.\ Scient.\ \'Ec.\ Norm.\ Sup., 4 s\'erie t. 26 (1993), 489--516.
\item{[CG]} L.\ Carleson, T.\ W.\ Gamelin, Complex dynamics.
  Springer-Verlag (1993).
\item{[FS1]} J.\ E.\ Forn{\ae}ss, N.\ Sibony, Complex dynamics in higher 
  dimension I. Ast\'e\-risque 222, 201--231 (1994).
\item{[FS2]} J.\ E.\ Forn{\ae}ss, N.\ Sibony, Complex dynamics in higher 
  dimension II. In Annals of Mathematics Studies 137, 135--182. Princeton
  University Press, 1995.
\item{[FS3]} J.\ E.\ Forn{\ae}ss, N.\ Sibony, Complex dynamics in 
  higher dimension, In P.\ M.\ Gauthier, G.\ Sabidussi, editors, Complex
  potential theory, pages 131--186. Kluwer Academic Publishers, 1994.
\item{[FS4]} J.\ E.\ Forn{\ae}ss, N.\ Sibony, Hyperbolic maps on $\px 2$.
  Math.\ Ann., to appear.
\item{[H]} S.\ Heinemann, Julia sets for holomorphic endomorphisms 
  of $\cx n$. Ergodic Theory Dynam. Systems 16, 1275--1296 (1996).
\item{[HP]} J.\ H.\ Hubbard, P.\ Papadopol, Superattractive fixed points 
  in $\cx n$. Indiana Univ. Math. J. 43, 321--365 (1994).
\item{[J1]} M.\ Jonsson, Sums of Lyapunov exponents for some polynomial 
  maps of $\cx 2$. To appear in Ergodic Theory Dynam. Systems.
\item{[J2]} M.\ Jonsson, Thesis, Royal Institute of Technology, 1997.
\item{[La]} S.\ Lang, Introduction to complex hyperbolic spaces.
  Springer-Verlag, 1987.
\item{[Le]} G.\ Levin, Disconnected Julia set and rotation sets, Ann.\
scient.\ \'Ec.\ Norm.\ Sup., t.~29, (1996) 1--22.
\item{[M]} J.\ Milnor, Dynamics in one complex variable: introductory 
  lectures. Pre\-print SUNY Stony Brook (1990).
\item{[Pe]} G.\ Peng. On the dynamics of non-degenerate polynomial
  endomorphisms of $\cx2$. PhD thesis, CUNY, 1997.
\item{[Pr]} F.\ Przytycki, Hausdorff dimension of harmonic measure on 
  the boundary of an attractive basin for a holomorphic map. Invent. Math.
  80, 161--179 (1985).
\item{[PS]} C.\ Pugh, M.\ Shub, Ergodic attractors. Trans. Amer. Math. Soc.
  312, 1--54 (1989).
\item{[Ro]} J.-P.\ Rosay, A remark on the paper by H. Alexander on
  Bennequin's problem. Invent. Math. 126, 625--627 (1996).
\item{[Ru]} D.\ Ruelle, Elements of differentiable dynamics and bifurcation
  theory. Academic Press, 1989.
\item{[T]} P.\ Tortrat, Aspects potentialistes de l'it\'eration 
  des polyn\^omes. S\'eminaire de Th\'eorie du Potentiel.
  Lecture Notes in Math 1235, 195--209 (1987).
\item{[U]} T.\ Ueda, Fatou sets in complex dynamics on projective spaces.
  J. Math. Soc. Japan 46, 545--555 (1994).
\item{[Y]} L.-S.\ Young, Ergodic theory of differentiable dynamical systems.
  In Branner{,} B.{,} Hjorth{,} P., editors, Real and complex 
  dynamical systems, pages 293--336, Kluwer Academic Publishers (1995).

\bigskip
\bigskip
\noindent Eric Bedford\par
\noindent Indiana University\par
\noindent Bloomington, IN 47405\par
\noindent bedford@indiana.edu
\bigskip
\noindent Mattias Jonsson\par
\noindent Royal Institute of Technology\par
\noindent 100 44 Stockholm, Sweden\par
\noindent mjo@math.kth.se

\bye

%% file: imsmark.tex
\def\SBIMSMark#1#2#3{
 \font\SBF=cmss10 at 10 true pt
 \font\SBI=cmssi10 at 10 true pt
 \setbox0=\hbox{\SBF Stony Brook IMS Preprint \##1}
 \setbox2=\hbox to \wd0{\hfil \SBI #2}
 \setbox4=\hbox to \wd0{\hfil \SBI #3}
 \setbox6=\hbox to \wd0{\hss
             \vbox{\hsize=\wd0 \parskip=0pt \baselineskip=10 true pt
                   \copy0 \break%
                   \copy2 \break%
                   \copy4 \break}}
 \dimen0=\ht6   \advance\dimen0 by \vsize \advance\dimen0 by 8 true pt
 \dimen2=\hsize \advance\dimen2 by .25 true in
 \ht6=0pt \dp6=0pt
 \setbox8=\vbox to \dimen0{\vfill \hbox to \dimen2{\hss \copy6}}
 \ht8=0pt \dp8=0pt \wd8=0pt
 \copy8
}